\documentclass[12pt]{article}

\usepackage{amsfonts,amsmath, epsfig,latexsym, subfigure, url}

\newtheorem{proposition}{Proposition}
\newtheorem{theorem}{Theorem}
\newtheorem{conjecture}{Conjecture}

\newtheorem{lemma}{Lemma}

\newtheorem{remark}{Remark}
\newcommand{\proof}{\noindent{\bf Proof.}\ \ }

\begin{document}

\title{$4$-valent plane graphs with $2$-, $3$- and $4$-gonal faces}

\author{Michel DEZA\\
CNRS/ENS, Paris and Institute of Statistical Mathematics, Tokyo,\\
\ Mathieu DUTOUR \\
ENS, Paris and Hebrew University, Jerusalem,
\footnote{Research of the second author was financed by EC's IHRP Programme, within the Research Training Network ``Algebraic Combinatorics in Europe,'' grant HPRN-CT-2001-00272.}\\
\ Mikhail SHTOGRIN 
\thanks{Third author acknowledges financial support of the Russian Foundation of Fundamental Research (grant 02-01-00803) and the Russian Foundation for Scientific Schools (grant 00-15-96011)}\\
Steklov Mathematical Institute, Moscow, Russia.
}
\date{\today}

\maketitle

\begin{abstract}
Call {\em $i$-hedrite} any $4$-valent $n$-vertex plane graph, whose 
faces are $2$-, $3$- and $4$-gons only and $p_2+p_3=i$. The edges of an $i$-hedrite, as of 
any Eulerian plane graph, are partitioned
by its {\em central circuits}, i.e. those, which are obtained by starting with an
edge and continuing at each vertex by the edge opposite the entering one. 
So, any $i$-hedrite is a projection of an alternating link, whose components
correspond to its central circuits.

Call an $i$-hedrite {\em irreducible}, if it has no 
{\em rail-road}, i.e. a 
circuit of $4$-gonal faces, in which every $4$-gon is adjacent to two of its 
neighbors on opposite edges.

We present the list of all $i$-hedrites with at most $15$ vertices. Examples of other results:\\[-6mm]
\begin{itemize}
\item[(i)] All $i$-hedrites, which are not $3$-connected, are identified.\\[-6mm]
\item[(ii)] Any irreducible $i$-hedrite has at most $i-2$ central circuits.\\[-6mm]
\item[(iii)] All $i$-hedrites without self-intersecting central circuits are listed.\\[-6mm]
\item[(iv)] All symmetry group of $i$-hedrites are listed.
\end{itemize}
\end{abstract}



\section{Introduction}

See \cite{Gr} for terms used here for plane graphs.
It is well-known that the p-vector of any $4$-valent plane graph satisfies to
$2p_2+p_3=8+ \sum_{i\geq 5} (i-4)p_i$.
Some examples of applications of plane $4$-valent graphs are {\em projections
of links}, {\em rectilinear embedding} in VLSI and {\em Gauss crossing 
problem} for plane graphs (see, for example, \cite{Liu}).

\vspace{2mm}

Call an {\em $i$-hedrite} any plane $2$-connected
$4$-valent graph, such that the number
$p_j$ of its $j$-gonal faces is zero for any $j$, different from 
$2,3$ and $4$, and such that $p_2=8-i$. So, 
an $n$-vertex $i$-hedrite has $(p_2, p_3, p_4)=(8-i, 2i-8, n+2-i)$.
Clearly, $(i;p_2,p_3)=(8;0,8)$, $(7;1,6)$, $(6;2,4)$,
$(5;3,2)$ and $(4;4,0)$ are all possibilities. 

An $8$-hedrite is called {\em octahedrite} in \cite{DSt}; in fact, this paper is a follow-up of \cite{DSt}.
In a way, this paper continues the program of Kirkman (\cite{Kirk} p.~282) of classification of projections of alternating links.

See on the Table below short presentation of $i$-hedrites and their symmetry groups. In the last column we indicate Goldberg-Coxeter operation $GC_{k,l}(G_0)$ starting from $i$-hedrite $G_0$ with smallest number of vertices (see \cite{Gold37}, \cite{Cox71} and \cite{DD03}).

For $3$-connected plane graphs without $2$-gons, the following Theorem of Mani is valid: the symmetry group of the graph can be realized as point group of a convex polyhedron having this graph as skeleton. But in the presence of $2$-gonal faces, one cannot speak of convex polyhedra but we expect that the Mani result still holds for this more general case.

\begin{center}
{\scriptsize
\begin{tabular}{||c|c|l|c|c|c||}
\hline
\hline
$i$   &\begin{tabular}{c}
$p_2$, $p_3$\\
$n$-vertex
\end{tabular}
&\begin{tabular}{l}
All possible point\\
symmetry groups $\Gamma$
\end{tabular}
&First $i$-hedrite $G_0$&\begin{tabular}{c}
First pure $i$-hedrite
\end{tabular}
&\begin{tabular}{c}
$GC_{k,l}(G_0)$\\
\end{tabular}\\\hline
$4$&\begin{tabular}{c}
$4$,$0$\\
$n\geq 2$,\\
$n$ even
\end{tabular}
&\begin{tabular}{l}
$5$: $D_2$, $D_{2d}$, $D_{2h}$,\\
$D_4$, $D_{4h}$,\\
i.e. all\\
$D_2\leq \Gamma\leq D_{4h}$
\end{tabular}
&\begin{tabular}{c}
\\[-1mm]
\epsfxsize=1.8cm
\epsffile{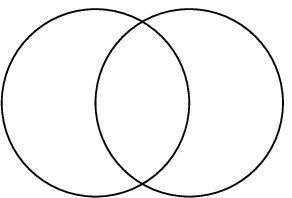}\\
{\bf Nr.2-1} \quad $D_{4h}$
\end{tabular}
&\begin{tabular}{c}
\\[-1mm]
\epsfxsize=1.8cm
\epsffile{4-hedrite2_1.eps}\\
{\bf Nr.2-1} \quad $D_{4h}$
\end{tabular}
&\begin{tabular}{c}
{\bf all}-$D_4$, $D_{4h}$\\[1mm]
$4$-hedrites,\\[1mm]
$n=2(k^2+l^2)$
\end{tabular}\\\hline
$5$&\begin{tabular}{c}
$3$,$2$\\
$n\geq 3$,\\
$n\not=4$
\end{tabular}
&\begin{tabular}{l}
$6$: $C_1$, $C_s$, $C_2$, $C_{2v}$\\
and $D_3$, $D_{3h}$
\end{tabular}
&\begin{tabular}{c}
\\[-1mm]
\epsfxsize=1.8cm
\epsffile{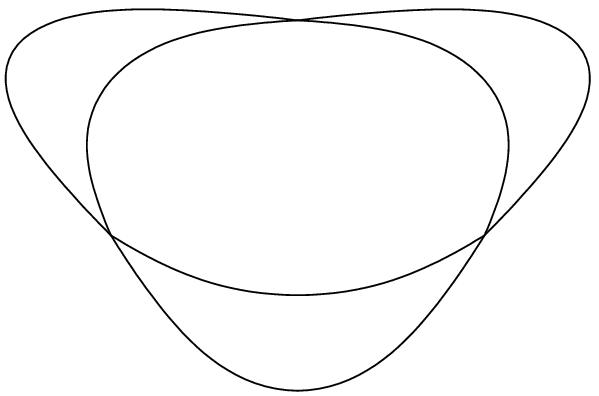}\\
{\bf Nr.3-1} \quad $D_{3h}$
\end{tabular}
&\begin{tabular}{c}
\\[-1mm]
\epsfxsize=1.8cm
\epsffile{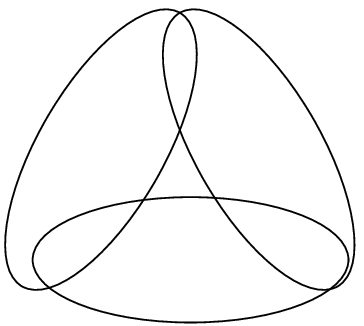}\\
{\bf Nr.6-2} \quad $D_{3h}$
\end{tabular}
&\begin{tabular}{c}
{\bf all}-$D_3$, $D_{3h}$\\[1mm]
$5$-hedrites,\\[1mm]
$n=3(k^2+l^2)$
\end{tabular}\\\hline
$6$&\begin{tabular}{c}
$2$,$4$\\
$n\geq 4$
\end{tabular}
&\begin{tabular}{l}
$9$: $C_1$, $C_s$, $C_2$, $C_{2v}$,\\
$C_i$, $C_{2h}$,\\
$D_2$, $D_{2d}$, $D_{2h}$,\\
i.e. all $\Gamma\leq D_{2h}$\\
or $\leq D_{2d}$
\end{tabular}
&\begin{tabular}{c}
\\[-1mm]
\epsfxsize=1.8cm
\epsffile{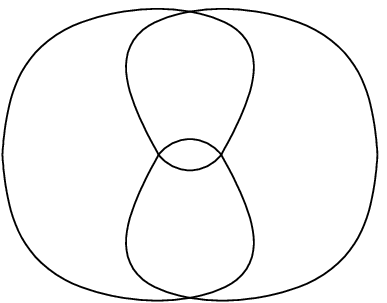}\\
{\bf Nr.4-1} \quad $D_{2d}$
\end{tabular}
&\begin{tabular}{c}
\\[-1mm]
\epsfxsize=1.8cm
\epsffile{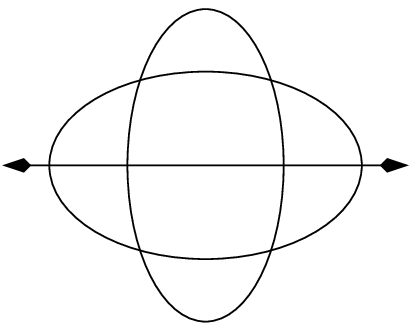}\\
{\bf Nr.8-5} \quad $D_{2h}$
\end{tabular}
&\begin{tabular}{c}
$D_2$, $D_{2d}$\\[1mm]
$6$-hedrites,\\[1mm]
$n=4(k^2+l^2)$
\end{tabular}\\\hline
$7$&\begin{tabular}{c}
$1$,$6$\\
$n\geq 7$
\end{tabular}
&\begin{tabular}{l}
$4$: $C_1$, $C_s$, $C_2$, $C_{2v}$,\\
i.e. all $\Gamma\leq C_{2v}$
\end{tabular}
&\begin{tabular}{c}
\\[-1mm]
\epsfxsize=1.8cm
\epsffile{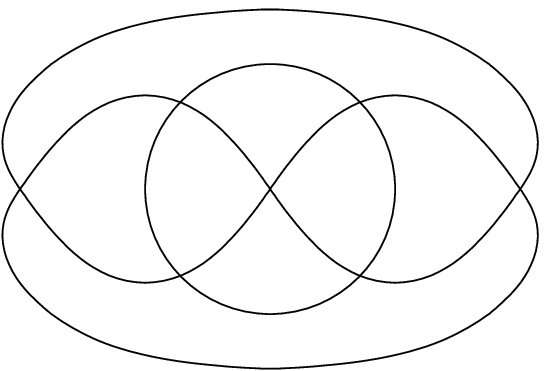}\\
{\bf Nr.7-1} \quad $C_{2v}$
\end{tabular}
&Not exist
&\begin{tabular}{c}
$C_2$, $C_{2v}$\\[1mm]
$7$-hedrites,\\[1mm]
$n=7(k^2+l^2)$
\end{tabular}\\\hline
$8$&\begin{tabular}{c}
$0$,$8$\\
$n\geq 6$,\\
$n\not= 7$
\end{tabular}
&\begin{tabular}{l}
$18$: all $13$ above\\
and $D_{3d}$, $D_{4d}$,\\
$O$, $O_{h}$, $S_4$
\end{tabular}
&\begin{tabular}{c}
\\[-1mm]
\epsfxsize=1.8cm
\epsffile{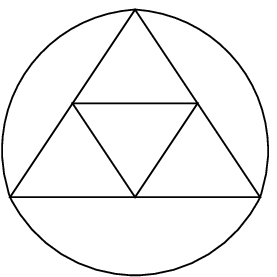}\\
{\bf Nr.6-1} \quad $O_h$
\end{tabular}
&\begin{tabular}{c}
\\[-1mm]
\epsfxsize=1.8cm
\epsffile{8-hedrite6-1.eps}\\
{\bf Nr.6-1} \quad $O_h$
\end{tabular}
&\begin{tabular}{c}
{\bf all}-$O$, $O_h$\\[1mm]
$8$-hedrites,\\[1mm]
$n=6(k^2+l^2)$
\end{tabular}\\\hline\hline
\end{tabular}
}
\end{center}

\section{Central circuits partition}

In this Section, we consider a connected plane graph $G$ with all vertices of 
even degree, i.e. an Eulerian graph. 
Call a circuit in $G$ {\it central} if it is obtained by starting with an
edge and continuing at each vertex by the edge opposite the entering one; such 
circuit is called also {\em traverse} 
(\cite{GK}), {\em straight ahead} (\cite{Ha}),  \cite{PTZ}), 
{\em straight Eulerian} (Chapter 17 of \cite{God}), 
{\em cut-through} (\cite{Je}),
{\em intersecting}, etc. Clearly, the edge-set of 
$G$ is partitioned by all its central circuits.

Such CC-partition can be considered (see, for example, \cite{Ha}) for any 
drawing on the plane of any Eulerian (in general, not planar) graph, so that 
edges are mapped into simple curves with at most one crossing point. 

Denote by 
$CC(G)=(...,a_i^{\alpha_i},...;...,b_j^{\beta_j},...)$ its {\it CC-vector}, 
where $...,a_i,...$ and  $...,b_j,...$ are increasing sequences of lengths of 
all its central circuits, simple ones and self-intersecting, respectively, 
and $\alpha_i, \beta_j$ are their respective multiplicities.
Clearly, $\sum_{i} a_i{\alpha_i}+ \sum_{j} b_j{\beta_j}=2n$, where 
$n$ is the number of vertices of $G$.

For a central circuit $C$, denote by $Int(C):=(c_0;...,c_k^{\gamma_k},...)$,
the {\em intersection vector of} $C$, where $c_0$ is
the number of self-intersections of the circuit $C$ and $...,c_k,...$ is
decreasing sequence of sizes of its intersection with other central
circuits, while the numbers $\gamma_k$ are respective multiplicities.

Two central circuits intersect in an even number number of vertices. 
The length of a central circuit is twice the number of its points of 
self-intersection plus the sum of its intersections with other circuits, 
so the length of a central circuit is even.

We will say that an $i$-hedrite is {\it pure} if any of its central circuits 
simple, i.e. has no self-intersections.
Easy to check that any pure $i$-hedrite has an even number $n$ of 
vertices. In fact, any vertex in this case belong to the intersection 
of exactly two central circuits.

Call an Eulerian graph $G$ {\it balanced}, if all its central circuits of
same length have the same intersection vector.
Any $8$-hedrite with $n \le 21$ is balanced, but there is unbalanced
$22$-vertex $8$-hedrite, which is $8$-hedrite {\bf 14-7} of Table 
\ref{tab:i-hedrite13_14} inflated along a central circuit of length $8$.
We do not find unbalanced $5$-hedrite or $7$-hedrite for $n\leq 15$. The
first unbalanced $6$-hedrite is {\bf 12-12}. Any $4$-hedrite is balanced
(Theorem \ref{Theorem-for-4-hedrite}).

For a plane graph $G$, denote by $G^*$ its plane dual and by {\em $Med(G)$} 
its {\em medial} graph. The vertices of $Med(G)$ are the edges of 
$G$, two of them being adjacent if the corresponding edges share a vertex and 
belong to the same face of the embedding of $G$ in the plane. 
So, $Med(G)$=$Med(G^*)$.

Clearly, $Med(G)$ is a $4$-valent plane graph and, for any $i$-hedrite $G$,
$Med(G)$ is an $i$-hedrite with twice the number of vertices of $G$ and
all $2$-, $3$-gonal faces being isolated. The medial of smallest
$8$-hedrite {\bf 6-1}, $7$-hedrite {\bf 7-1}, $6$-hedrite {\bf 4-1},
$5$-hedrite {\bf 3-1}, $4$-hedrite {\bf 3-1} are, respectively, $8$-hedrite
{\bf 12-4}, $7$-hedrite {\bf 14-9}, $6$-hedrite {\bf 8-3},
$5$-hedrite {\bf 6-2}, $4$-hedrite {\bf 4-1}.
The operation of taking the medial is a particular case of the {\em Goldberg-Coxeter construction} for the parameters $(k,l)=(1,1)$ (\cite{Gold37}, \cite{Cox71}, \cite{DD03}).

\section{Intersection of central circuits}

The following Theorem is a local version (for ``parts'' of the sphere) of
the Euler formula
$2p_2+p_3=8+ \sum_{i\geq 5} (i-4)p_i$ for $p$-vector of any $4$-valent 
plane $3$-connected graph $P$.

For any $4$-valent planar graph $P$, a {\em patch} $A$ is a region of $P$
bounded by $q$ arcs (paths of edges) belonging to central circuits 
(different or coinciding), such that all $q$ arcs form together a circle. 
A patch can be seen as a $q$-gon; we admit also 
$0$-gonal $A$, i.e. just the interior of a simple central circuit. 
Suppose that the patch $A$ is {\em regular}, i.e.
the continuation of any of bounding arc (on the central circuit to which it
belongs) lies outside of the patch.
See below two examples of patch.

\begin{center}
\epsfxsize=60mm
\epsfbox{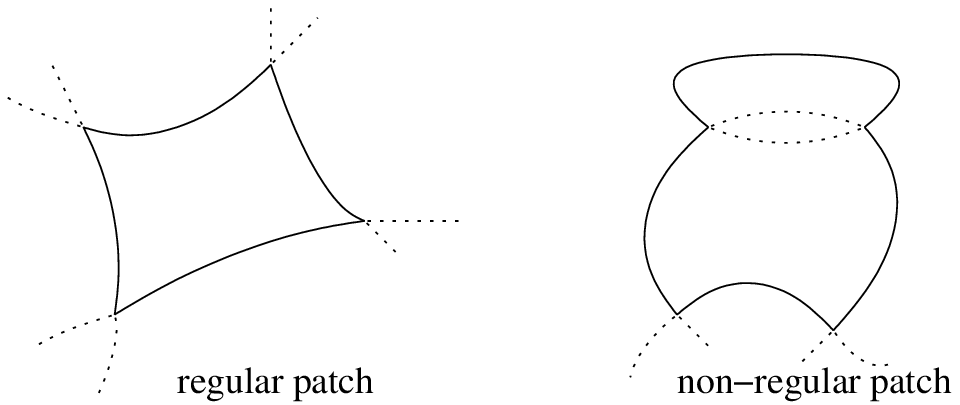}
\end{center}

Let $p'(A):=p'_1,...$ be the $p$-vector enumerating the faces of the 
patch $A$. The {\em curvature of the patch} $A$ is defined as
$c(A)=\sum_{k\geq 1} (4-k)p'_k$. So, the $k$-gon can be seen as,
respectively, positively curved, flat, or negatively curved, if
$k<4$, $k=4$, or $k>4$.

\begin{proposition}\label{Local-Euler-Formula}
(i) If $A$ be a regular patch, then $c(A)=4-q$, moreover:

(i.1) $c(A)=4$ if and only if $A$ is bounded by a simple central circuit.

(i.2) $c(A)=0$ if and only if $A$ is a rectangle formed by $4$-gons put 
together.

(ii) Any patch $A$ is the union of regular patches $A_1, \dots, A_p$; one has  $c(A)=(4-q_1)+\dots+(4-q_p)$, where each patch $A_i$ is bounded by $q_i$ arcs.

(iii) If a graph $G$ is the union of patches $A_1, \dots, A_p$, then $8=c(A_1)+\dots+c(A_p)$.
\end{proposition}
\proof (i) is a restatement in our terms of Theorem 1 of \cite{DSt}, (i1) and (i2) are easy consequences. The properties (ii) and (iii) follow from the definition of the curvature of a patch and from Euler formula.

\begin{proposition}
(i) Any $4$-valent plane graph whose faces are $k$-gons with $k$ even has central circuits with no self-intersection vertices.

(ii) At least one central-circuit of a $7$-hedrite self-intersects.
\end{proposition}
\proof In fact, if a central circuit of a $4$-hedrite self-intersects, then we have an $1$-gonal regular patch. The equality of above Theorem becomes $\sum_{i} (4-i)p'_i=3$, an impossibility since the left hand side is even.

Take a central circuit containing an edge of the unique $2$-gon, then the sequence (possibly empty) of adjacent $4$-gons will necessarily finish by a $3$-gon, or this $2$-gon; both cases yield a self-intersection.



Let us call {\em graph of curvatures} of an $i$-hedrite $G$, the
graph (possibly, with loops and multiple edges) having as vertex-set
all $2$-gons and $3$-gons of $G$. 
Two vertices (say, $2$- or $3$-gonal faces $F$ and $F'$ of $G$) of 
this $i$-vertex graph are adjacent if there is a pseudo-road 
connecting them. A {\em pseudo-road} is sequence of $4$-gons, say,
$F_1, \dots,F_l$, such that putting $F_0=F$ and $F_{l+1}=F'$, 
we have that any $F_k$ with $1\leq k\leq l$ is adjacent to $F_{k-1}$ 
and $F_{k+1}$ on opposite edges (cf. the definition of
a {\em rail-road} in next Section). Clearly, in the graph of curvatures, 
the vertices corresponding to $2$- and $3$-gons, have degree $2$ and $3$, 
respectively.

\begin{proposition}\label{intersec}
Let $C_1$, $C_2$ be any two central circuits of an $i$-hedrite. Then 
they are disjoint if and only if they are simple and there exist a 
ring of $4$-gons separating them.

\end{proposition}
\proof In fact, if both $C_1$ and $C_2$ are simple circuits, Theorem 
is evident: the curvature of the interior of a patch is $4$ and so, two
circuits are separated by $4$-gons only. Suppose that $C_1$ is 
self-intersecting. Then it has at least three regular patches and each 
of them has curvature at most $3$.
The circuit $C_2$, being disjoint with $C_1$, lies entirely inside one 
of those patches, say, $A$. So, all its $3$-gons and $2$-gons, except, 
possibly, those from its exterior patch, lie in $A$. So, $c(A)\geq 5$,
since the exterior patch of $C_2$ has curvature at most $3$.
It contradicts to the fact that $A$ has curvature at most $3$. 

\begin{remark}
Consider a $4$-valent plane graph $G$ having only one central circuit;
then, the set of faces of $G$ can be partitioned into two classes
${\cal C}_1$, ${\cal C}_2$ in chess manner. Every vertex $v$ is 
contained in two faces $F$ and $F'$. Also, the unique central circuit 
can be given an orientation, which induces an orientation on the set
of edges.

The vertex $v$ is incident to two edges of $F$, $e_1$ and $e_2$, 
and to two edges of $F'$, $e'_1$ and $e'_2$. If $e_1$, $e_2$ 
have both arrow pointing to the vertex or both arrow pointing out 
of the vertex, then the same is true for $e'_1$ and $e'_2$, and then,
we say that $v$ {\em belongs to Class I}.
Class I and its complement, Class II, form a bipartition of the
set of vertices of the knot; reversing orientation of the central
circuit or interchanging ${\cal C}_1$ and ${\cal C}_2$ does not 
change the bipartition. 

If the graph consists of $p$, $p\geq 2$, central circuits $C_1, \dots, C_p$,
then, one can put orientations on every central circuit and get
a bipartition of the set of vertices. But in that case the 
bipartition will depend on the chosen orientations.
\end{remark}

\section{Adding and removal central circuits}

%
%
%
%
%

The {\em deleting of a central circuit $C$} in an $i$-hedrite $G$ consists 
of removal of all edges and vertices contained in $C$. It produces a 
$4$-valent plane graph $P'$ having only $k$-gonal faces 
with $k \leq 4$. But since cases $k=0,1$ are possible, we do not
always obtain an $i$-hedrite.

The {\em cutting} of an $i$-hedrite $G$ consists of adding another central 
circuit to it. The faces of the new $i'$-hedrite $G'$ with $8\geq i'\geq i$ 
comes from the cutting of faces of $G$. This operation is only partially
defined, since arbitrary cutting can produce $k$-gons with $k>4$. The 
cutting of a $4$-gon in several $4$-gons (two, if the face
is traversed only once) is possible only if the $4$-gon is traversed 
on opposite edges. This corresponds to the notion of {\it shore-zone} 
in \cite{DSt}.
A cutting changes CC-partition of an $i$-hedrite only in the following 
way: new central circuit $C$ is added and all others central circuits 
remain unchanged, except that the length of each of them increases by one 
for any intersection with $C$.

Call a {\em rail-road} a circuit of $4$-gons, possibly self-intersecting, 
in which every $4$-gon is adjacent to two of its neighbors on opposite
edges. A rail-road is bounded by two ``parallel'' central circuits.
The deleting of  one of those central circuits (in other words, collapsing
rail-road into one central circuit) is called {\em reduction}.
The cutting produces a rail-road if and only if it is an {\em inflation 
along a central circuit $C$}, i.e. replacing it by (thin enough) 
rail-road. A {\em $t$-inflation along a central circuit $C$} is replacing
this central circuit by $t-1$ parallel (thin enough) rail-roads.
A {\em $t$-inflation of an $i$-hedrite} is new $i$-hedrite
obtained from original one by simultaneous $t$-inflation along all
of its central circuits. 
A $t$-inflation of $G$ is $G$ if $t=1$, and it is just inflation of $G$ 
if $t=2$.


An $i$-hedrite is called {\em irreducible} if it contains no 
rail-road. It is called {\em maximal irreducible} if it cannot be
obtained from another $i'$-hedrite by a cutting.


\begin{remark}

Let $C$ be a central circuit of $G$ with 
$CC(G)=(...,a_i^{\alpha_i},...;...,b_j^{\beta_j},...)$, and let 
$Int(C)=(c_0;c_1^{\gamma_1},...,c_r^{\gamma_r})$. The $t$-inflation of $G$
denoted by $G^t$ has 
$CC(G^t)=(...,ta_i^{t\alpha_i},...;...,tb_j^{t\beta_j},...)$; if $C'$
is one of $t$ parallel copies of $C$, then 
$Int(C')=(c_0;c_1^{t\gamma_1},...,c_r^{t\gamma_r}, (2c_0)^{t-1})$.
\end{remark}

\section{Connectivity of $i$-hedrites}
For any integer $m\geq 2$ denote:

by $I_{6,2m}$ the $2m$-vertex $6$-hedrite, such that each $2$-gon is adjacent to two $3$-gons;

by $I_{5,2m+1}$ the $(2m+1)$-vertex $5$-hedrite, such that 
two $2$-gons share a vertex and remaining $2$-gon is adjacent
to two $3$-gons;

by $I_{4,2m+2}$ the $(2m+2)$-vertex $4$-hedrite, such that 
four $2$-gons are organized into two pairs sharing a vertex; 

by $J_{4,2m}$ the $m$-inflation of only one central circuit of
$4$-hedrite {\bf 2-1}; they are projections of {\em composite}
alternating links $2^2_1\#2^2_1\#2^2_1\dots\#2^2_1$ ($m$ times),
which we denote by $m\times 2^2_1$.

See in Table \ref{FundamentalInfo} the first occurrences 
(for $2\leq m\leq 5$) of those graphs, followed by their symmetry 
groups and CC-vectors.

\begin{table}
\centering
\epsfxsize=100mm
\epsfbox{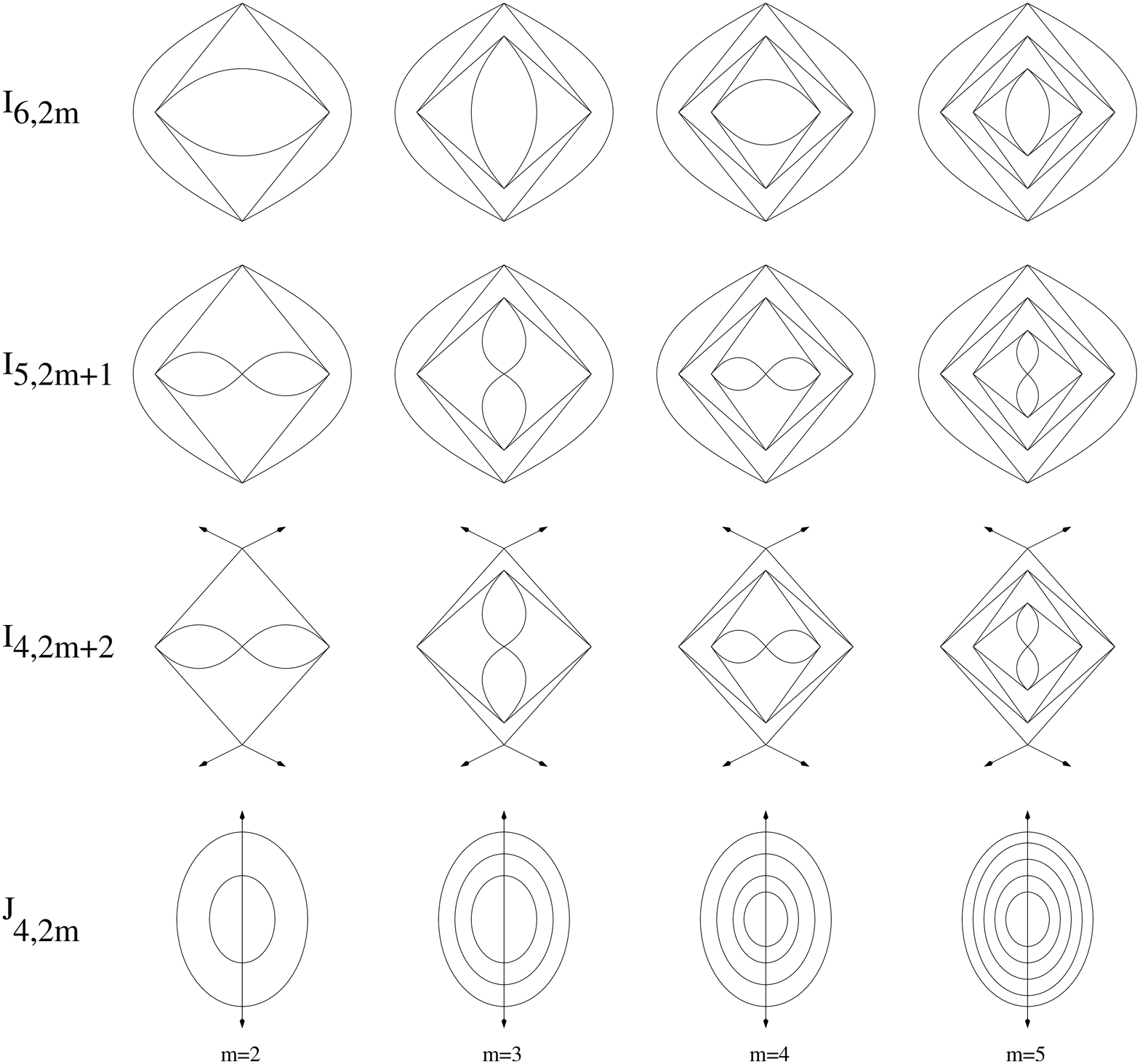}\\[2mm]
\begin{tabular}{||c|c|c||}
\hline\hline
$i$-hedrite            &Group      &CC-vector\\\hline\hline
$I_{6,2m}$, $m$ even   &$D_{2d}$   &$4m$\\\hline
$I_{6,2m}$, $m$ odd    &$D_{2h}$   &$(2m)^2$\\\hline
$I_{5,2m+1}$           &$C_{2v}$   &$4m+2$\\\hline
$I_{4,2m+2}$, $m$ even &$D_{2d}$   &$(2m+2)^2$\\\hline
$I_{4,2m+2}$, $m$ odd  &$D_{2h}$   &$(2m+2)^2$\\\hline
$J_{4,2m}$             &$D_{2h}$   &$2^m, 2m$\\\hline\hline
\end{tabular}
\caption{All $i$-hedrites, which are {\em not} $3$-connected}
\label{FundamentalInfo}
\end{table}

\begin{lemma}
Any $i$-hedrite is $2$-connected.
\end{lemma}
\proof Let $G$ be an $i$-hedrite and assume that there is one vertex $v$,
such that $G-\{v\}$ is disconnected in two components $C_1$ and $C_2$.
Then two edges from $v$ will connect to a vertex $w$ of $C_1$ and two
edges from $v$ will connect to a vertex $w'$ of $C_2$, because, otherwise,
the exterior face is $m$-gonal with $m>4$. See below the corresponding
drawing.

\begin{center}
\epsfxsize=60mm
\epsfbox{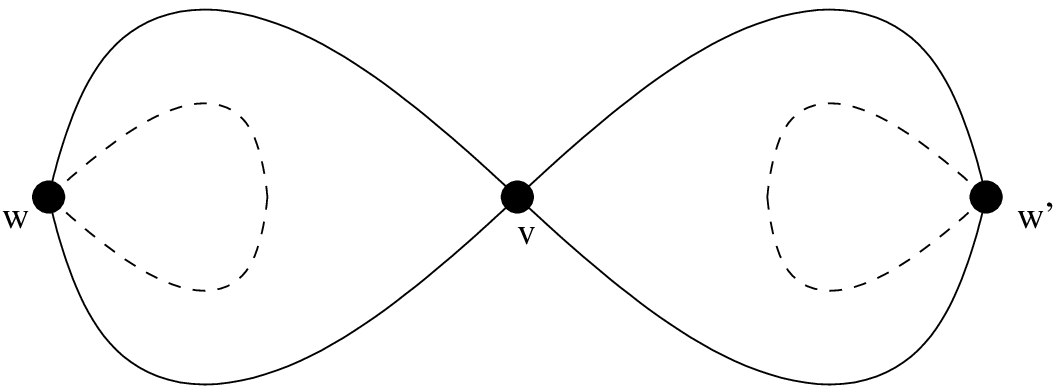}
\end{center}

But the vertex $w$ will disconnect the graph and so, iterating the construction, we obtain an infinite sequence $v_1, \dots, v_n$ of vertices that disconnects $G$. 
This is impossible since $G$ is finite.

Any $i$-hedrite (moreover, any Eulerian graph) has at least one 
Eulerian circuit of edges; so, there is no cut-edge.
But a cut-vertex appears already for some Eulerian ($3$-vertex)
and $4$-valent ($4$-vertex) plane graphs.

\begin{theorem}\label{3-connectedness}
Any $i$-hedrite, which is not $3$-connected, is one of $I_{6,2m}$, $I_{5, 2m+1}$, $I_{4, 2m+2}$, $J_{4, 2m}$ for some $m\geq 2$.

\end{theorem}

\proof Let $G$ be an $i$-hedrite and assume that it is not $3$-connected. 
Then there are two vertices, say, $v$ and $v'$, such that 
$G-\{v, v'\}$ is disconnected in two components, say, $C_1$ and $C_2$.
Amongst the $4$ edges from $v$ (respectively $v'$), the edges 
$\{e_1,\dots, e_{s}\}$ (respectively $\{e'_1,\dots, e'_{s'}\}$) go
to $C_1$. Two numbers $s$ and $s'$ can takes values $1$, $2$ or $3$; we
will consider all possible cases.

\begin{figure}
\centering
\mbox{
\subfigure[The case $s=1$, $s'=1$]{\epsfig{figure=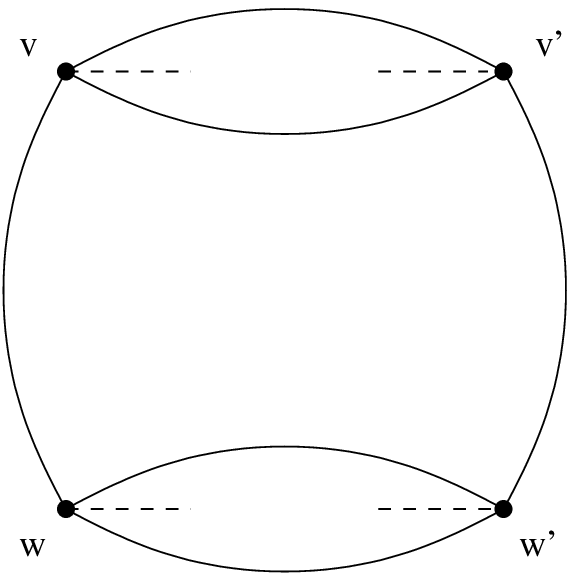,width=.25\textwidth}}
\subfigure[The case $s=2$, $s'=2$]{\epsfig{figure=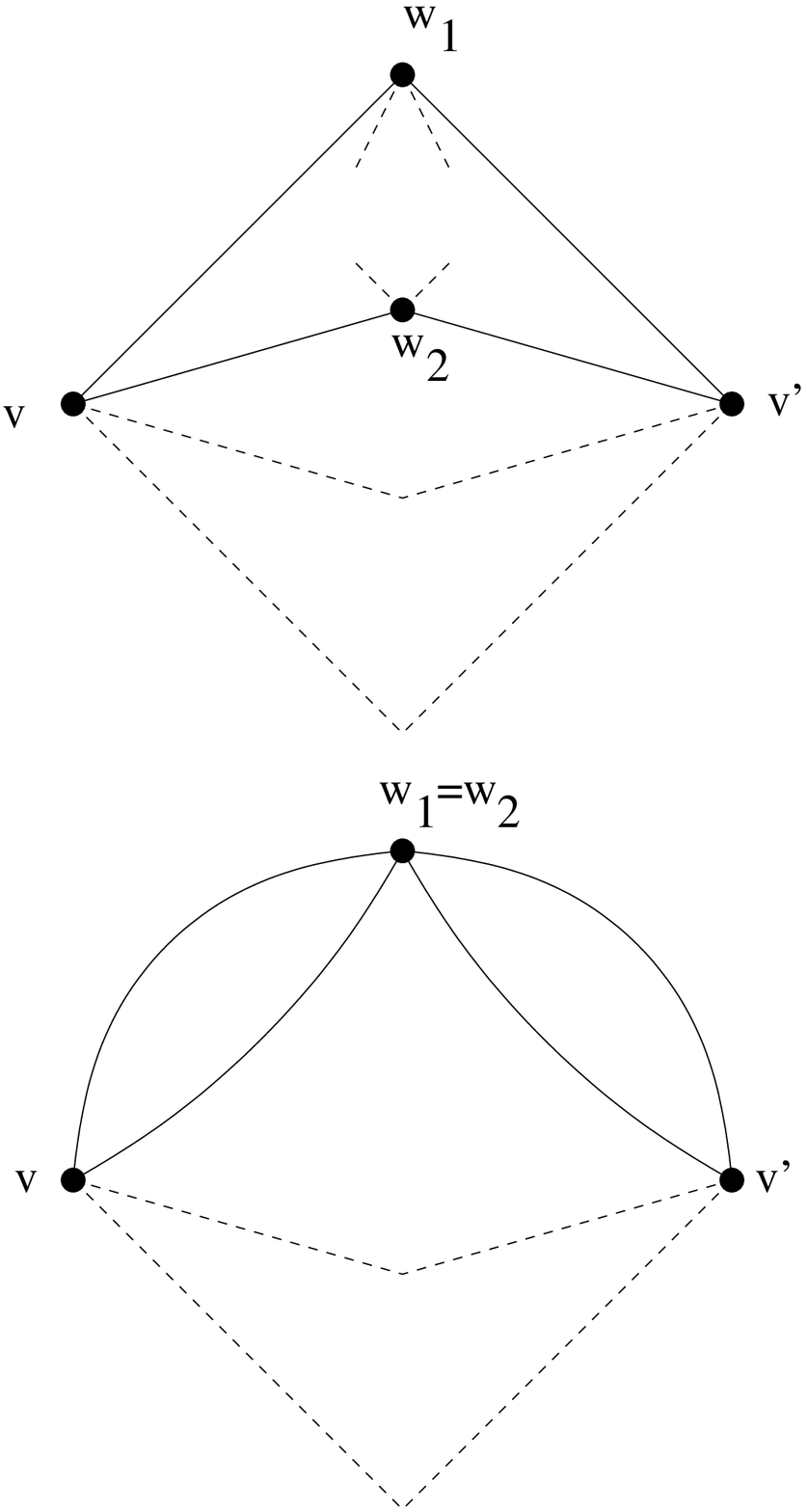,width=.25\textwidth}}
\subfigure[The case $s=2$, $s'=1$]{\epsfig{figure=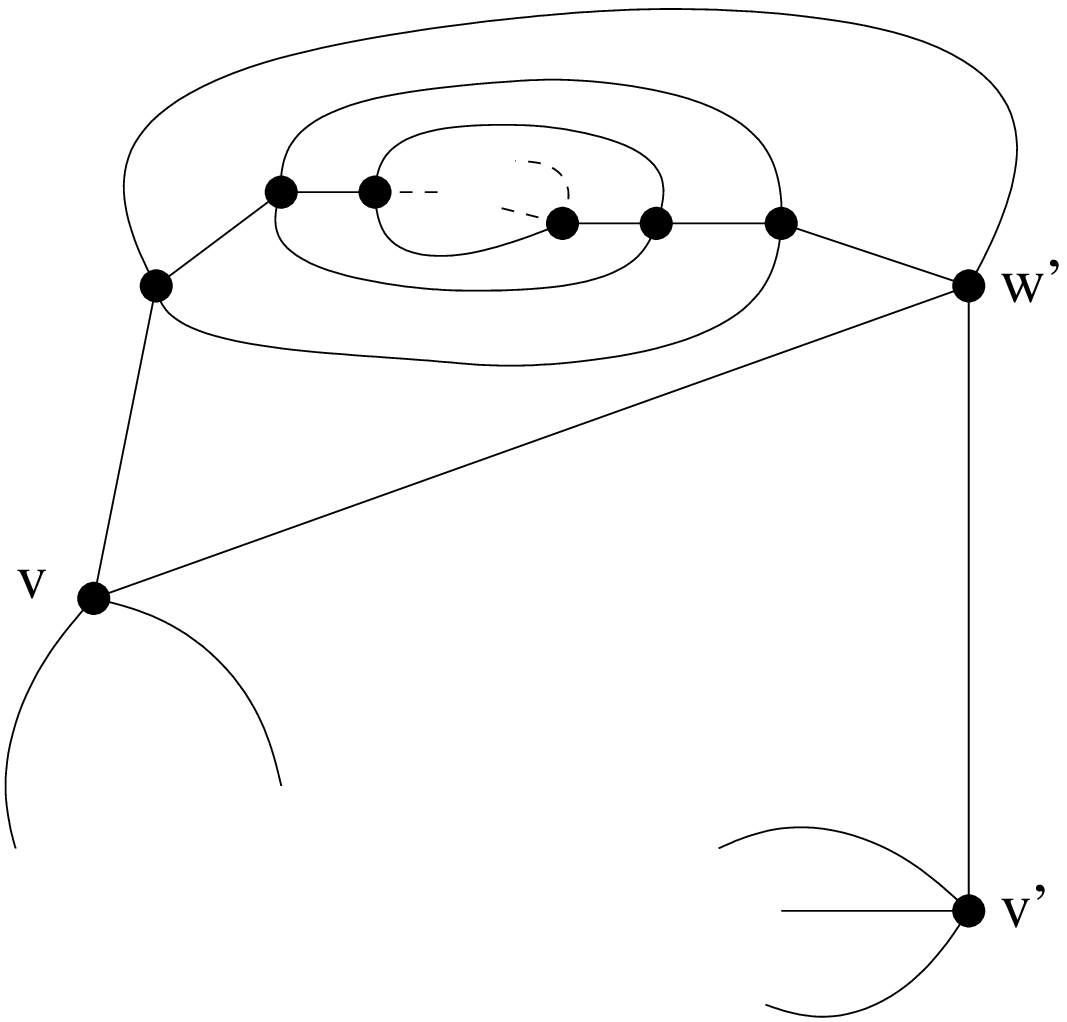,width=.30\textwidth}}
}
\caption{The three cases of Theorem \ref{3-connectedness}}
\label{fig:TheThreeCases}
\end{figure}

Assume that $s=1$ and $s'=1$, then the edges $e$ and $e'$ must be
distinct, since, otherwise, $C_1$ is the empty graph. Moreover, $e$ and
$e'$ have no common vertices, since, otherwise, $G$ would not be $2$-connected.
So, $v$ and $v'$ are connected by $e$ and $e'$ to a vertex $w$ and $w'$,
respectively. Since face-size is at most $4$, the vertices $v$ and $v'$,
(respectively, $w$ and $w'$) are linked by two edges
(see Figure \ref{fig:TheThreeCases}).
Two points $w$ and $w'$ can either be connected by two edges and
we are done, or disconnect the graph. In the latter case, we can iterate
the construction. Since the graph is finite, the construction eventually
finish and we get a graph $J_{4,2m}$. If $s=1$ and
$s'=3$, then by a similar reasoning, one gets again a graph $J_{4,2m}$.

Assume that $s=2$ and $s'=2$. One has $\{e_1, e_2\}\cap \{e'_1, e'_2\}=\emptyset$, since, otherwise, one can attribute an edge to $C_2$ and get the case
$s=1$ and $s'=1$. So, one has, say, $e_1\cap e'_1=\{w_1\}$ and 
$e_2\cap e'_2=\{w_2\}$, and the following two possibilities 
(see Figure \ref{fig:TheThreeCases}): either $w_1=w_2$ (this
corresponds to $\{e_1, e_2\}$ and $\{e'_1, e'_2\}$ being the edges
of two $2$-gons) and we are done, or $w_1\not= w_2$.
Assume now that $w_1\not= w_2$; two points $w_1$ and $w_2$ can 
either be connected by two edges and we are done, or disconnect the
graph. In the latter case, we can iterate the construction. Since the 
graph is finite, the construction eventually finish. If we do the same
construction on the other side, then we get a similar structure and the 
graph is of the form $I_{4,2m+2}$, $I_{5, 2m+1}$ or $I_{6,2m}$ with 
$m\geq 2$.

Assume now that $s=2$ and $s'=1$. The edges $e_1$, $e_2$ and $e'_1$ are
all distinct, since, otherwise, the vertex $v$ disconnects the graph.
So, $v'$ is connected by $e'_1$ to a vertex $w'$. Now, either $v'$ or $w'$
is connected to $v$, since, otherwise, we would have a $5$-gonal face.
If $w'$ is connected to $v$, then the pair $\{w', v\}$ disconnects the graph.
This construction is infinite (see Figure \ref{fig:TheThreeCases}); so,
we get a contradiction.

\begin{theorem}
(i) If $G$ is an $i$-hedrite with two adjacent $2$-gons, then 
this is a $4$-hedrite {\bf 2-1} or a $J_{4,2m}$ with $m\geq 2$.

(ii) If $G$ is an $i$-hedrite with two $2$-gons sharing a vertex, then 
it is either $4$-hedrite {\bf 4-1}, or an $I_{4,2m+2}$, or $5$-hedrite 
{\bf 3-1}, or an $I_{5,2m+1}$ with $m\geq 2$.

\end{theorem}

\proof The proof for (i), (ii) is similar to the cases $(s,s')$=$(1,1)$, $(2,2)$ of Theorem \ref{3-connectedness}.

\begin{theorem}

(i) an $n$-vertex $4$-hedrites exists if and only if $n\geq 2$, even.

(ii) an $n$-vertex $5$-hedrites exists if and only if $n\geq 3$, $n\not= 4$.

(iii) an $n$-vertex $6$-hedrites exists if and only if $n\geq 4$.

(iv) an $n$-vertex $7$-hedrites exists if and only if $n\geq 7$.

(v) an $n$-vertex $8$-hedrites exists if and only if $n\geq 6$, $n\not= 7$.

\end{theorem}
\proof The case (v) is proven in \cite{Gr}, page 282.
The case (i) is trivial; take, for example, the serie $J_{4,2m}$.

The series $I_{6,2m}$, $I_{5,2m+1}$ for any $m \ge 2$, give $6$-hedrites
and $5$-hedrites with even and, respectively, odd number of vertices.

For $5$-, $6$- and $7$-hedrites, we 
get by $(t-1)$-inflation
along a central circuit of length $4$ in corresponding
$i$-hedrite {\bf 6-1}, {\bf 5-1} and {\bf 7-1},
series with $4t+2$, $4t+1$ and $4t+3$ vertices. 
By $(t-2)$-inflation
along such central circuit in $7$-hedrite {\bf 8-1}, we get a serie of
$7$-hedrites with $4t$ vertices.
By $(t-1)$-inflation along central circuit in $6$-hedrite {\bf 11-4}, we
get serie of $6$-hedrites with $8t+3$ vertices.
By $t$-inflation along central circuit of length $8$ in 
$6$-hedrite {\bf 15-10}, we get serie of $6$-hedrites with $8t+7$ vertices.

Inscribing consecutively $4$-gons in the $4$-gon, which is adjacent only to
$3$-gons, in $7$-hedrite {\bf 9-1} and {\bf 10-2}, we get series of 
$7$-hedrites for the remaining cases of $4t+1$ and $4t+2$ vertices.

In case of $5$-hedrites, it remains to prove existence for the 
case $n=4t>1$. We obtain existence in the sub-cases 
$n=m\times 4^b$, where $m\geq 3$ and not divisible by $4$, 
$n=8\times 4^b$, and $n=16\times 4^b$, respectively:

by $b$-inflation of some $m$-vertex $5$-hedrite; we showed their existence,

by $b$-inflation of $5$-hedrites {\bf 8-1},

by $b$-inflation of any $16$-vertex $5$-hedrite (for example, one coming from {\bf 10-2} by inflation along central circuit of length $6$).

Our computation (see the last Section) present all $i$-hedrites with at most $15$ vertices.


\section{Irreducible $i$-hedrites}

\begin{theorem}\label{irre}
Any irreducible $i$-hedrite has at most $i-2$ central circuits and equality is attained for each $i$, $4\leq i\leq 8$.
\end{theorem}
\proof For $i=8$, the Theorem is proved in \cite{DSt}. We will show, using suitable cutting, that this result implies the Theorem for $i<8$.

Let us start with the simplest case of $7$-hedrites. Consider a simple circuit $S$ in its curvature graph, which contains the vertex corresponding to unique $2$-gon. Remind that the vertices in the curvature graph correspond to $2$- or $3$-gons, while edges correspond to pseudo-roads. So, the simple circuit $S$ corresponds to the circuit of faces of $G$, containing our $2$-gon, some $4$-gons and, possibly, some of six $3$-gons; see the picture below.

\begin{center}
\epsfxsize=40mm
\epsfbox{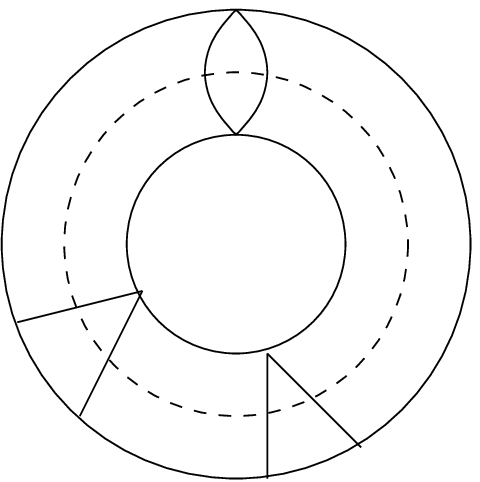}
\end{center}

Suppose that the $7$-hedrite $G$ is irreducible and has $k$ central circuits. By adding the central circuit $C$, which is shown by dotted lines on picture above, we produce an $8$-hedrite (since, the $2$-gon is cut by $C$ in two $3$-gons), which is still irreducible and has $k+1$ central circuits. So, $k+1\leq 6$, by Theorem 3 of \cite{DSt}.

For remaining cases of $i$-hedrites with $i=4,5,6$, the proof is similar. In each case, we consider all possible distribution of $2$-gons by simple circuits in the graph of curvatures and, for each such circuit, we add suitable number of new central circuits.

\begin{figure}
\centering
\mbox{\subfigure[The two cases for $6$-hedrites]{\epsfig{figure=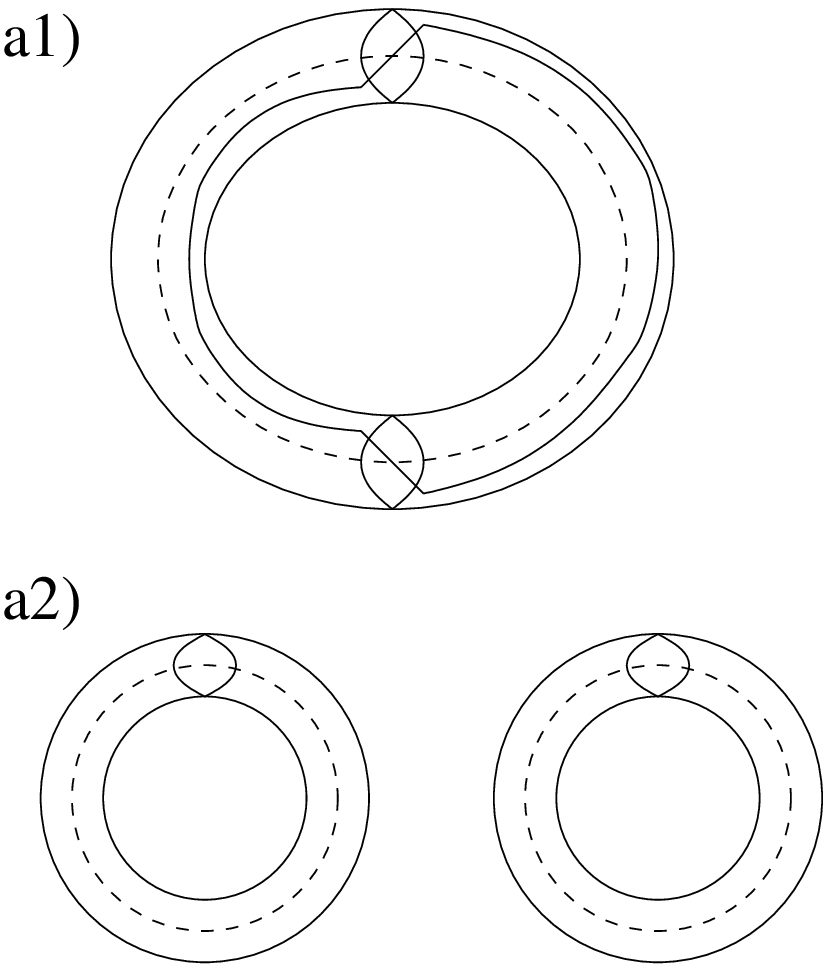,width=.25\textwidth}}\quad
\subfigure[The three cases for $5$-hedrites]{\epsfig{figure=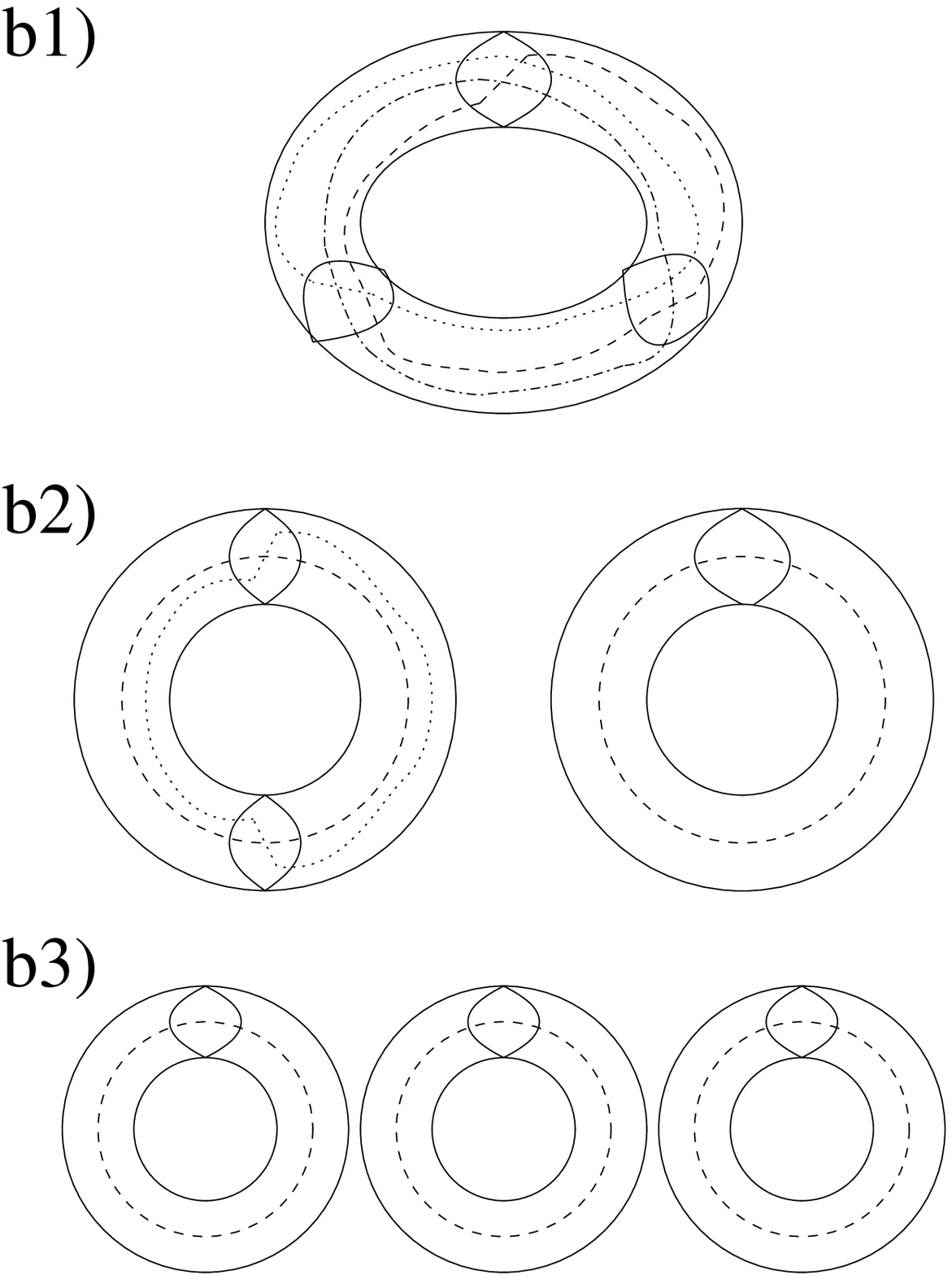,width=.25\textwidth}}\quad
\subfigure[The two cases for $4$-hedrites]{\epsfig{figure=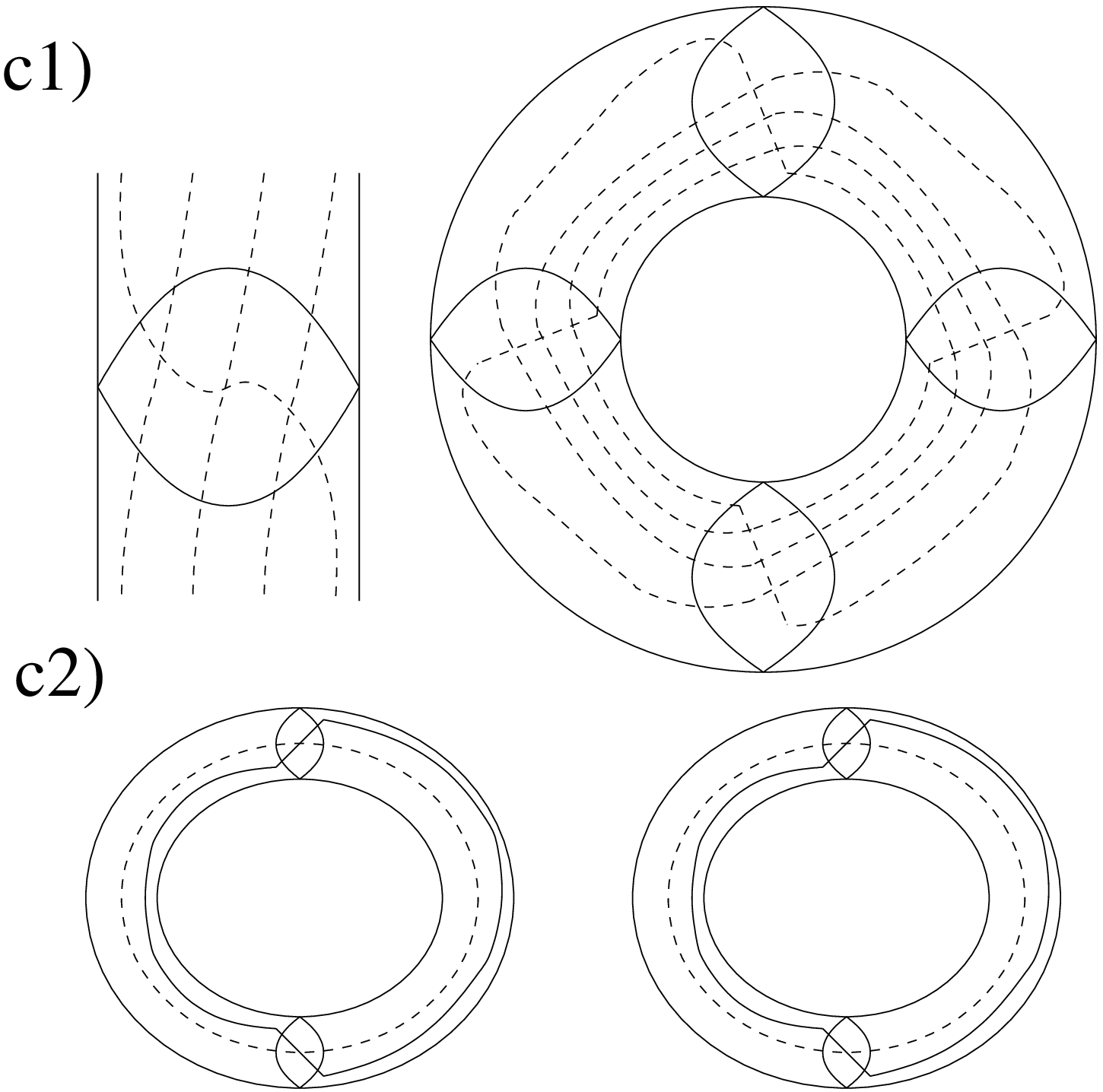,width=.30\textwidth}}}
\caption{The construction of Theorem \ref{irre} for $6$-, $5$-, $4$-hedrites}
\label{The456hedriteCases}
\end{figure}

All possibilities are presented on Figure \ref{The456hedriteCases}:
two for $i=6$, three for $i=5$ and two for $i=4$. In the last case, there
are no $3$-gons and so, simple circuits in the curvature graph contain
only even number of $2$-gons by local Euler formula of Theorem 
\ref{Local-Euler-Formula}. Note that the case of $4$-hedrites is
obvious by Theorem 5 of \cite{DSt}.


\begin{figure}
\centering
\mbox{
\subfigure[The first case]{\epsfig{figure=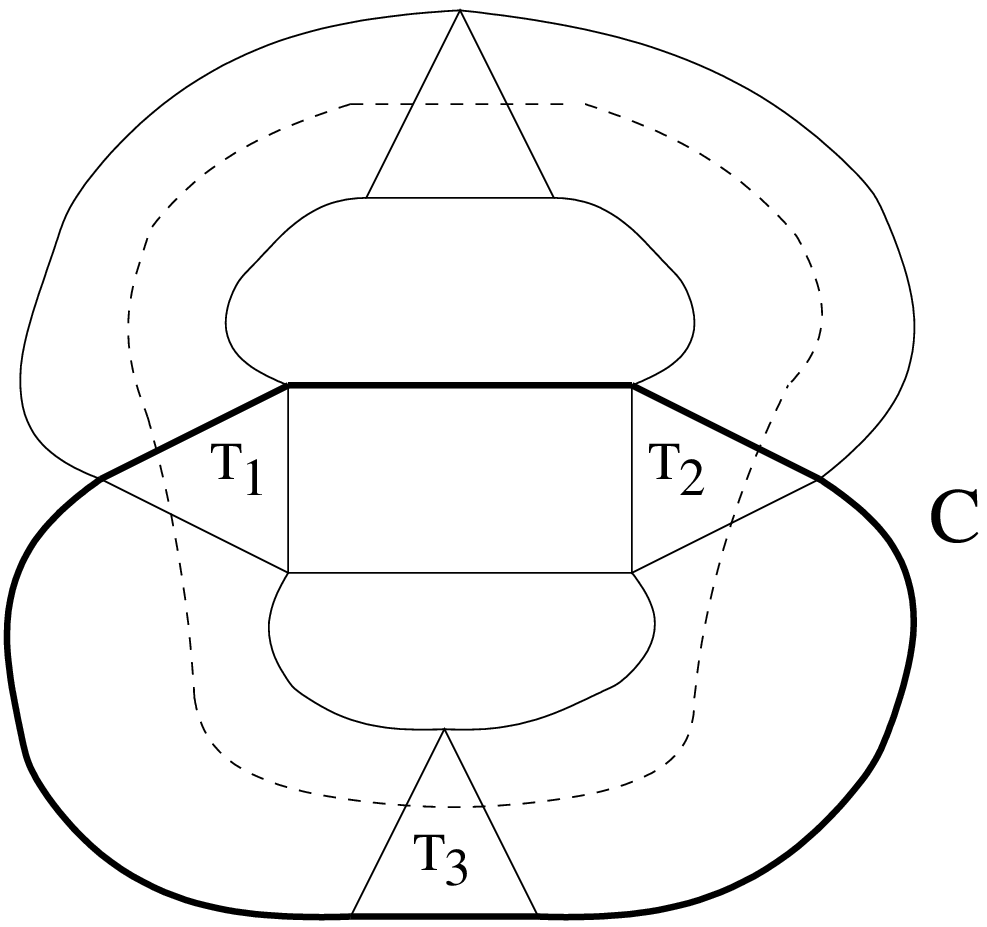,width=.25\textwidth}}
\subfigure[The second case]{\epsfig{figure=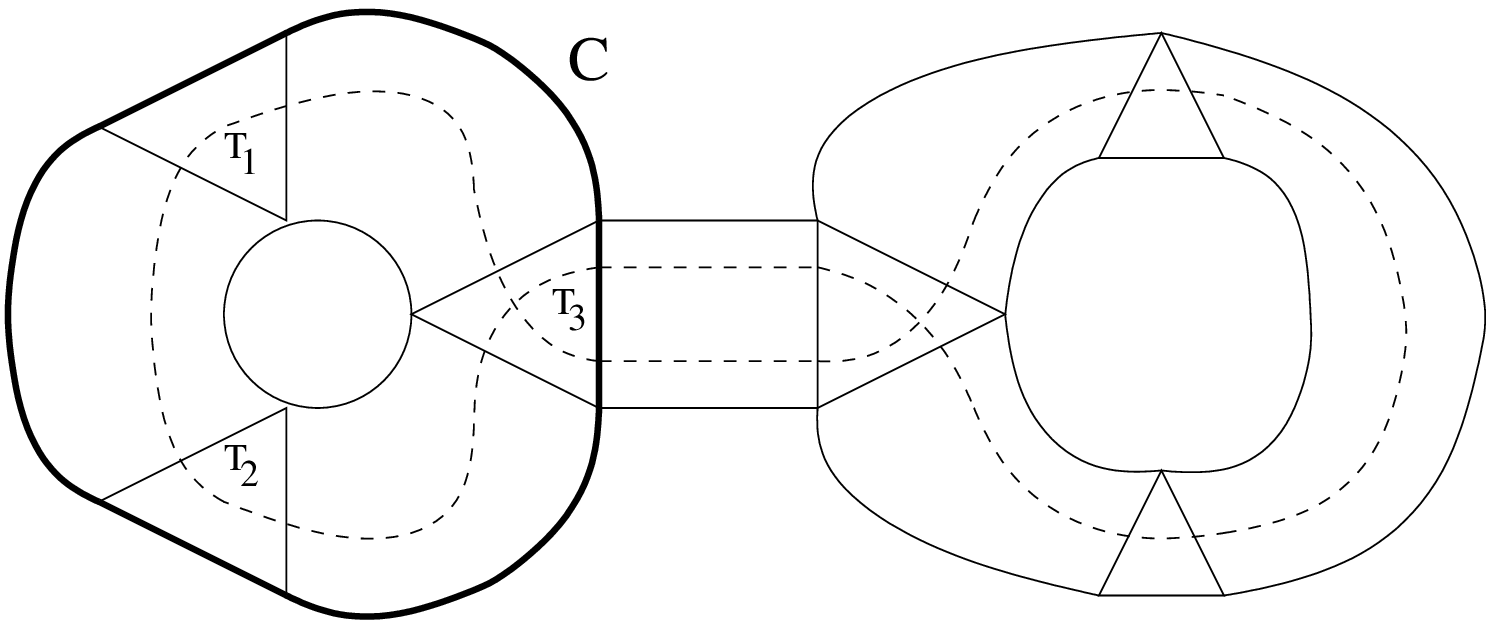,width=.50\textwidth}}
}
\caption{The two cases of cutting of irreducible $8$-hedrite}
\label{TwoCasesCutting}
\end{figure}

\begin{lemma}
Let $G$ be an irreducible $8$-hedrite and $C$ a central circuit,
which is incident to three $3$-gons on one side. Then one can add another central circuit to $G$,  so that the resulting graph is still irreducible.
\end{lemma}
\proof From every one of the three $3$-gons, say, $T_1$, $T_2$, $T_3$, one can define two pseudo-roads from the sides of the $3$-gon, which do not belong to the central circuit.
Each such pseudo-road defines an edge, say, $e_1$, $e_2$, $e_3$ in the graph 
of curvatures and so, a triangle in that graph. Then, either two triangles 
$T_i$ and $T_j$ are linked by a path, which does not involve the edges $e_k$,
or they are not linked by such a patch. In both cases we can cut the 
$8$-hedrite according to Figure \ref{TwoCasesCutting} and obtain 
another $8$-hedrite, which is still irreducible.

See below example of an irreducible $7$-hedrite (its CC-vector is $(10^2, 12^2; 20)$, its symmetry group is $C_{2v}$) with the maximum number of central circuits.

\begin{center}
\epsfxsize=60mm
\epsfbox{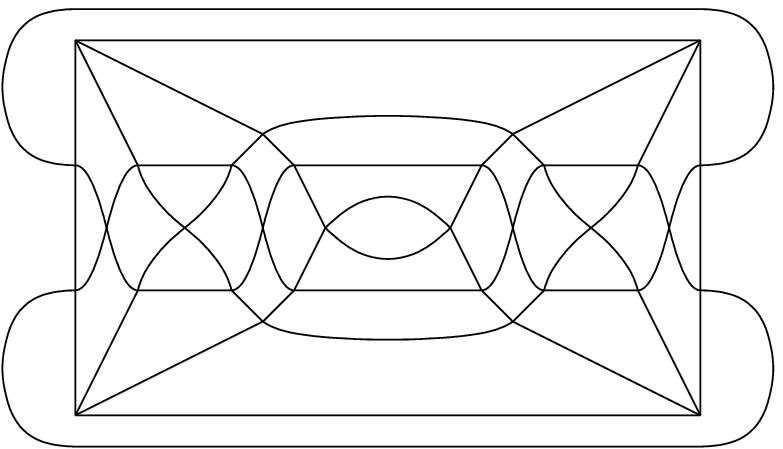}
\end{center}

For all other $i$, there is an example of irreducible $i$-hedrite with $i-2$ central circuits (see 
Theorem \ref{TheOneWithSimpleCentralCircuit}), which is, moreover, pure.

\begin{conjecture}
An irreducible $i$-hedrite is maximal irreducible if and only if it has $i-2$ central circuits.
\end{conjecture}

\section{Classification of pure irreducible $i$-hedrites}
The easiest case, $i=4$, of i-hedrites admits following complete
characterization:

\begin{theorem}\label{Theorem-for-4-hedrite}
(i) Any $4$-hedrite can be obtained from some $4$-hedrite with two central
circuits by simultaneous $t_1$- and $t_2$-inflation along those circuits; it is
irreducible if and only if $t_1=t_2=1$.

(ii) Any $4$-hedrite with two central circuits is defined by
its number of vertices $n$ and by {\em shift} $j$, $0 \le j \le n/4$
with $gcd(n/2, j)=1$,
vertices between the pair of boundary $2$-gons on the horizontal
circuit (see, for example, $4$-hedrite {\bf 8-1}) and the remaining
pair of $2$-gons. Remark that several different values of shift can yield
the same graph.

(iii) Any $4$-hedrite is balanced.

\end{theorem}
\proof (i) and (ii) are proved in \cite{DSt}, 
while (iii) is obvious for $4$-hedrite with two central circuits
and remain true under $t_1$- and $t_2$-inflation.

The shift $j=0$ corresponds to {\bf 2-1} and its only-on-one-circuit
$m$-inflations $J_{4,n=2m}$. The shift $j=1$ corresponds to {\bf 4-1} and 
$I_{4,n=2m+2}$.
Denote by $K_{4,4m}$, for any $m \ge 2$, any $4m$-vertex $4$-hedrite obtained
from {\bf 4-1} by $m$-inflation of only one its central circuit; so, its 
CC-vector is $(4^m,4m)$, its symmetry is $D_{2d}$ and it is reducible.
Clearly, any $K_{4,n=4m}$ has the maximal shift $j=n/4$.

%

\begin{figure}
{\small
\setlength{\unitlength}{1cm}
\begin{minipage}[t]{3.5cm}
\centering
\epsfxsize=2.5cm
\epsffile{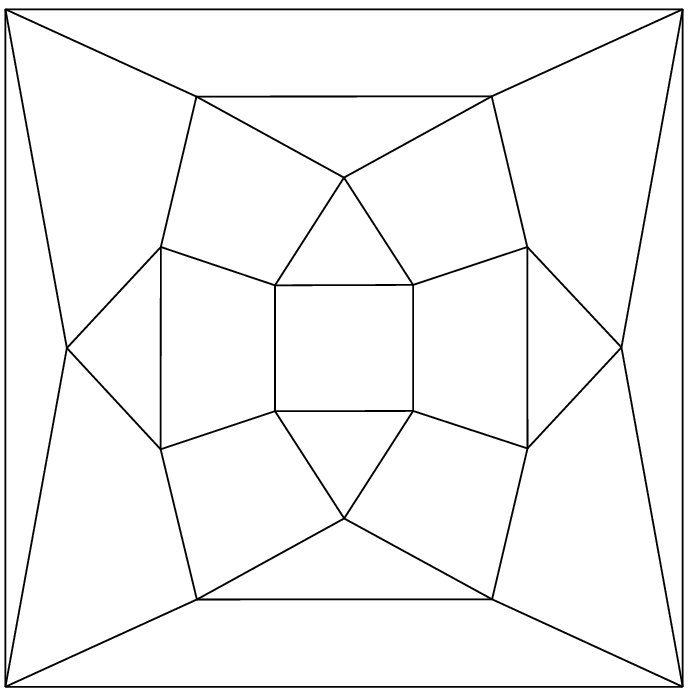}\par
{{\bf Nr.20-1} \quad $D_{2d}$ \\ $(8^5)$ \\ }
\end{minipage}
\begin{minipage}[t]{3.5cm}
\centering
\epsfxsize=2.5cm
\epsffile{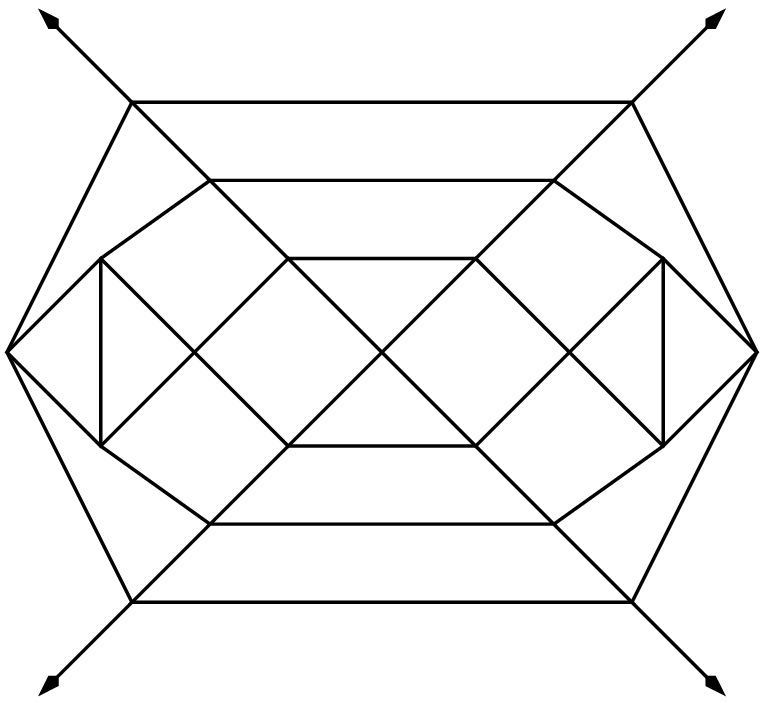}\par
{{\bf Nr.22-1} \quad $D_{2h}$ \\ $(8^3,10^2)$ \\ }
\end{minipage}
\setlength{\unitlength}{1cm}
\begin{minipage}[t]{3.5cm}
\centering
\epsfxsize=2.3cm
\epsffile{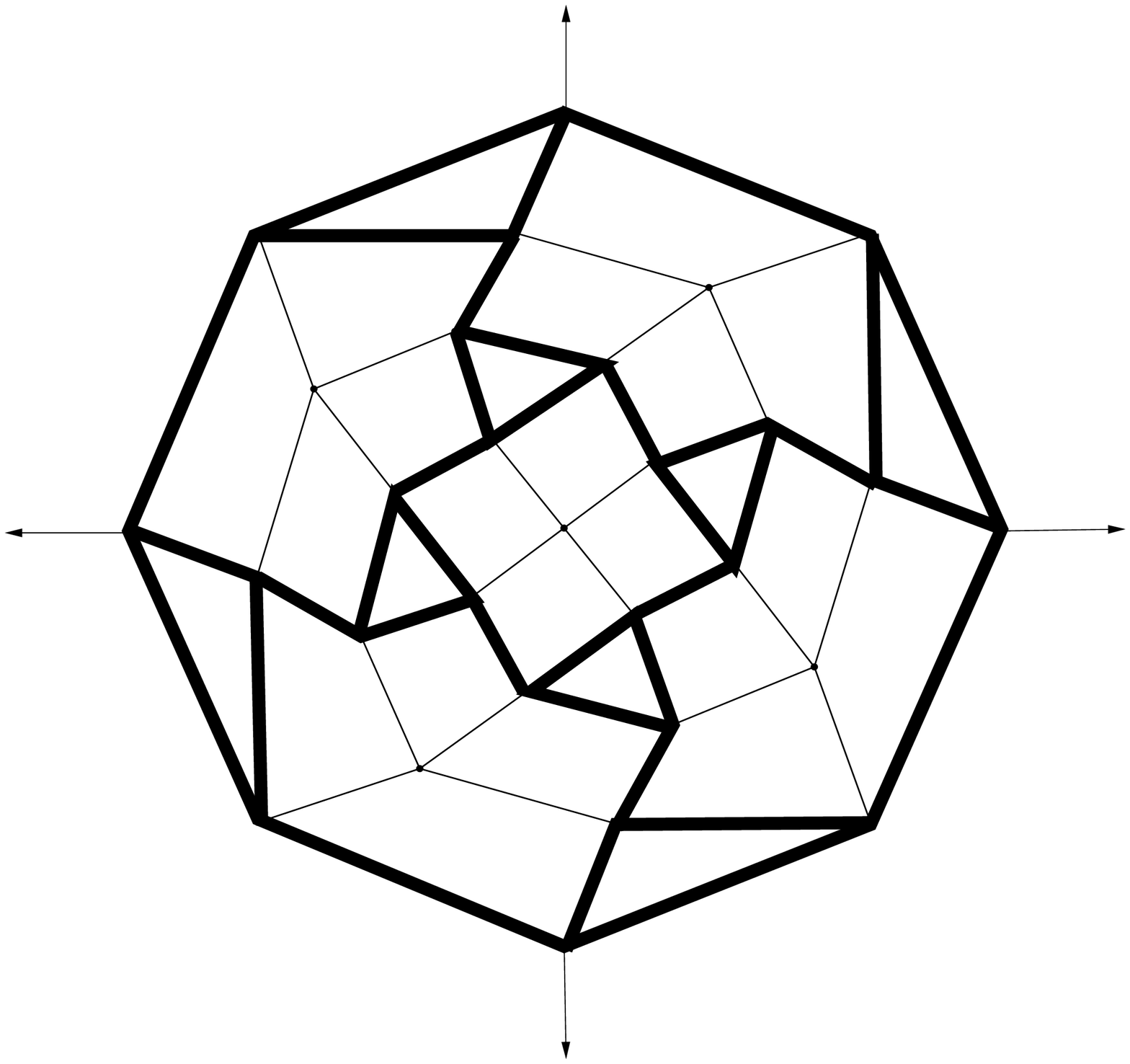}\par
{{\bf Nr.30-1} \quad $O$ \\ $(10^6)$ \\}
\end{minipage}
\setlength{\unitlength}{1cm}
\begin{minipage}[t]{3.5cm}
\centering
\epsfxsize=2.3cm
\epsffile{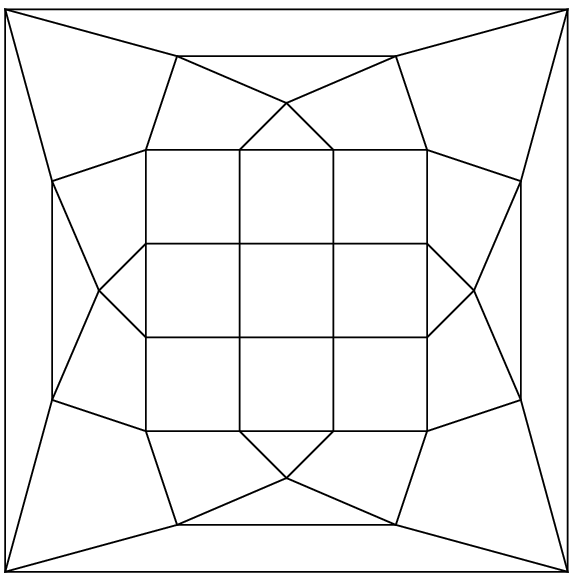}\par
{{\bf Nr.32-1} \quad $D_{4h}$ \\ $(10^4,12^2)$ \\}
\end{minipage}
}
\caption{Pure irreducible $i$-hedrites with $5$ or $6$ central circuits}
\label{ThePureIrreducibleOctahedriteWith56CC}
\end{figure}

\begin{theorem}\label{TheOneWithSimpleCentralCircuit}

Any pure irreducible $i$-hedrite is, either any $4$-hedrite with two
central circuits, or a $5$-hedrite {\bf 6-2}, or one of $6$-hedrites
{\bf 8-6}, {\bf 14-20}, or one of the following eight $8$-hedrites:
{\bf 6-1}, {\bf 12-4}, {\bf 12-5}, {\bf 14-7}, and (see Figure 
\ref{ThePureIrreducibleOctahedriteWith56CC}) {\bf 20-1}, {\bf 22-1},
{\bf 30-1}, {\bf 32-1}.

\end{theorem}

\proof Let $G$ be a pure irreducible $i$-hedrite having $r$ central circuits. 
If one deletes a central circuit, then, in general, $1$-gon can appear. It 
does not happen for $G$, since it would imply a self-intersection of a central 
circuit. So, the result of deletion of a central circuit from $G$ 
produces an pure irreducible $i$-hedrite with $r-1$ central circuits.

First, if $r=2$, then the Theorem 5 from \cite{DSt} gives that such $G$ are 
exactly $4$-hedrites with two central circuits; all of them are classified 
in Theorem \ref{Theorem-for-4-hedrite}.

We prove the Theorem by systematic analysis of all possible ways to add 
to $G$ (for $r=2,3,4,5$) a central circuit, in order to get a pure 
irreducible $i$-hedrite with $r+1$ central circuits.

\begin{figure}
\centering
\epsfxsize=55mm
\epsfbox{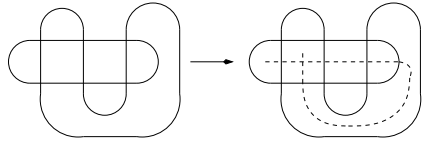}
\caption{No {\em pure} irreducible $i$-hedrite can be obtained by cutting of the above $4$-hedrite}
\label{Cutting4hedrite}
\end{figure}

Let $r=2$. Then $G$ can be only one of two smallest $4$-hedrites. In 
fact, if $G$ is another $4$-hedrite, then, because of classification Theorem 
\ref{Theorem-for-4-hedrite}, it has a form as in Figure \ref{Cutting4hedrite}.

New central circuit should cut both $2$-gons on opposite edges, since, otherwise, there is a rail-road. But Figure \ref{Cutting4hedrite} shows, on example for $n=6$, that a self-intersection appears if two central circuits intersect in more than four vertices.

So, the only possible $4$-hedrites with {\em two} central circuits,
which can be cut in order to produce irreducible pure $i$-hedrite
are $4$-hedrites {\bf 2-1} and {\bf 4-1}. All cases are indicated below. 

\begin{center}
\epsfxsize=120mm
\epsfbox{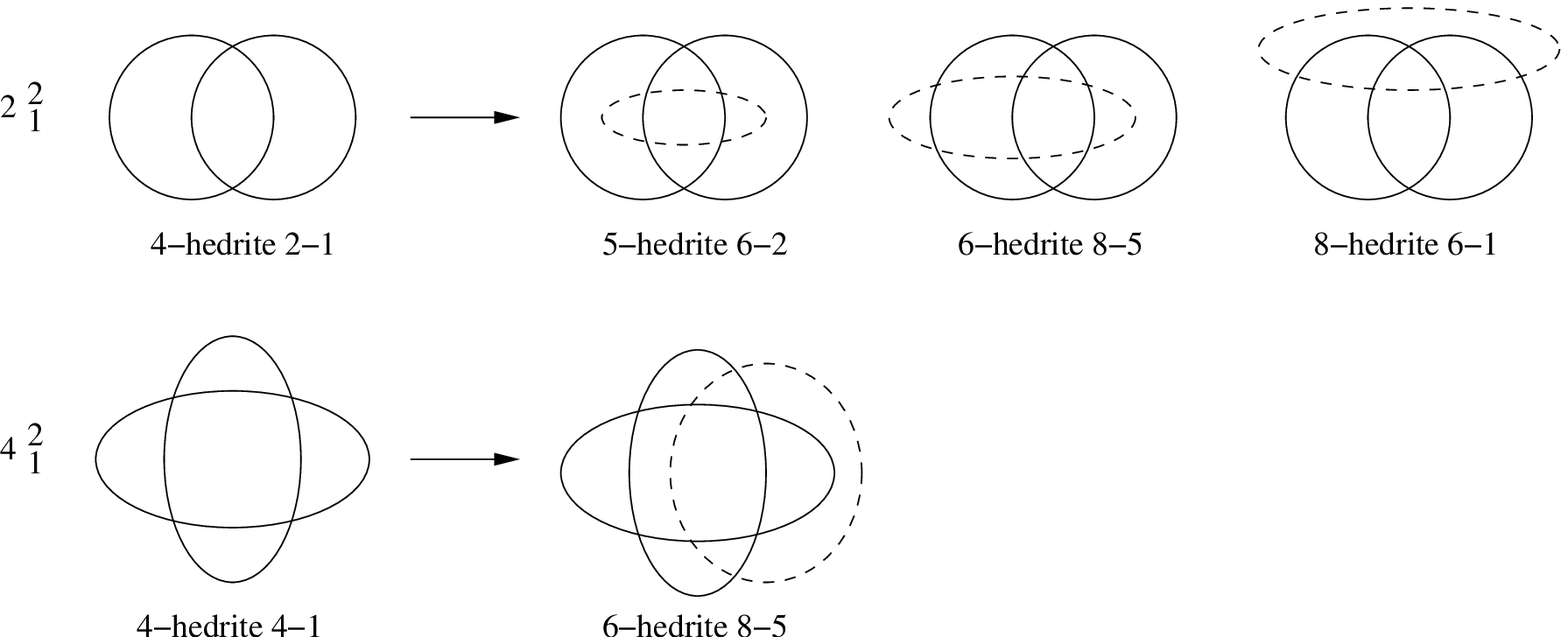}
\end{center}

Now, all irreducible pure $i$-hedrites with {\em three} central circuits are $5$-hedrite {\bf 6-2}, $6$-hedrite {\bf 8-5} and $8$-hedrite {\bf 6-1} (i.e. the projections of links $6^3_1$, $8^3_6$ and $6^3_2$).
In fact, we apply the same procedure to those three $i$-hedrites; see picture below:

\begin{center}
\epsfxsize=150mm
\epsfbox{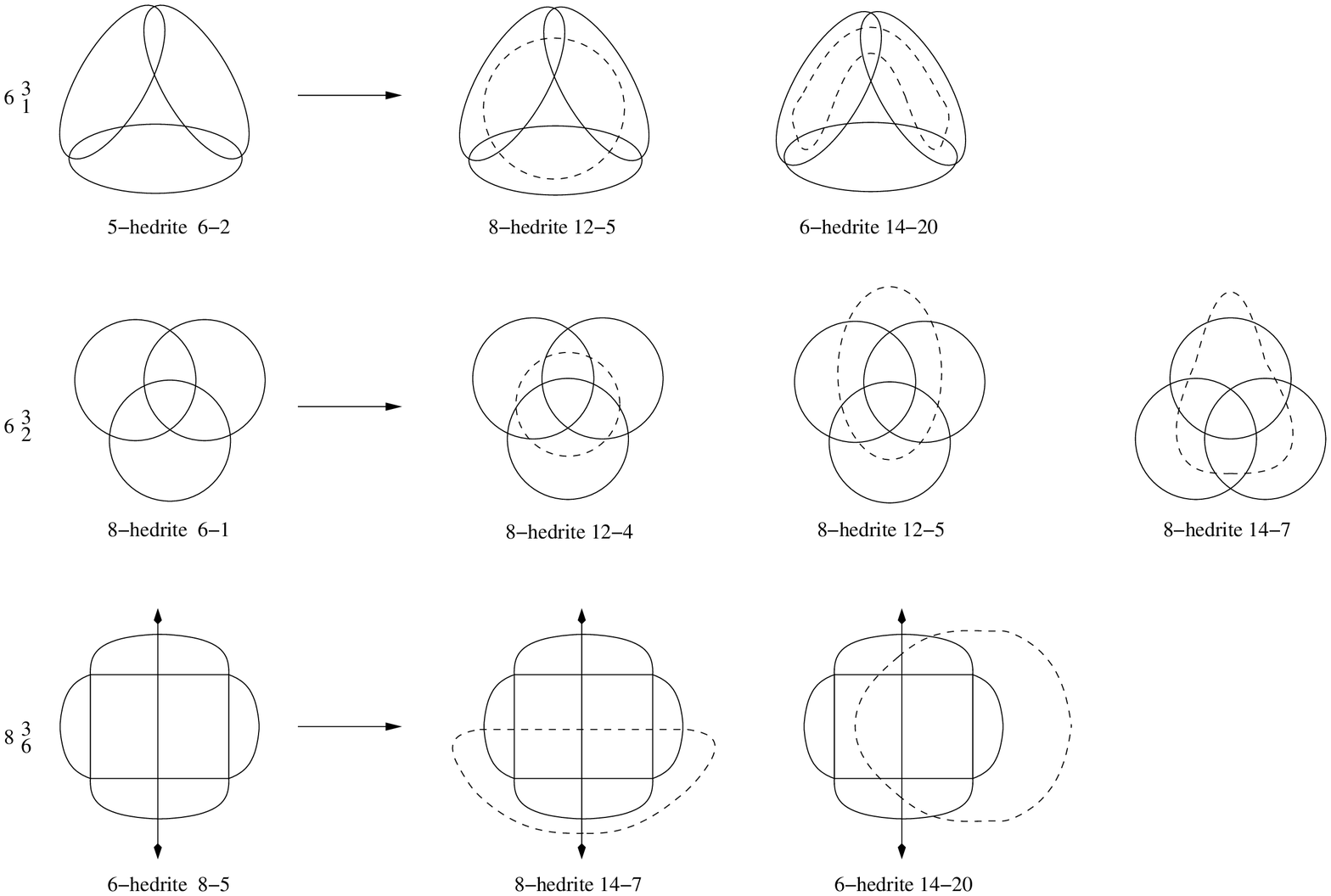}
\end{center}

Next, all irreducible pure $i$-hedrites with {\em four} central circuits are $8$-hedrites {\bf 12-4}, {\bf 12-5}, {\bf 14-7} and $6$-hedrite {\bf 14-20}.

By the same method, one can see that there are exactly two pure irreducible $i$-hedrites with {\em five} central circuits and two with {\em six} central circuits (see Figure \ref{ThePureIrreducibleOctahedriteWith56CC}).

\begin{remark}

Any pure $i$-hedrite comes from a pure {\em irreducible} $i$-hedrite 
with, say, $j$ central circuits by simultaneous $t_1-, \dots, t_j$-inflation along those circuits; it is irreducible if and only if $t_1=\dots=t_j=1$.

\end{remark}

\section{Symmetry groups of $i$-hedrites}
We consider below the maximal symmetry groups of plane graphs; these groups are identified with the corresponding point groups.

%
%
%

\begin{theorem}
We indicate here the lists of symmetry groups of $i$-hedrites,
together with the smallest number of vertices, for which they
appear:
\begin{itemize}
\item[(i)] The only symmetry groups of $4$-hedrites are point subgroups of $D_{4h}$, which contain $D_{2}$ as point subgroup, i.e. $D_{4h}$($n=2$), $D_4$($n=10$), $D_{2h}$($n=4$), $D_{2d}$($n=6$), $D_2$($n=12$).

\item[(ii)] The only symmetry groups of $5$-hedrites are: $C_1$($n=10$), $C_2$($n=8$), $C_s$($n=7$), $C_{2v}$($n=5$), $D_3$($n=15$), $D_{3h}$($n=3$).

\item[(iii)] The only symmetry groups of $6$-hedrites are: $D_{2d}$($n=4$), $D_{2h}$($n=6$) and all their point subgroups, i.e. $D_{2}$($n=12$), $C_{2h}$($n=10$), $C_{2v}$($n=5$), $C_i$($16\le n\le 30$), $C_{2}$($n=6$), $C_{s}$($n=9$), $C_{1}$($n=9$).

\item[(iv)] The only symmetry groups of $7$-hedrites are point subgroups of $C_{2v}$, i.e. $C_{2v}$($n=7$), $C_{2}$($n=11$), $C_{s}$($n=8$), $C_{1}$($n=11$).

\item[(v)] The only symmetry groups of $8$-hedrites are: $C_{1}$($n=16$), $C_s$($n=14$), $C_2$($n=12$), $C_{2v}$($n=11$), $C_i$($22\leq n\leq 46$), $C_{2h}$($22\leq n \leq 26$), $S_4$($22\leq n\leq 60$), $D_2$($n=10$), $D_{2d}$($n=14$), $D_{2h}$($n=22$), $D_3$($n=18$), $D_{3d}$($n=12$), $D_{3h}$($n=9$), $D_4$($n=18$), $D_{4d}$($n=8$), $D_{4h}$($n=10$), $O$($n=30$), $O_h$($n=6$).

\end{itemize}

\end{theorem}
\proof For $4$-hedrites, see \cite{DSt}. Any transformation
stabilizing a $2$-gon, can interchange its two edges and two vertices. So, the stabilizing point subgroup of a $2$-gon can be $C_{2v}$, $C_s$, $C_2$ or $C_1$ only.

The unique $2$-gon of a $7$-hedrite has to be preserved by the symmetry
group; so, all possibilities are: $C_{2v}$, $C_s$, $C_2$, $C_1$.

Every symmetry of an $i$-hedrite induces a symmetry on its $2$-gons and $3$-gons. Since the stabilizer of a $2$-gon, $3$-gon has maximal size $4$, $6$, this imply that the order of the symmetry group of an $i$-hedrite is bounded from above by $4|Sym(8-i)|=4(8-i)!$ and $6|Sym(2i-8)|=6(2i-8)!$.

So, in particular, the order of symmetry group of an $6$-hedrite is at
most $8$. If $f$ is an element of order three, then it fixes each of two
$2$-gons. Since the stabilizer of a $2$-gon does not contain an element
of order three, no such $f$ exists. If $f$ is a rotation of
order $4$, then $f^2$ is a rotation of order $2$ stabilizing each $2$-gon;
so, the axis of $f$ goes through the two $2$-gons. This is a contradiction.
By a search in the Tables of the groups one can see that the only possibilities are: $C_1$, $C_s$, $C_2$, $C_i$, $C_{2v}$, $C_{2h}$, $D_2$, $D_{2h}$, $D_{2d}$. But there exists a $6$-hedrite for any of such symmetries, in Figure \ref{special-i-hedrites} and subsection \ref{subsection-6-hedrites} below.

For $5$-hedrites, since there are two $3$-gons, the maximal order of the group is $12$. The oddness of the number of $2$-gons excludes central symmetry, axis of order $4$, and groups $D_2$, $D_{2h}$, $D_{2d}$.
If $G$ is a $5$-hedrite with a $3$-fold axis then this axis goes through
the two $3$-gons, say, $T_1$ and $T_2$. If one consider a belt of $4$-gons
around $T_1$, then, after a number $p$ of steps, one will encounter a
$2$-gon and so, by symmetry, three $2$-gons. So, we will have the 
following possibilities for $p=1$:

\begin{center}
\epsfxsize=90mm
\epsfbox{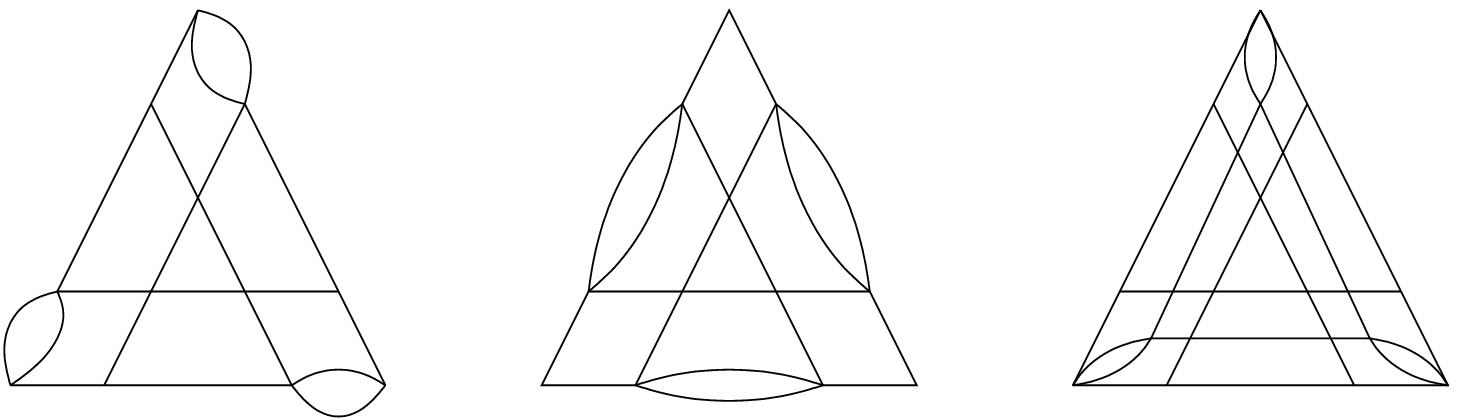}
\end{center}

There is only one way to extend this graph to a $5$-hedrite and the obtained extension has symmetry at least $D_3$. So, the group is $D_{3}$, $D_{3h}$ or $D_{3d}$.

An $8$-hedrite $G$ has $k$-fold axis of rotation with $k=2$, $3$ or $4$. 
If $k=3$, then the axis of the rotation goes through two $3$-gons,
say $T_1$ and $T_8$. If one consider, around $3$-gon $T_1$, a belt of 
$4$-gons, then, after some number $p$ of steps, one will encounter a 
$3$-gon and so, by symmetry, three $3$-gons, say, $\{T_2,T_3,T_4\}$. 
Adding, if necessary, belts of $4$-gons, one will encounter the last 
three $3$-gons, say, $\{T_5,T_6,T_7\}$. The patch formed by the six
triangles $T_2, \dots, T_7$ has symmetry $D_3$ at least and so, 
$G$ has also this symmetry. Consequently, the symmetry of $G$ is $D_3$,
$D_{3h}$ or $D_{3d}$.

If $k=4$, then the axis of the rotation goes through a vertex or
a $4$-gon. Assume, for simplicity, that this axis goes through a
vertex; then, by repeating above reasoning, one obtains two orbits of
$3$-gons, say, $\{T_1, T_2, T_3, T_4\}$ and $\{T_5, T_6, T_7, T_8\}$,
under $4$-fold symmetry and the symmetry group is $D_4$, $D_{4h}$,
$D_{4d}$, $O$ or $O_h$.

So, one obtains the above list of $18$ possible point groups. All 
these groups appear in subsection \ref{subsection-8-hedrites} (groups 
$O_h$, $D_{4d}$, $D_{3h}$, $D_{2}$, $D_{4h}$, $C_{2v}$, $D_{3d}$, $C_2$),
in Figure \ref{special-i-hedrites} (groups $D_{2d}$, $C_s$, $C_1$, $D_4$, 
$D_{3}$, $C_{2h}$, $C_{i}$, $S_4$) and in Figure 
\ref{ThePureIrreducibleOctahedriteWith56CC} (groups $D_{2h}$, $O$); for
all groups (except, possibly, $C_i$, $C_{2h}$ and $S_4$) those examples
of $8$-hedrites have smallest number of vertices.

\begin{figure}
{\small
\begin{minipage}[t]{4cm}
\centering
\epsfxsize=4cm
\epsffile{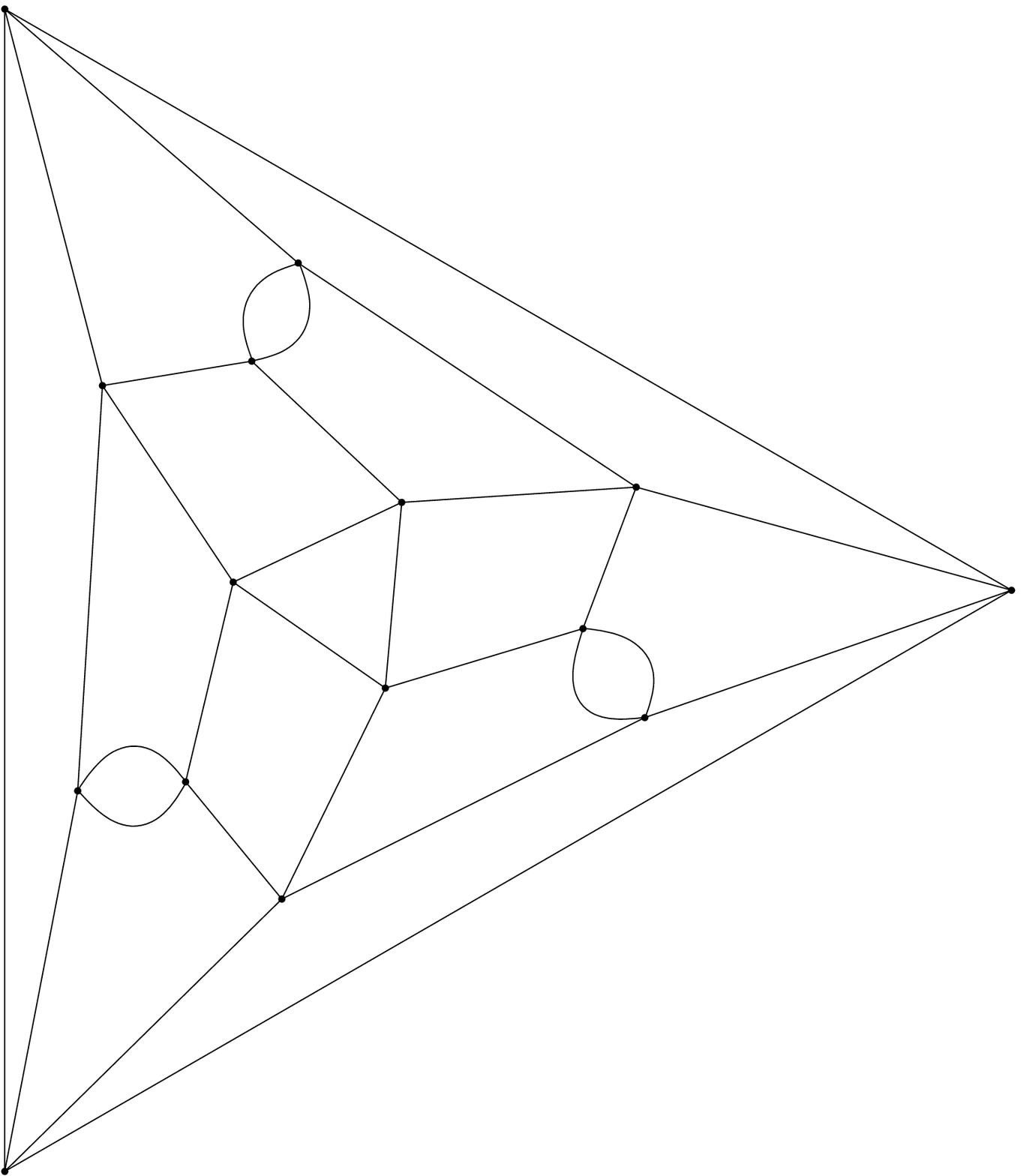}\par
{$5$-hedrite {\bf 15-1} $D_3$ smallest $(30)$}
\end{minipage}
\begin{minipage}[t]{4cm}
\centering
\epsfxsize=4cm
\epsffile{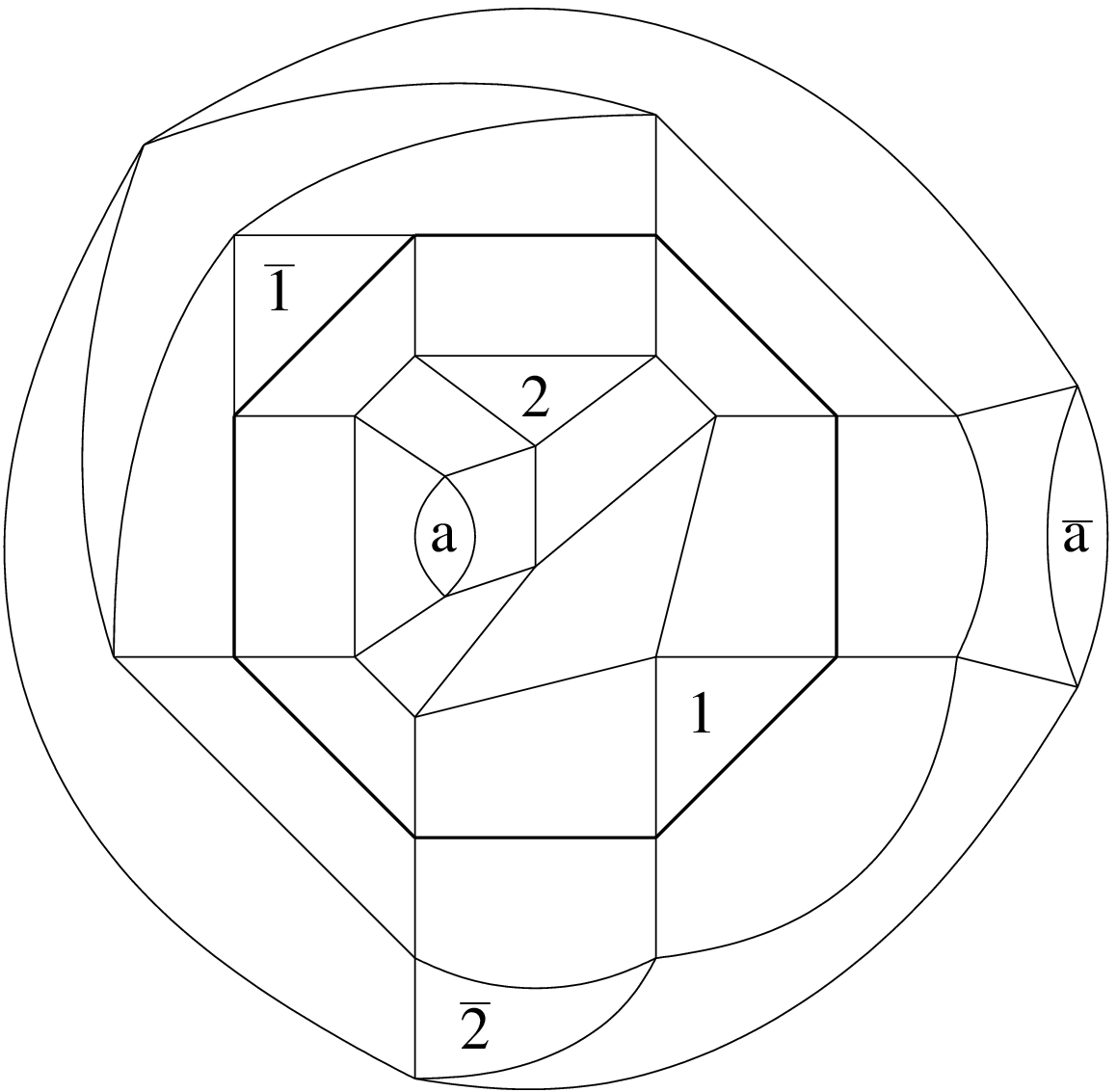}\par
{$6$-hedrite {\bf 30-1} $C_i$ $(8;26^2)$}
\end{minipage}
\begin{minipage}[t]{4cm}
\centering
\epsfxsize=4cm
\epsffile{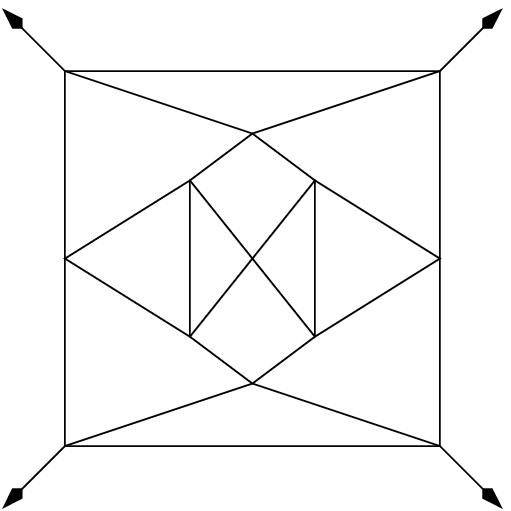}\par
{$8$-hedrite {\bf 14-4} $D_{2d}$ smallest $(14^2)$}
\end{minipage}
\begin{minipage}[t]{4cm}
\centering
\epsfxsize=4cm
\epsffile{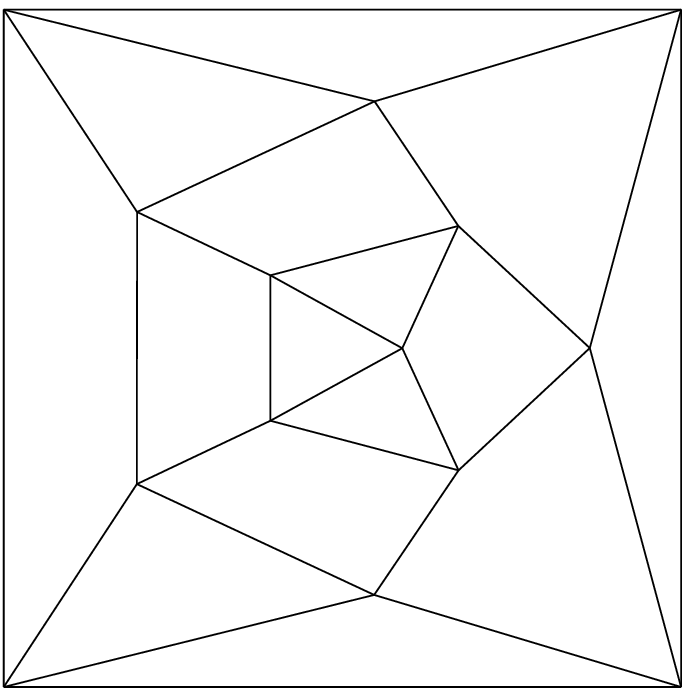}\par
{$8$-hedrite {\bf 14-2} $C_s$ smallest $(6;22)$}
\end{minipage}
\begin{minipage}[t]{4cm}
\centering
\epsfxsize=4cm
\epsffile{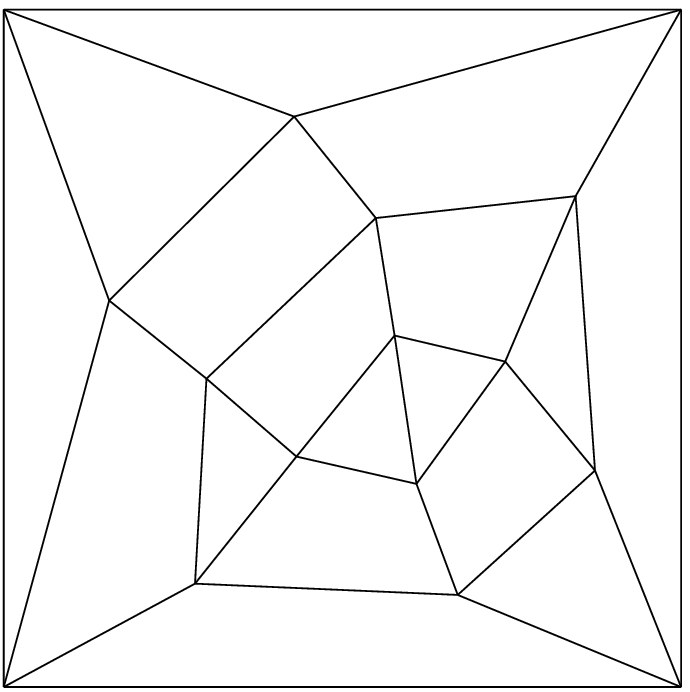}\par
{$8$-hedrite {\bf 16-1} $C_1$ smallest $(6, 8; 18)$}
\end{minipage}
\begin{minipage}[t]{4cm}
\centering
\epsfxsize=4cm
\epsffile{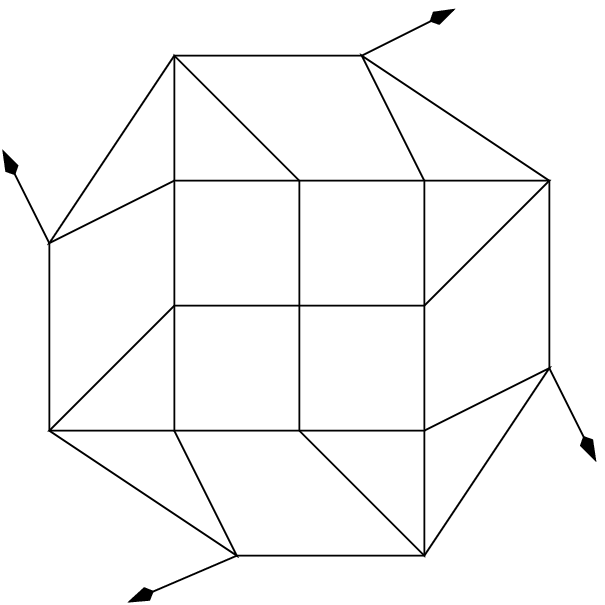}\par
{$8$-hedrite {\bf 18-1} $D_4$ smallest $(18^2)$}
\end{minipage}
\begin{minipage}[t]{4cm}
\centering
\epsfxsize=4cm
\epsffile{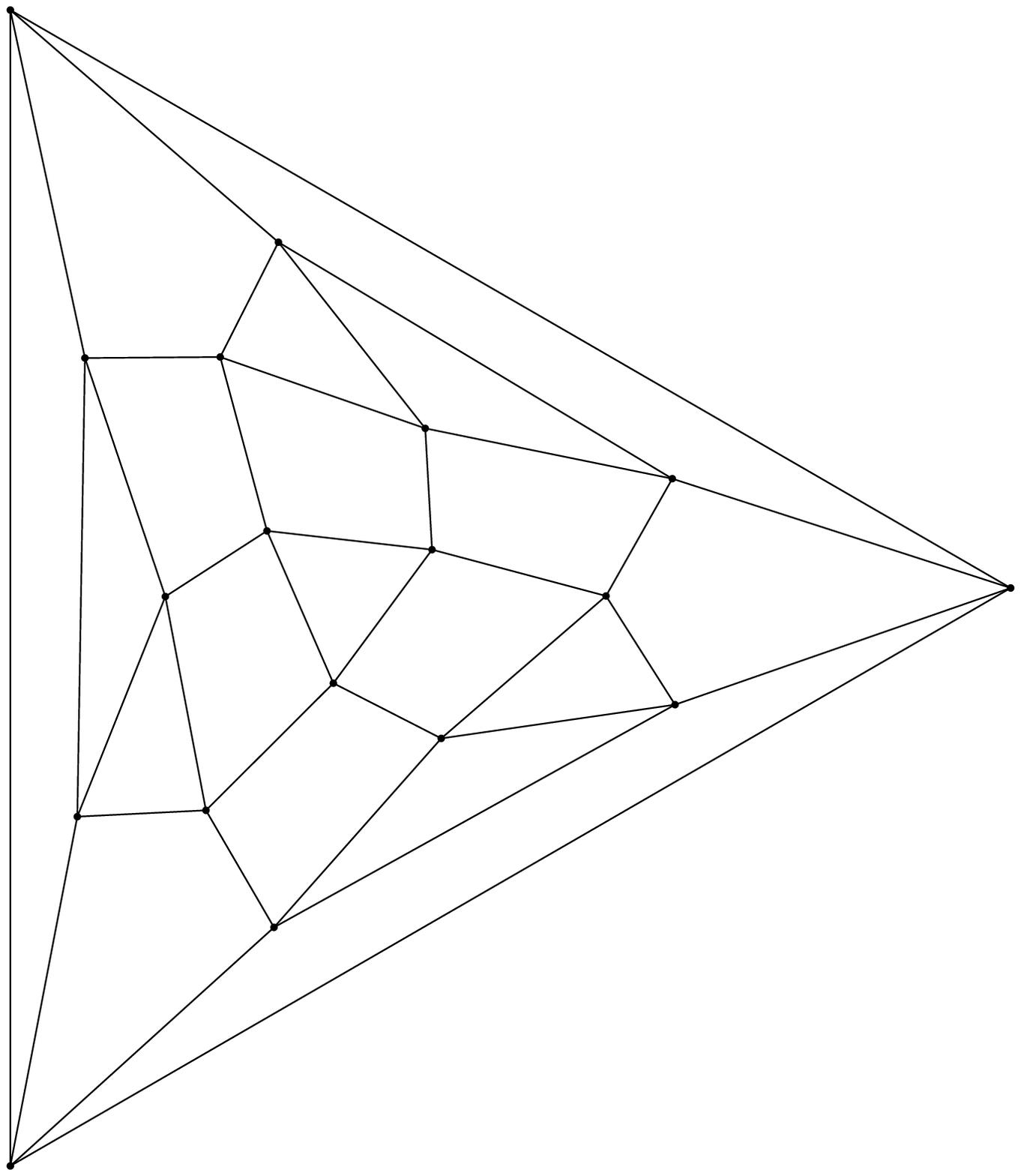}\par
{$8$-hedrite {\bf 18-2} $D_3$ smallest $(36)$}
\end{minipage}
\begin{minipage}[t]{4cm}
\centering
\epsfxsize=4cm
\epsffile{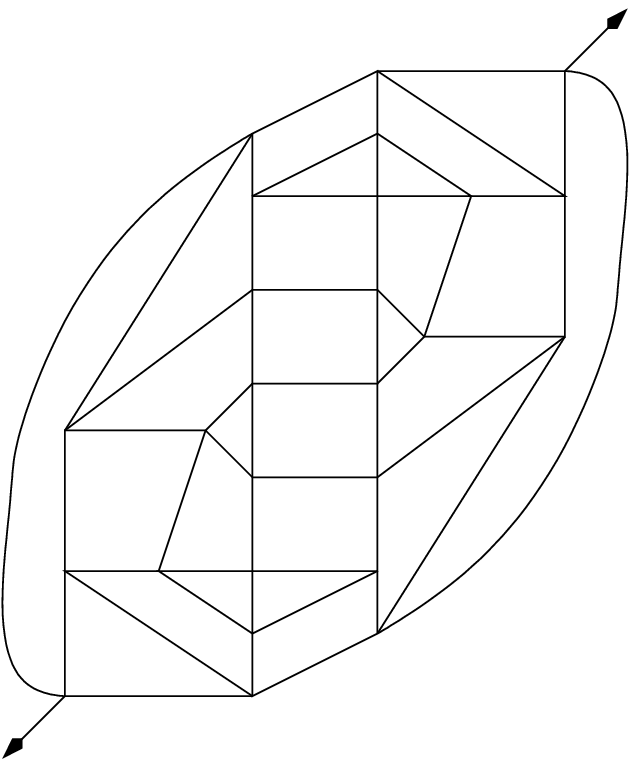}\par
{$8$-hedrite {\bf 26-1} $C_{2h}$ $(8^2; 36)$}
\end{minipage}
\setlength{\unitlength}{1cm}
\begin{minipage}[t]{4cm}
\centering
\epsfxsize=4cm
\epsffile{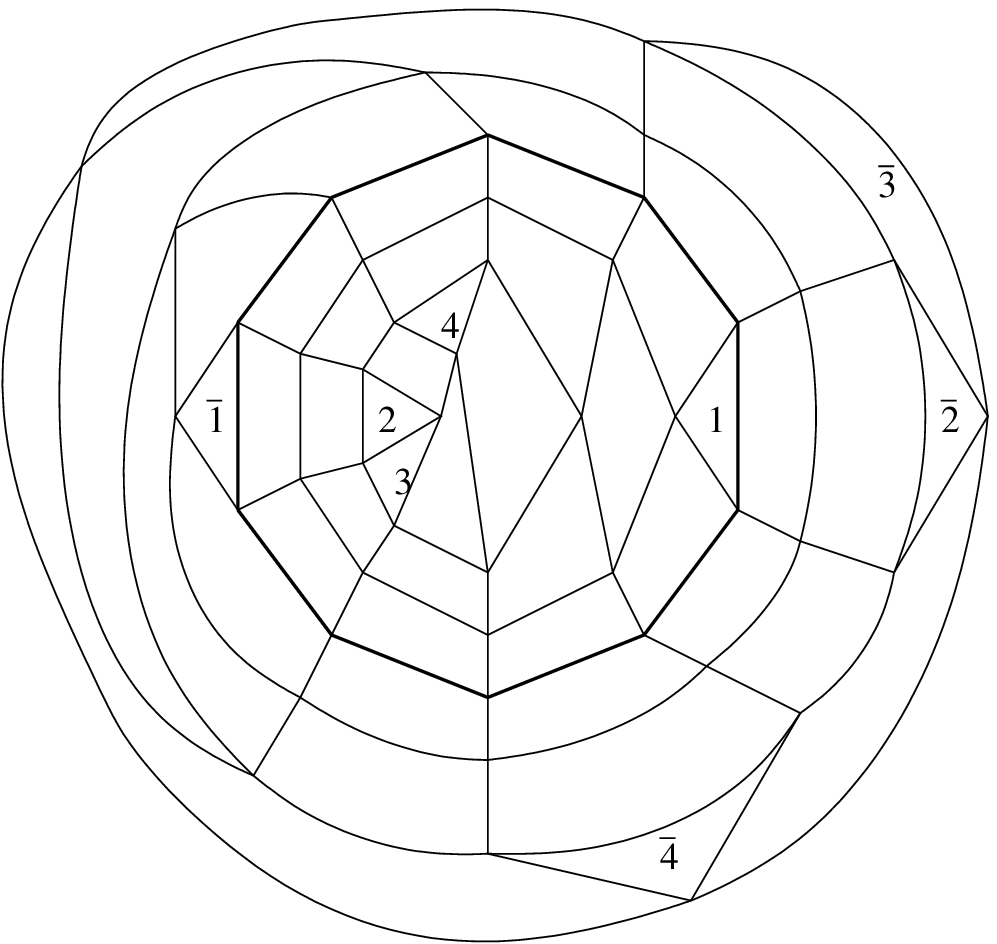}\par
{$8$-hedrite {\bf 46-1} $C_i$ $(10; 82)$}
\end{minipage}
\setlength{\unitlength}{1cm}
\begin{minipage}[t]{4cm}
\centering
\epsfxsize=4cm
\epsffile{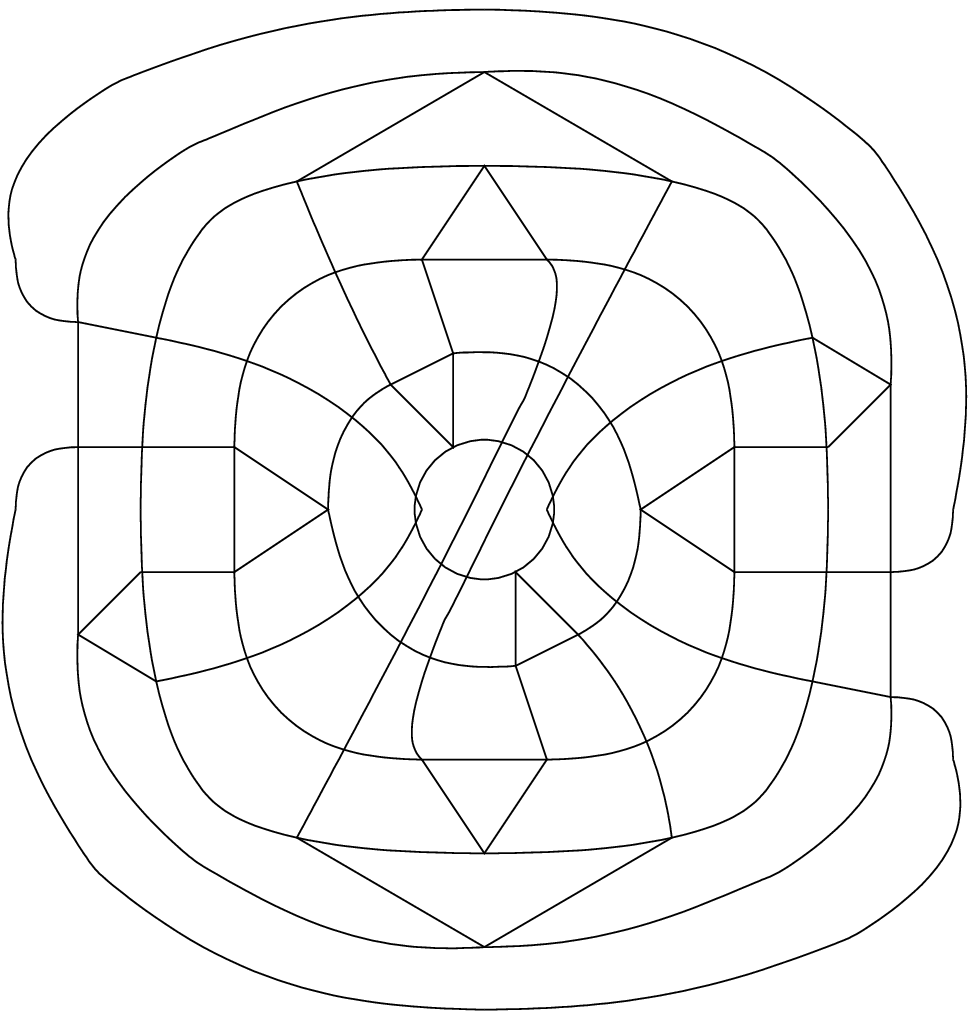}\par
{$8$-hedrite {\bf 60-1} $S_4$ $(16;26^4)$}
\end{minipage}
}
\caption{Some $i$-hedrites with special symmetry groups}
\label{special-i-hedrites}
\end{figure}


\begin{remark}
The simplest case, $i=4$, of $i$-hedrites admits following characterization for each of
its five possible groups.
It has symmetry $D_{2h}$  (respectively, $D_{2d}$) if and only if it is an
$t$-inflation, for some $t \ge 1, m \ge 2$, of $J_{4,2m}$ or $I_{4,2m+2}$
with even $m$ (respectively, of $K_{4,4m}$ or $I_{4,2m+2}$ with odd $m$).
Any $4$-hedrite with group $D_{4}$ or $D_{4h}$ has $2(k^2+l^2)$ vertices for
some $k \ge l \ge 0$ (the group $D_{4h}$ corresponds to the case $l=k$ or $0$);
it comes from the smallest $4$-hedrite {\bf 2-1} by Goldberg-Coxeter
construction (see \cite{Gold37}, \cite{Cox71} and \cite{DD03}).
All other 4-hedrites have symmetry $D_2$; in shift terms, they are exactly those, for which interchange of central circuits changes the value of shift.
\end{remark}

We expect, that a $7$-hedrite with the highest symmetry $C_{2v}$ 
exists for any $n\geq 10$.


\section{Small $i$-hedrites}

Here and below all links are given 
in Rolfsen's notation (see the table in \cite{Rolf} and also,  
for example, \cite{Kaw}) for links with at most 9 
crossings and knots with 10 crossings, or, otherwise, in
Dowker-Thistlewhaite's numbering (see \cite{T}), if any.
We write $\sim$ if the projection in the pictures and Table below 
is different from the one given in corresponding cases above.

We give on the pictures below all $i$-hedrites with at most $12$ vertices,
indicating under 
picture of each its symmetry, CC-vector and corresponding alternating link.
If an $i$-hedrite is $2$-connected but not $3$-connected, then we add
a symbol ${\bf *}$ just after the number. If an $i$-hedrite is reducible,
i.e. has a rail-road, then we add mention ``red.''. All $i$-hedrites
with $13$, $14$ and $15$ vertices are listed in Table
\ref{tab:i-hedrite13_14}.

Only three reducible $i$-hedrites with $n \leq 15$ have self-intersecting
railroad: $5$-hedrites {\bf 12-3}, {\bf 14-6} and $6$-hedrite {\bf 13-11}.

On the pictures below, in order to express better the (maximal)
symmetry of an $i$-hedrite, we put:

(i) a double arrow, in order to represent an edge passing at infinity,

(ii) a quadruple arrow, in order to represent a vertex at infinity.

\subsection{$4$-hedrites}\label{subsection-4-hedrites}
{\small
\setlength{\unitlength}{1cm}
\begin{minipage}[t]{3.5cm}
\centering
\epsfxsize=2.5cm
\epsffile{4-hedrite2_1.eps}\par
{{\bf Nr.2-1} \quad $D_{4h}$\\ $2^2_1$ \quad $(2^2)$\\\vspace{3mm} }
\end{minipage}
\setlength{\unitlength}{1cm}
\begin{minipage}[t]{3.5cm}
\centering
\epsfxsize=2.5cm
\epsffile{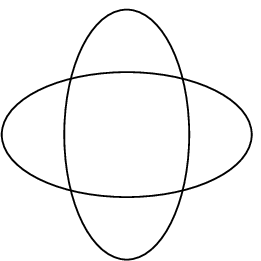}\par
{{\bf Nr.4-1} \quad $D_{4h}$\\ $4^2_1$ \quad $(4^2)$\\\vspace{3mm} }
\end{minipage}
\setlength{\unitlength}{1cm}
\begin{minipage}[t]{3.5cm}
\centering
\epsfxsize=2.5cm
\epsffile{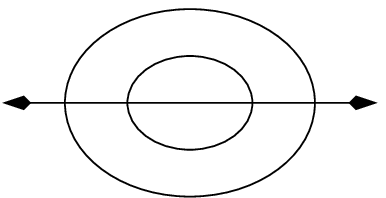}\par
{{\bf Nr.4-2${}^*$} \quad $D_{2h}$\\ $2\times 2^2_1$ \quad $(2^2,4)$ red.\\\vspace{3mm} }
\end{minipage}
\setlength{\unitlength}{1cm}
\begin{minipage}[t]{3.5cm}
\centering
\epsfxsize=2.5cm
\epsffile{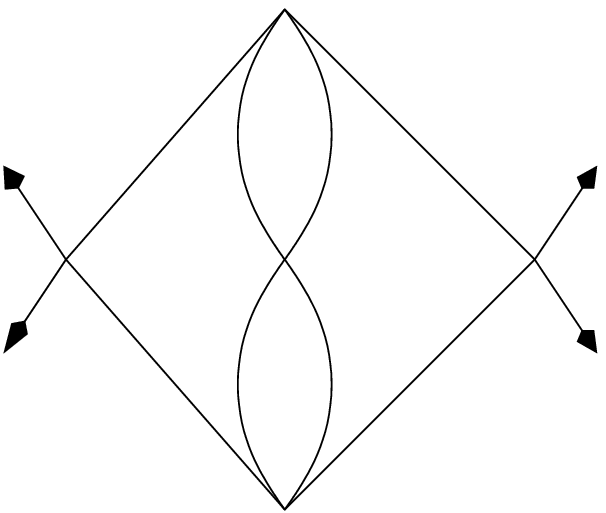}\par
{{\bf Nr.6-1${}^*$} \quad $D_{2d}$\\ $6^2_2$ \quad $(6^2)$\\\vspace{3mm} }
\end{minipage}
\setlength{\unitlength}{1cm}
\begin{minipage}[t]{3.5cm}
\centering
\epsfxsize=2.5cm
\epsffile{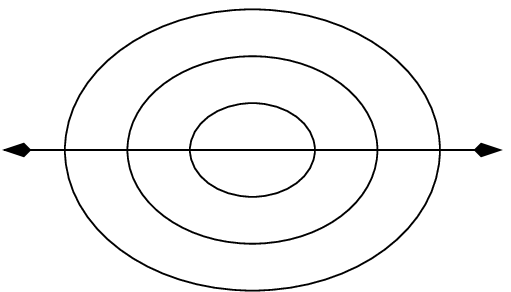}\par
{{\bf Nr.6-2${}^*$} \quad $D_{2h}$\\ $3\times 2^2_1$ \quad $(2^3,6)$ red.\\\vspace{3mm} }
\end{minipage}
\setlength{\unitlength}{1cm}
\begin{minipage}[t]{3.5cm}
\centering
\epsfxsize=2.5cm
\epsffile{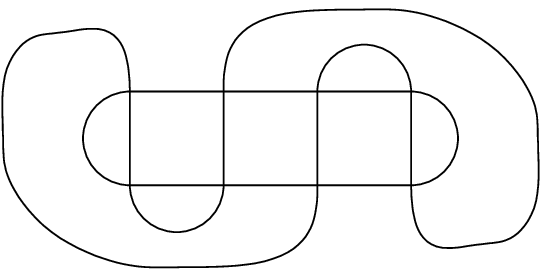}\par
{{\bf Nr.8-1${}^*$} \quad $D_{2h}$\\ $\sim 8^2_4$ \quad $(8^2)$\\\vspace{3mm} }
\end{minipage}
\setlength{\unitlength}{1cm}
\begin{minipage}[t]{3.5cm}
\centering
\epsfxsize=2.5cm
\epsffile{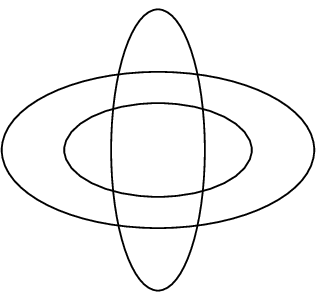}\par
{{\bf Nr.8-2} \quad $D_{2d}$\\ $8^3_4$ \quad $(4^2,8)$ red.\\\vspace{3mm} }
\end{minipage}
\setlength{\unitlength}{1cm}
\begin{minipage}[t]{3.5cm}
\centering
\epsfxsize=2.5cm
\epsffile{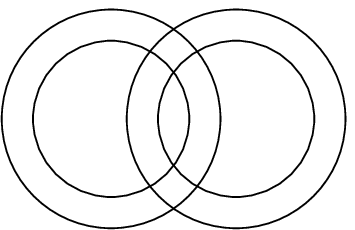}\par
{{\bf Nr.8-3} \quad $D_{4h}$\\ $8^4_{1}$ \quad $(4^4)$ red.\\\vspace{3mm} }
\end{minipage}
\setlength{\unitlength}{1cm}
\begin{minipage}[t]{3.5cm}
\centering
\epsfxsize=2.5cm
\epsffile{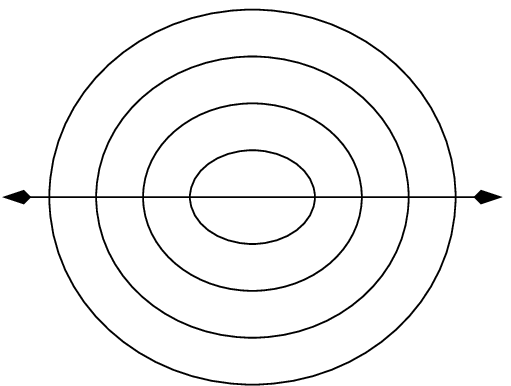}\par
{{\bf Nr.8-4${}^*$} \quad $D_{2h}$\\ $4\times 2^2_1$ \quad $(2^4,8)$ red.\\\vspace{3mm} }
\end{minipage}
\setlength{\unitlength}{1cm}
\begin{minipage}[t]{3.5cm}
\centering
\epsfxsize=2.5cm
\epsffile{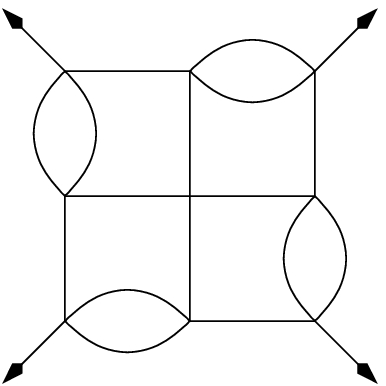}\par
{{\bf Nr.10-1} \quad $D_4$\\ $10^2_{121}$ \quad $(10^2)$\\\vspace{3mm} }
\end{minipage}
\setlength{\unitlength}{1cm}
\begin{minipage}[t]{3.5cm}
\centering
\epsfxsize=2.5cm
\epsffile{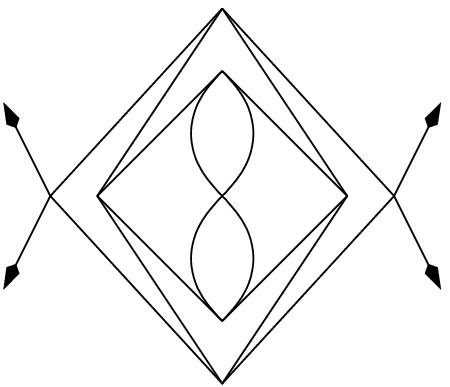}\par
{{\bf Nr.10-2${}^*$} \quad $D_{2d}$\\ $\sim 10^2_{120}$ \quad $(10^2)$\\\vspace{3mm} }
\end{minipage}
\setlength{\unitlength}{1cm}
\begin{minipage}[t]{3.5cm}
\centering
\epsfxsize=2.5cm
\epsffile{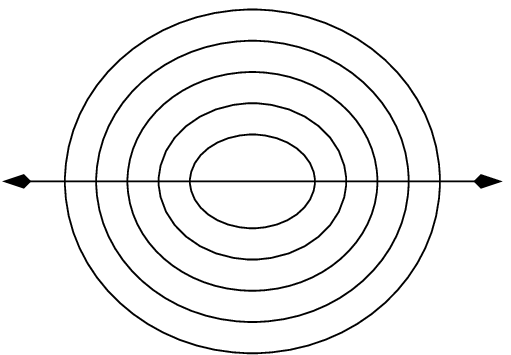}\par
{{\bf Nr.10-3${}^*$} \quad $D_{2h}$\\ $5\times 2^2_1$ \quad $(2^5,10)$ red.\\\vspace{3mm} }
\end{minipage}
\setlength{\unitlength}{1cm}
\begin{minipage}[t]{3.5cm}
\centering
\epsfxsize=2.5cm
\epsffile{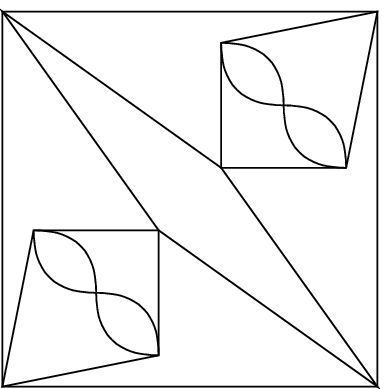}\par
{{\bf Nr.12-1${}^*$} \quad $D_{2h}$\\ $????$ \quad $(12^2)$\\\vspace{3mm} }
\end{minipage}
\setlength{\unitlength}{1cm}
\begin{minipage}[t]{3.5cm}
\centering
\epsfxsize=2.5cm
\epsffile{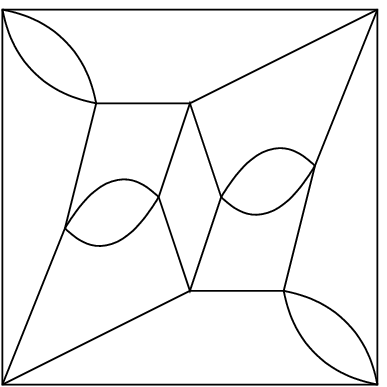}\par
{{\bf Nr.12-2} \quad $D_2$\\ $????$ \quad $(6^2,12)$ red.\\\vspace{3mm} }
\end{minipage}
\setlength{\unitlength}{1cm}
\begin{minipage}[t]{3.5cm}
\centering
\epsfxsize=2.5cm
\epsffile{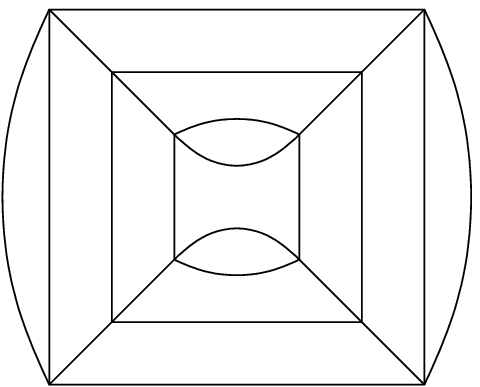}\par
{{\bf Nr.12-3} \quad $D_{2d}$\\ $????$ \quad $(4^3,12)$ red.\\\vspace{3mm} }
\end{minipage}
\setlength{\unitlength}{1cm}
\begin{minipage}[t]{3.5cm}
\centering
\epsfxsize=2.5cm
\epsffile{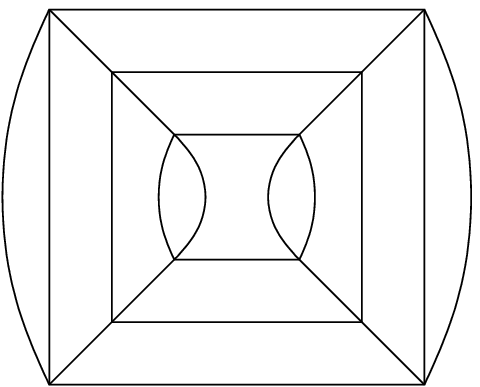}\par
{{\bf Nr.12-4} \quad $D_{2h}$\\ $????$ \quad $(4^3,6^2)$ red.\\\vspace{3mm} }
\end{minipage}
\setlength{\unitlength}{1cm}
\begin{minipage}[t]{3.5cm}
\centering
\epsfxsize=2.5cm
\epsffile{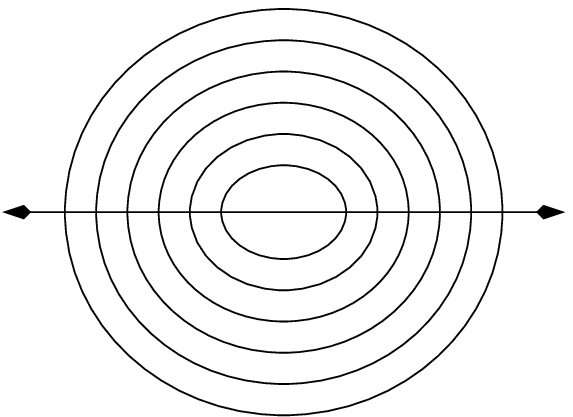}\par
{{\bf Nr.12-5${}^*$} \quad $D_{2h}$\\ $6\times 2^2_1$ \quad $(2^6,12)$ red.\\\vspace{3mm} }
\end{minipage}
}

\subsection{$5$-hedrites}\label{subsection-5-hedrites}
{\small
\setlength{\unitlength}{1cm}
\begin{minipage}[t]{3.5cm}
\centering
\epsfxsize=2.5cm
\epsffile{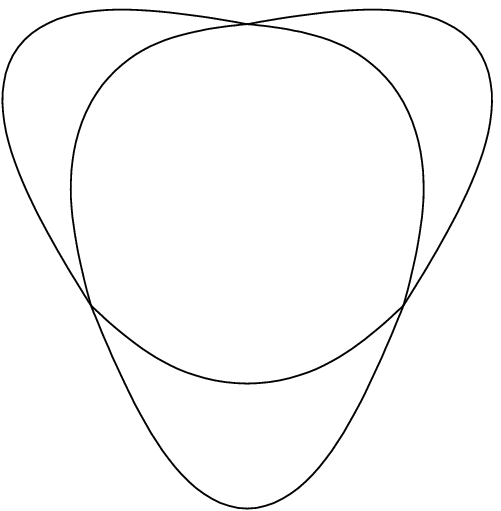}\par
{{\bf Nr.3-1} \quad $D_{3h}$\\ $3_{1}$ \quad $(6)$\\\vspace{3mm} }
\end{minipage}
\setlength{\unitlength}{1cm}
\begin{minipage}[t]{3.5cm}
\centering
\epsfxsize=2.5cm
\epsffile{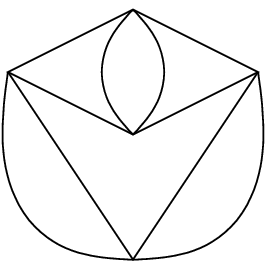}\par
{{\bf Nr.5-1${}^*$} \quad $C_{2v}$\\ $5_{2}$ \quad $(10)$\\\vspace{3mm} }
\end{minipage}
\setlength{\unitlength}{1cm}
\begin{minipage}[t]{3.5cm}
\centering
\epsfxsize=2.5cm
\epsffile{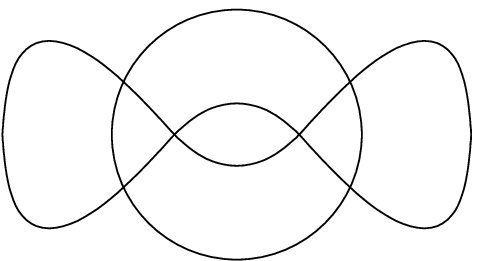}\par
{{\bf Nr.6-1} \quad $C_{2v}$\\ $6^2_{3}$ \quad $(4;8)$\\\vspace{3mm} }
\end{minipage}
\setlength{\unitlength}{1cm}
\begin{minipage}[t]{3.5cm}
\centering
\epsfxsize=2.5cm
\epsffile{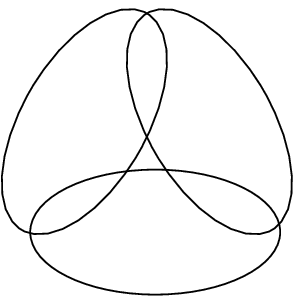}\par
{{\bf Nr.6-2} \quad $D_{3h}$\\ $6^3_{1}$ \quad $(4^3)$\\\vspace{3mm} }
\end{minipage}
\setlength{\unitlength}{1cm}
\begin{minipage}[t]{3.5cm}
\centering
\epsfxsize=2.5cm
\epsffile{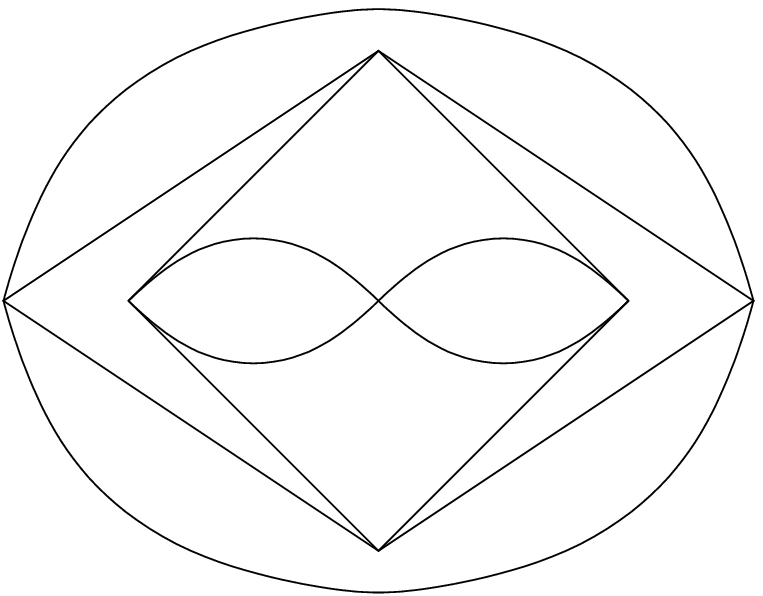}\par
{{\bf Nr.7-1${}^*$} \quad $C_{2v}$\\ $\sim 7_{5}$ \quad $(14)$\\\vspace{3mm} }
\end{minipage}
\setlength{\unitlength}{1cm}
\begin{minipage}[t]{3.5cm}
\centering
\epsfxsize=2.5cm
\epsffile{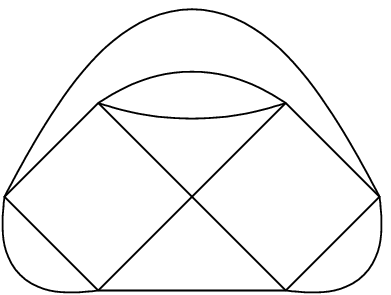}\par
{{\bf Nr.7-2} \quad $C_{s}$\\ $7^2_{5}$ \quad $(4;10)$\\\vspace{3mm} }
\end{minipage}
\setlength{\unitlength}{1cm}
\begin{minipage}[t]{3.5cm}
\centering
\epsfxsize=2.5cm
\epsffile{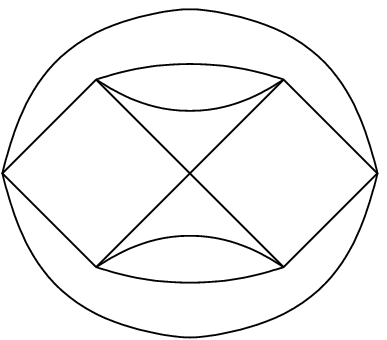}\par
{{\bf Nr.7-3} \quad $C_{2v}$\\ $7^3_1$ \quad $(4^2;6)$\\\vspace{3mm} }
\end{minipage}
\setlength{\unitlength}{1cm}
\begin{minipage}[t]{3.5cm}
\centering
\epsfxsize=2.5cm
\epsffile{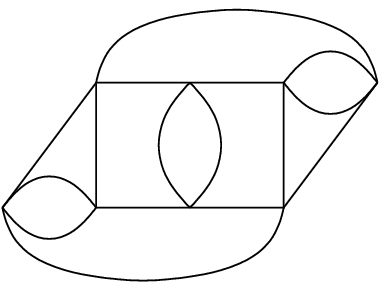}\par
{{\bf Nr.8-1} \quad $C_{2}$\\ $\sim 8_{15}$ \quad $(16)$\\\vspace{3mm} }
\end{minipage}
\setlength{\unitlength}{1cm}
\begin{minipage}[t]{3.5cm}
\centering
\epsfxsize=2.5cm
\epsffile{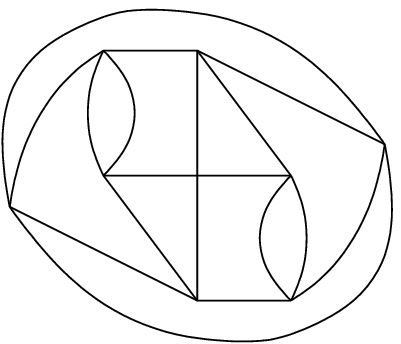}\par
{{\bf Nr.9-1} \quad $C_{2}$\\ $9_{38}$ \quad $(18)$\\\vspace{3mm} }
\end{minipage}
\setlength{\unitlength}{1cm}
\begin{minipage}[t]{3.5cm}
\centering
\epsfxsize=2.5cm
\epsffile{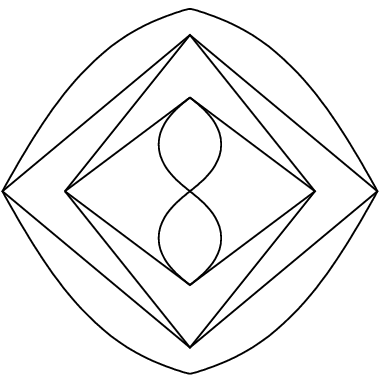}\par
{{\bf Nr.9-2${}^*$} \quad $C_{2v}$\\ $\sim 9_{18}$ \quad $(18)$\\\vspace{3mm} }
\end{minipage}
\setlength{\unitlength}{1cm}
\begin{minipage}[t]{3.5cm}
\centering
\epsfxsize=2.5cm
\epsffile{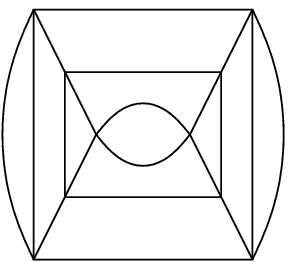}\par
{{\bf Nr.10-1} \quad $C_{2v}$\\ $10^3_{155}$ \quad $(4^2;12)$ red.\\\vspace{3mm} }
\end{minipage}
\setlength{\unitlength}{1cm}
\begin{minipage}[t]{3.5cm}
\centering
\epsfxsize=2.5cm
\epsffile{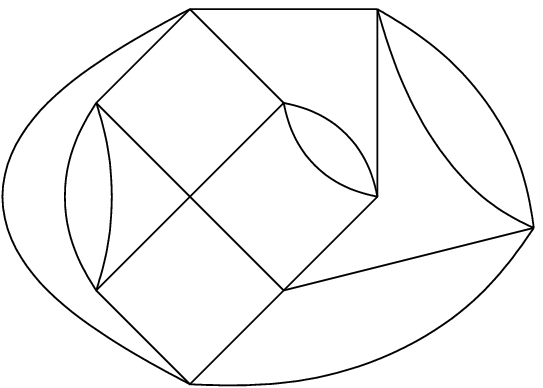}\par
{{\bf Nr.10-2} \quad $C_{1}$\\ $10^2_{85}$ \quad $(6;14)$\\\vspace{3mm} }
\end{minipage}
\setlength{\unitlength}{1cm}
\begin{minipage}[t]{3.5cm}
\centering
\epsfxsize=2.5cm
\epsffile{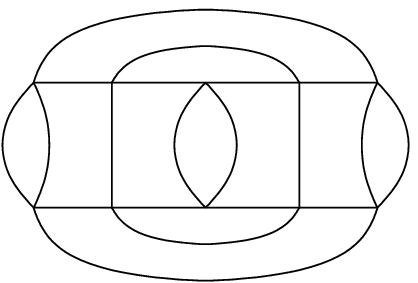}\par
{{\bf Nr.10-3} \quad $C_{2v}$\\ $10^4_{173}$ \quad $(4^2,6^2)$ red.\\\vspace{3mm} }
\end{minipage}
\setlength{\unitlength}{1cm}
\begin{minipage}[t]{3.5cm}
\centering
\epsfxsize=2.5cm
\epsffile{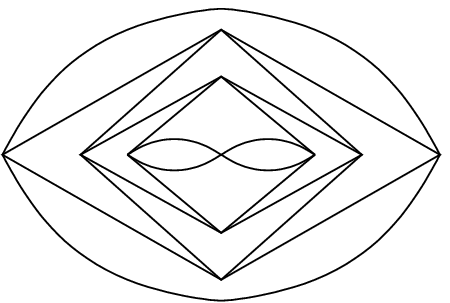}\par
{{\bf Nr.11-1${}^*$} \quad $C_{2v}$\\ $\sim 11_{236}$ \quad $(22)$\\\vspace{3mm} }
\end{minipage}
\setlength{\unitlength}{1cm}
\begin{minipage}[t]{3.5cm}
\centering
\epsfxsize=2.5cm
\epsffile{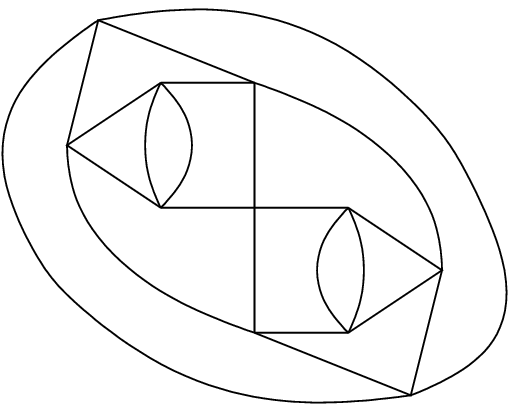}\par
{{\bf Nr.11-2} \quad $C_{2}$\\ $\sim 11_{124}$ \quad $(22)$\\\vspace{3mm} }
\end{minipage}
\setlength{\unitlength}{1cm}
\begin{minipage}[t]{3.5cm}
\centering
\epsfxsize=2.5cm
\epsffile{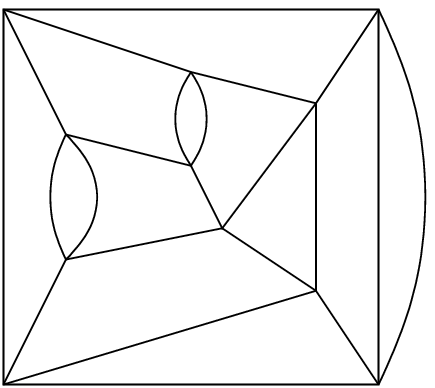}\par
{{\bf Nr.11-3} \quad $C_{1}$\\ $11^2_{226}$ \quad $(6;16)$\\\vspace{3mm} }
\end{minipage}
\setlength{\unitlength}{1cm}
\begin{minipage}[t]{3.5cm}
\centering
\epsfxsize=2.5cm
\epsffile{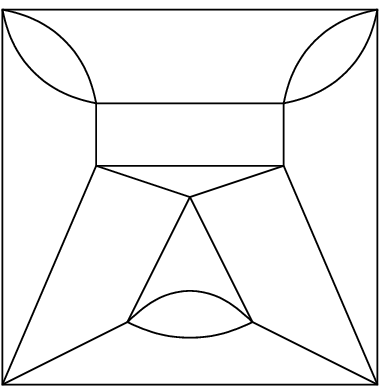}\par
{{\bf Nr.11-4} \quad $C_{s}$\\ $11^3_{500}$ \quad $(4^2;14)$ red.\\\vspace{3mm} }
\end{minipage}
\setlength{\unitlength}{1cm}
\begin{minipage}[t]{3.5cm}
\centering
\epsfxsize=2.5cm
\epsffile{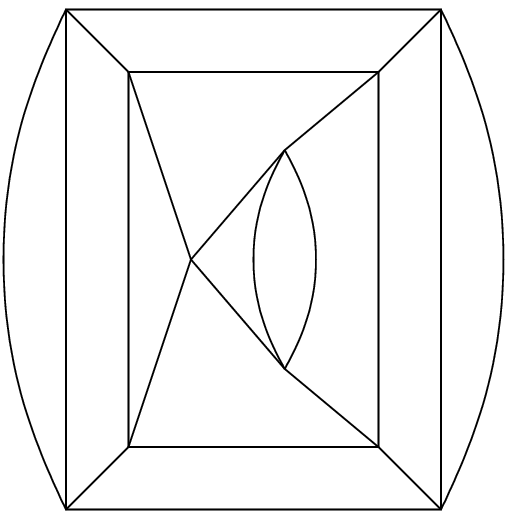}\par
{{\bf Nr.11-5} \quad $C_{s}$\\ $11^4_{547}$ \quad $(4^2,6;8)$ red.\\\vspace{3mm} }
\end{minipage}
\setlength{\unitlength}{1cm}
\begin{minipage}[t]{3.5cm}
\centering
\epsfxsize=2.5cm
\epsffile{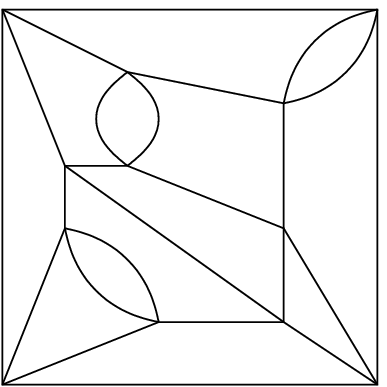}\par
{{\bf Nr.12-1} \quad $C_{1}$\\  $\sim 12_{431}$ \quad $(24)$\\\vspace{3mm} }
\end{minipage}
\setlength{\unitlength}{1cm}
\begin{minipage}[t]{3.5cm}
\centering
\epsfxsize=2.5cm
\epsffile{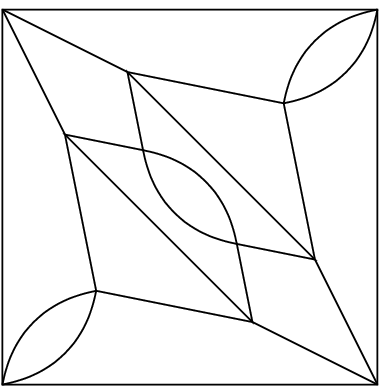}\par
{{\bf Nr.12-2} \quad $C_{2v}$\\ $????$ \quad $(12^2)$\\\vspace{3mm} }
\end{minipage}
\setlength{\unitlength}{1cm}
\begin{minipage}[t]{3.5cm}
\centering
\epsfxsize=2.5cm
\epsffile{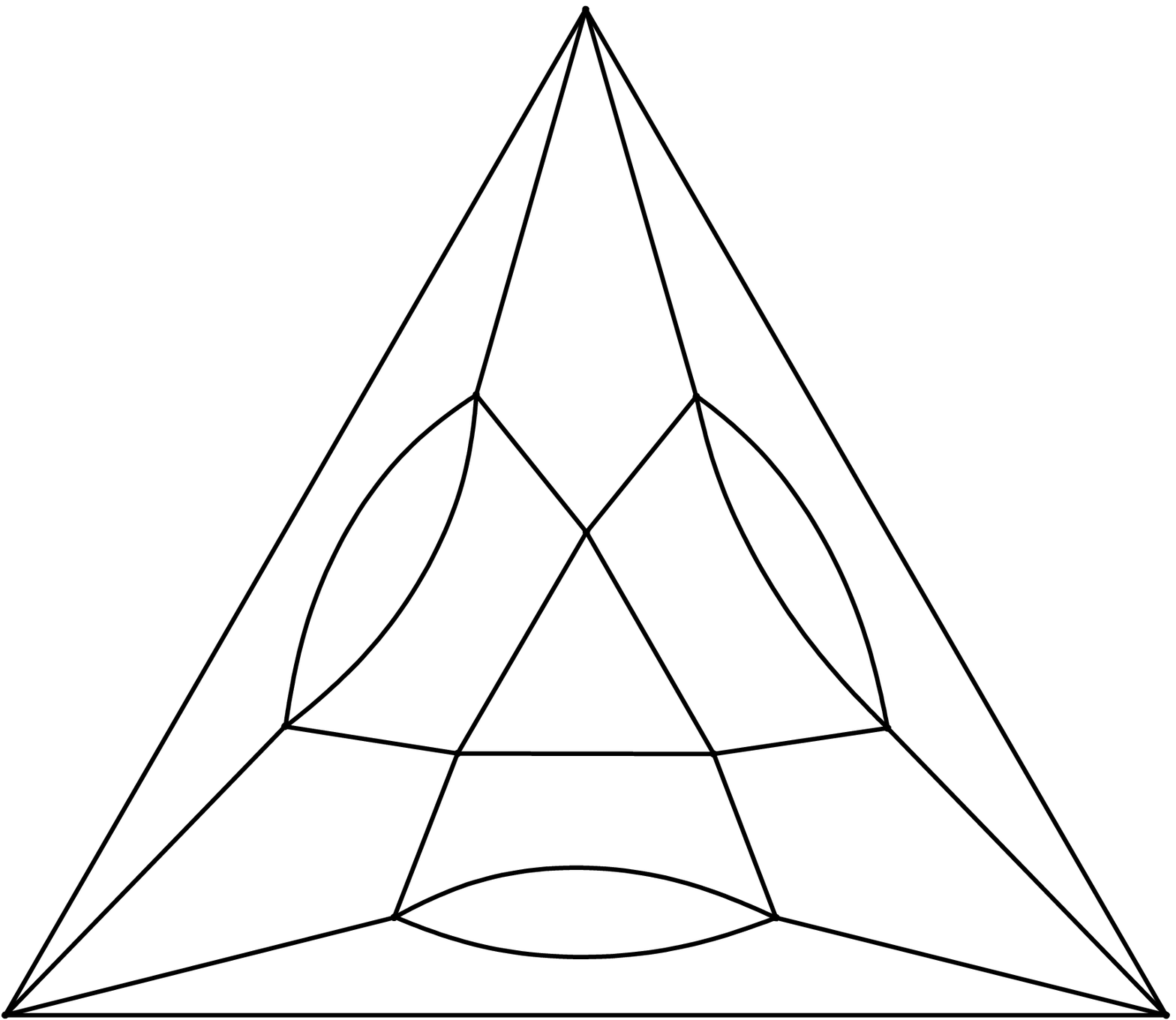}\par
{{\bf Nr.12-3} \quad $D_{3h}$\\ $????$ \quad $(12^2)$ red.\\\vspace{3mm} }
\end{minipage}
}

\subsection{$6$-hedrites}\label{subsection-6-hedrites}
{\small
\setlength{\unitlength}{1cm}
\begin{minipage}[t]{3.5cm}
\centering
\epsfxsize=2.5cm
\epsffile{6-hedrite4_1.eps}\par
{{\bf Nr.4-1} \quad $D_{2d}$\\ $4_{1}$ \quad $(8)$\\\vspace{3mm} }
\end{minipage}
\setlength{\unitlength}{1cm}
\begin{minipage}[t]{3.5cm}
\centering
\epsfxsize=2.5cm
\epsffile{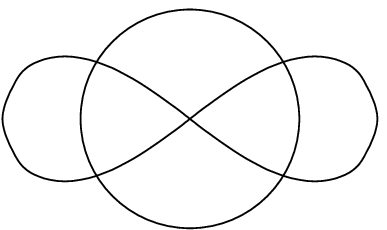}\par
{{\bf Nr.5-1} \quad $C_{2v}$\\ $5^2_{1}$ \quad $(4;6)$\\\vspace{3mm} }
\end{minipage}
\setlength{\unitlength}{1cm}
\begin{minipage}[t]{3.5cm}
\centering
\epsfxsize=2.5cm
\epsffile{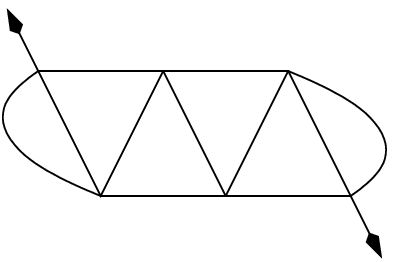}\par
{{\bf Nr.6-1} \quad $C_{2}$\\ $6_{3}$ \quad $(12)$\\\vspace{3mm} }
\end{minipage}
\setlength{\unitlength}{1cm}
\begin{minipage}[t]{3.5cm}
\centering
\epsfxsize=2.5cm
\epsffile{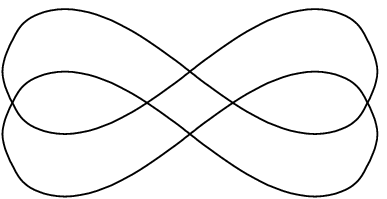}\par
{{\bf Nr.6-2${}^*$} \quad $D_{2h}$\\ $\sim 6^2_{3}$ \quad $(6^2)$\\\vspace{3mm} }
\end{minipage}
\setlength{\unitlength}{1cm}
\begin{minipage}[t]{3.5cm}
\centering
\epsfxsize=2.5cm
\epsffile{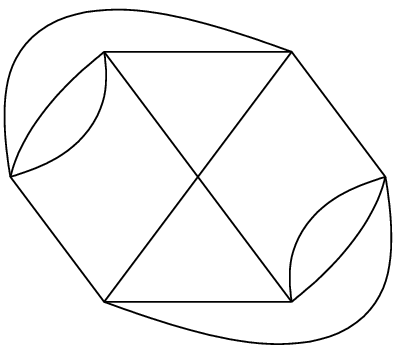}\par
{{\bf Nr.7-1} \quad $C_{2}$\\ $\sim 7_{7}$ \quad $(14)$\\\vspace{3mm} }
\end{minipage}
\setlength{\unitlength}{1cm}
\begin{minipage}[t]{3.5cm}
\centering
\epsfxsize=2.5cm
\epsffile{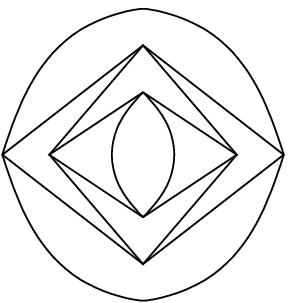}\par
{{\bf Nr.8-1${}^*$} \quad $D_{2d}$\\ $\sim 8_{12}$ \quad $(16)$\\\vspace{3mm} }
\end{minipage}
\setlength{\unitlength}{1cm}
\begin{minipage}[t]{3.5cm}
\centering
\epsfxsize=2.5cm
\epsffile{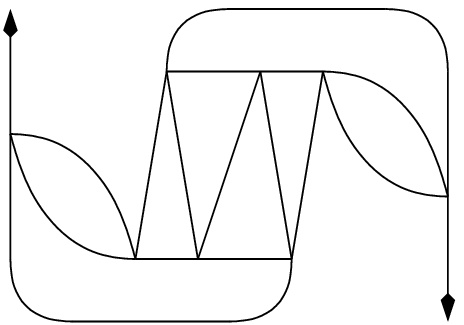}\par
{{\bf Nr.8-2} \quad $C_{2}$\\ $8_{17}$ \quad $(16)$\\\vspace{3mm} }
\end{minipage}
\setlength{\unitlength}{1cm}
\begin{minipage}[t]{3.5cm}
\centering
\epsfxsize=2.5cm
\epsffile{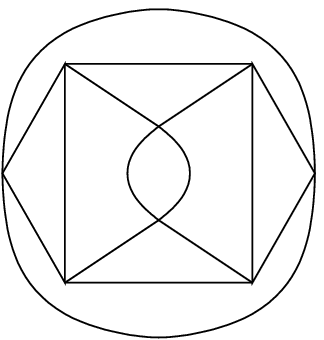}\par
{{\bf Nr.8-3} \quad $D_{2d}$\\ $8^2_{14}$ \quad $(4;12)$\\\vspace{3mm} }
\end{minipage}
\setlength{\unitlength}{1cm}
\begin{minipage}[t]{3.5cm}
\centering
\epsfxsize=2.5cm
\epsffile{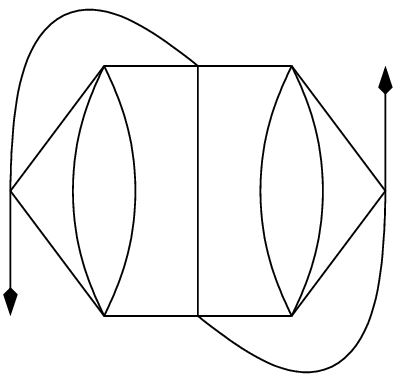}\par
{{\bf Nr.8-4} \quad $C_{2}$\\ $\sim 8^2_{8}$ \quad $(6;10)$\\\vspace{3mm} }
\end{minipage}
\setlength{\unitlength}{1cm}
\begin{minipage}[t]{3.5cm}
\centering
\epsfxsize=2.5cm
\epsffile{6-hedrite8_5.eps}\par
{{\bf Nr.8-5} \quad $D_{2h}$\\ $8^3_{6}$ \quad $(4,6^2)$\\\vspace{3mm} }
\end{minipage}
\setlength{\unitlength}{1cm}
\begin{minipage}[t]{3.5cm}
\centering
\epsfxsize=2.5cm
\epsffile{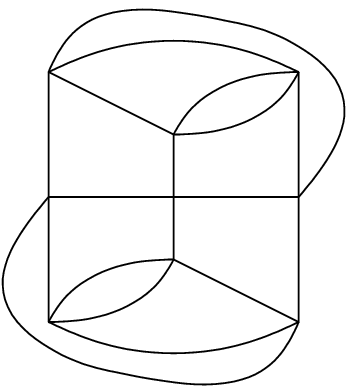}\par
{{\bf Nr.9-1} \quad $C_{2}$\\ $\sim 9_{31}$ \quad $(18)$\\\vspace{3mm} }
\end{minipage}
\setlength{\unitlength}{1cm}
\begin{minipage}[t]{3.5cm}
\centering
\epsfxsize=2.5cm
\epsffile{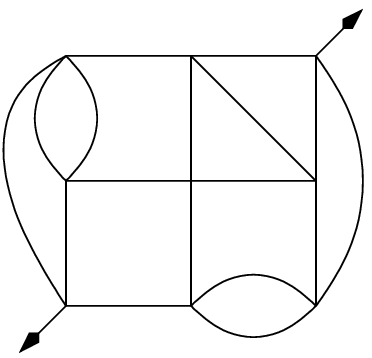}\par
{{\bf Nr.9-2} \quad $C_{1}$\\ $9_{33}$ \quad $(18)$\\\vspace{3mm} }
\end{minipage}
\setlength{\unitlength}{1cm}
\begin{minipage}[t]{3.5cm}
\centering
\epsfxsize=2.5cm
\epsffile{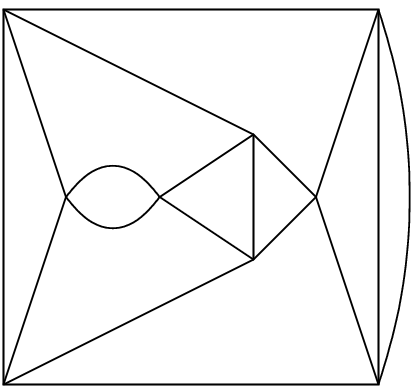}\par
{{\bf Nr.9-3} \quad $C_{s}$\\ $9^2_{38}$ \quad $(4;14)$\\\vspace{3mm} }
\end{minipage}
\setlength{\unitlength}{1cm}
\begin{minipage}[t]{3.5cm}
\centering
\epsfxsize=2.5cm
\epsffile{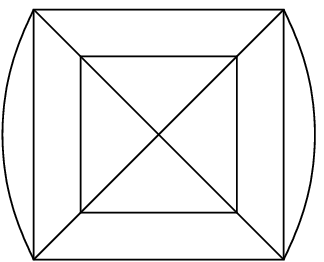}\par
{{\bf Nr.9-4} \quad $C_{2v}$\\ $9^3_{12}$ \quad $(4^2;10)$ red.\\\vspace{3mm} }
\end{minipage}
\setlength{\unitlength}{1cm}
\begin{minipage}[t]{3.5cm}
\centering
\epsfxsize=2.5cm
\epsffile{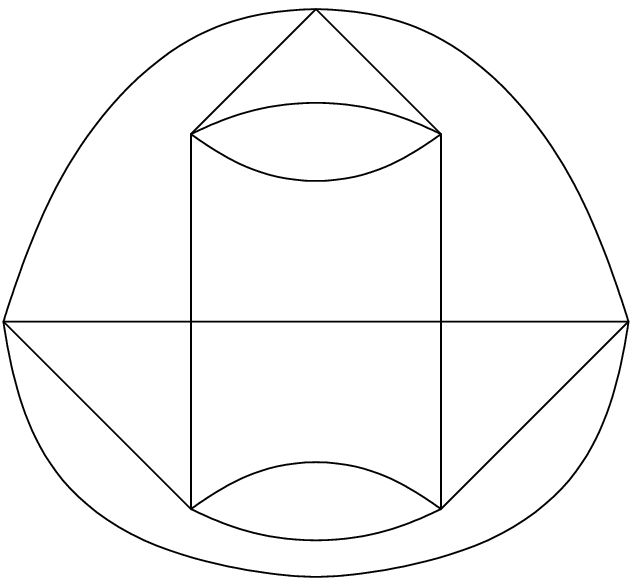}\par
{{\bf Nr.9-5} \quad $C_{s}$\\ $9^3_{11}$ \quad $(4,6;8)$\\\vspace{3mm} }
\end{minipage}
\setlength{\unitlength}{1cm}
\begin{minipage}[t]{3.5cm}
\centering
\epsfxsize=2.5cm
\epsffile{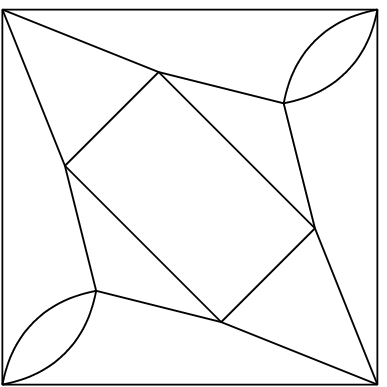}\par
{{\bf Nr.10-1} \quad $C_{2v}$\\ $10_{120}$ \quad $(20)$\\\vspace{3mm} }
\end{minipage}
\setlength{\unitlength}{1cm}
\begin{minipage}[t]{3.5cm}
\centering
\epsfxsize=2.5cm
\epsffile{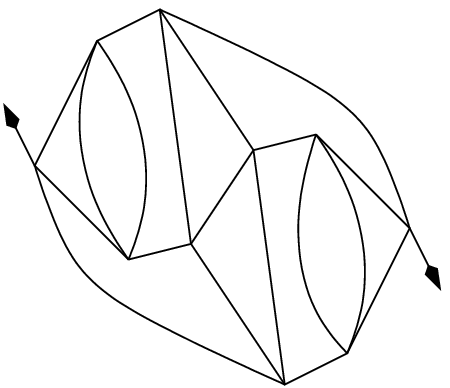}\par
{{\bf Nr.10-2} \quad $C_{2}$\\ $\sim 10_{88}$ \quad $(20)$\\\vspace{3mm} }
\end{minipage}
\setlength{\unitlength}{1cm}
\begin{minipage}[t]{3.5cm}
\centering
\epsfxsize=2.5cm
\epsffile{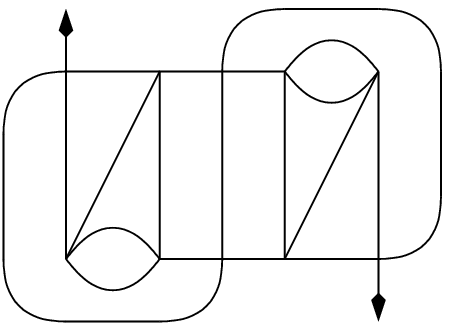}\par
{{\bf Nr.10-3} \quad $C_{2}$\\ $\sim 10_{45}$ \quad $(20)$\\\vspace{3mm} }
\end{minipage}
\setlength{\unitlength}{1cm}
\begin{minipage}[t]{3.5cm}
\centering
\epsfxsize=2.5cm
\epsffile{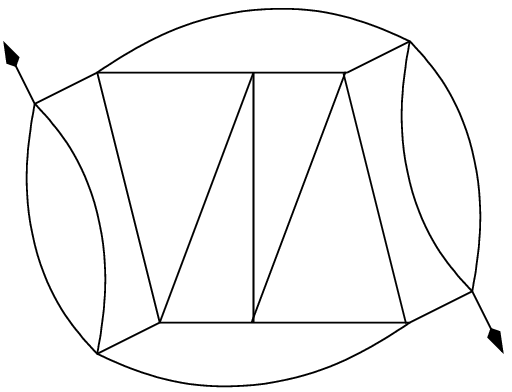}\par
{{\bf Nr.10-4} \quad $C_{2}$\\ $10_{115}$ \quad $(20)$\\\vspace{3mm} }
\end{minipage}
\setlength{\unitlength}{1cm}
\begin{minipage}[t]{3.5cm}
\centering
\epsfxsize=2.5cm
\epsffile{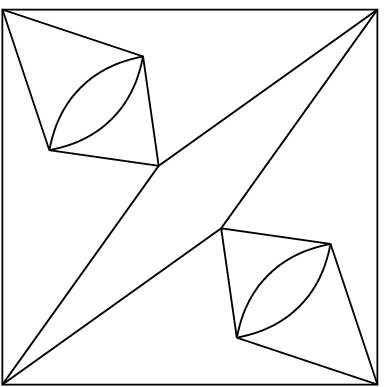}\par
{{\bf Nr.10-5${}^*$} \quad $D_{2h}$\\ $\sim 10^2_{87}$ \quad $(10^2)$\\\vspace{3mm} }
\end{minipage}
\setlength{\unitlength}{1cm}
\begin{minipage}[t]{3.5cm}
\centering
\epsfxsize=2.5cm
\epsffile{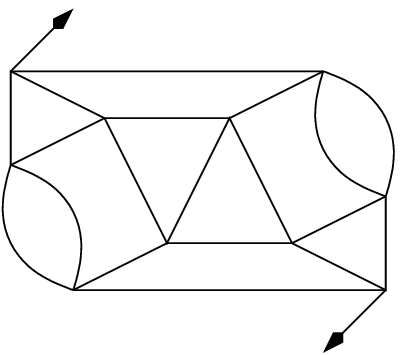}\par
{{\bf Nr.10-6} \quad $C_{2}$\\ $10^2_{86}$ \quad $(6;14)$\\\vspace{3mm} }
\end{minipage}
\setlength{\unitlength}{1cm}
\begin{minipage}[t]{3.5cm}
\centering
\epsfxsize=2.5cm
\epsffile{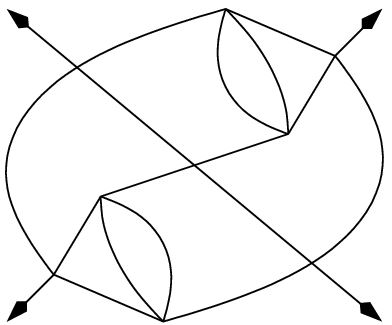}\par
{{\bf Nr.10-7} \quad $C_{2}$\\ $10^2_{43}$ \quad $(4;16)$\\\vspace{3mm} }
\end{minipage}
\setlength{\unitlength}{1cm}
\begin{minipage}[t]{3.5cm}
\centering
\epsfxsize=2.5cm
\epsffile{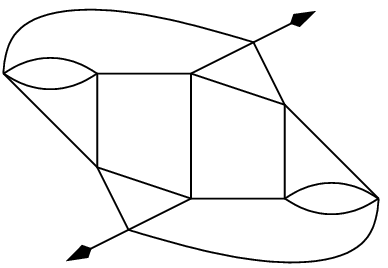}\par
{{\bf Nr.10-8} \quad $C_{2h}$\\ $\sim 10^3_{136}$ \quad $(4;8^2)$\\\vspace{3mm} }
\end{minipage}
\setlength{\unitlength}{1cm}
\begin{minipage}[t]{3.5cm}
\centering
\epsfxsize=2.5cm
\epsffile{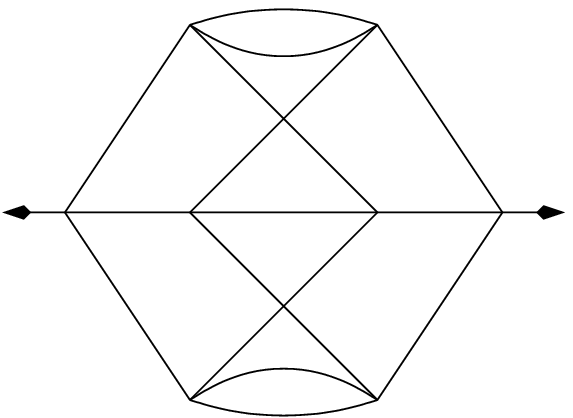}\par
{{\bf Nr.10-9} \quad $C_{2v}$\\ $10^3_{136}$ \quad $(4,6;10)$\\\vspace{3mm} }
\end{minipage}
\setlength{\unitlength}{1cm}
\begin{minipage}[t]{3.5cm}
\centering
\epsfxsize=2.5cm
\epsffile{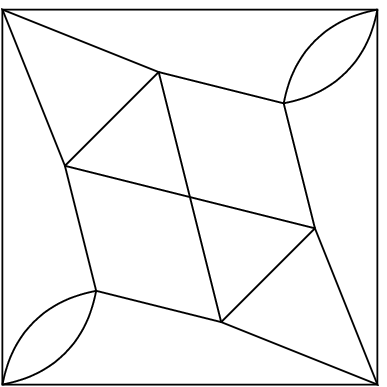}\par
{{\bf Nr.11-1} \quad $C_{2v}$\\ $11_{332}$ \quad $(22)$\\\vspace{3mm} }
\end{minipage}
\setlength{\unitlength}{1cm}
\begin{minipage}[t]{3.5cm}
\centering
\epsfxsize=2.5cm
\epsffile{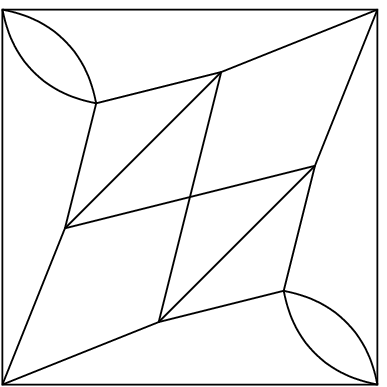}\par
{{\bf Nr.11-2} \quad $C_{2v}$\\ $11_{297}$ \quad $(22)$\\\vspace{3mm} }
\end{minipage}
\setlength{\unitlength}{1cm}
\begin{minipage}[t]{3.5cm}
\centering
\epsfxsize=2.5cm
\epsffile{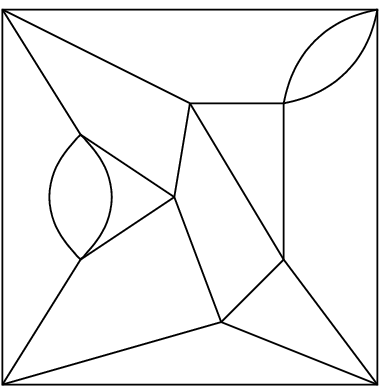}\par
{{\bf Nr.11-3} \quad $C_{1}$\\ $\sim 11_{125}$ \quad $(22)$\\\vspace{3mm} }
\end{minipage}
\setlength{\unitlength}{1cm}
\begin{minipage}[t]{3.5cm}
\centering
\epsfxsize=2.5cm
\epsffile{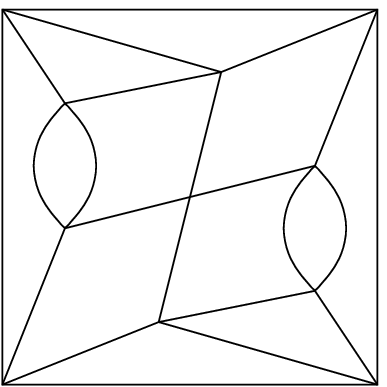}\par
{{\bf Nr.11-4} \quad $C_{2}$\\ $11^2_{317}$ \quad $(8;14)$\\\vspace{3mm} }
\end{minipage}
\setlength{\unitlength}{1cm}
\begin{minipage}[t]{3.5cm}
\centering
\epsfxsize=2.5cm
\epsffile{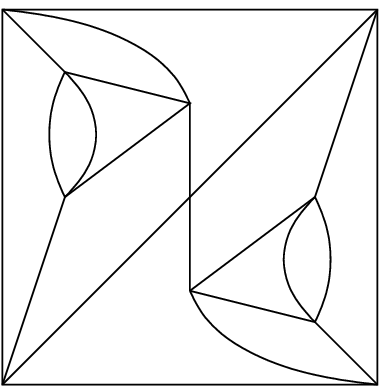}\par
{{\bf Nr.11-5} \quad $C_{2}$\\ $????$ \quad $(8;14)$\\\vspace{3mm} }
\end{minipage}
\setlength{\unitlength}{1cm}
\begin{minipage}[t]{3.5cm}
\centering
\epsfxsize=2.5cm
\epsffile{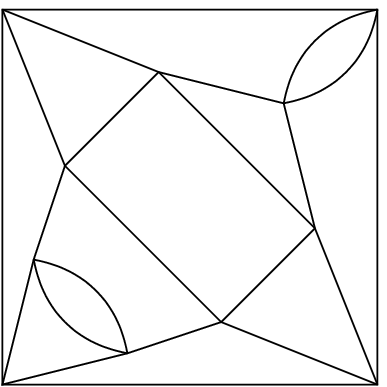}\par
{{\bf Nr.11-6} \quad $C_{s}$\\ $????$ \quad $(8,14)$\\\vspace{3mm} }
\end{minipage}
\setlength{\unitlength}{1cm}
\begin{minipage}[t]{3.5cm}
\centering
\epsfxsize=2.5cm
\epsffile{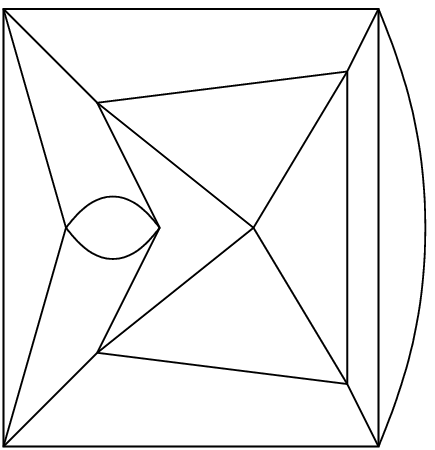}\par
{{\bf Nr.11-7} \quad $C_{s}$\\ $11^2_{351}$ \quad $(10,12)$\\\vspace{3mm} }
\end{minipage}
\setlength{\unitlength}{1cm}
\begin{minipage}[t]{3.5cm}
\centering
\epsfxsize=2.5cm
\epsffile{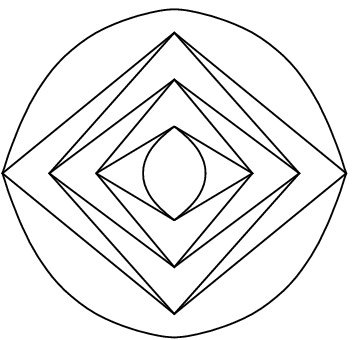}\par
{{\bf Nr.12-1${}^*$} \quad $D_{2d}$\\ $\sim 12_{477}$ \quad $(24)$\\\vspace{3mm} }
\end{minipage}
\setlength{\unitlength}{1cm}
\begin{minipage}[t]{3.5cm}
\centering
\epsfxsize=2.5cm
\epsffile{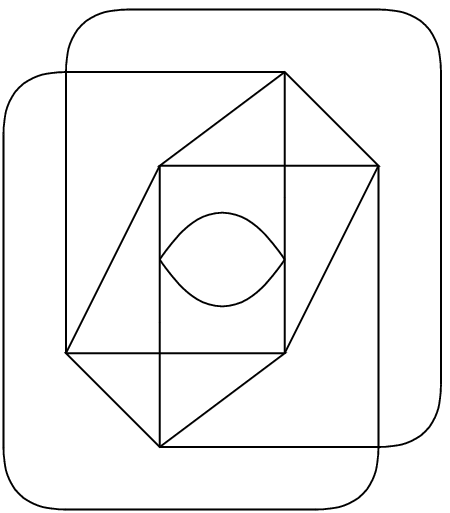}\par
{{\bf Nr.12-2} \quad $D_{2}$\\ $12_{1152}$ \quad $(24)$\\\vspace{3mm} }
\end{minipage}
\setlength{\unitlength}{1cm}
\begin{minipage}[t]{3.5cm}
\centering
\epsfxsize=2.5cm
\epsffile{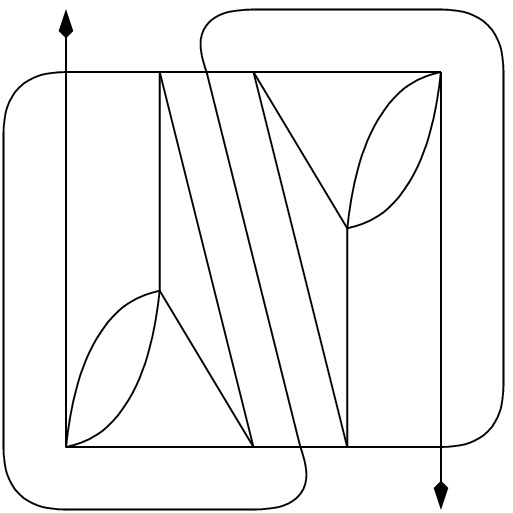}\par
{{\bf Nr.12-3} \quad $C_{2}$\\ $\sim 12_{499}$ \quad $(24)$\\\vspace{3mm} }
\end{minipage}
\setlength{\unitlength}{1cm}
\begin{minipage}[t]{3.5cm}
\centering
\epsfxsize=2.5cm
\epsffile{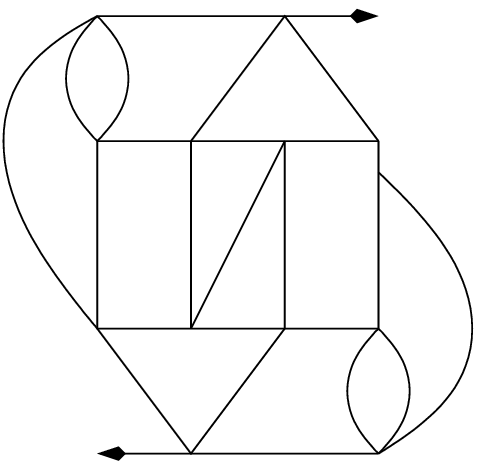}\par
{{\bf Nr.12-4} \quad $C_{2}$\\ $\sim 12_{458}$ \quad $(24)$\\\vspace{3mm} }
\end{minipage}
\setlength{\unitlength}{1cm}
\begin{minipage}[t]{3.5cm}
\centering
\epsfxsize=2.5cm
\epsffile{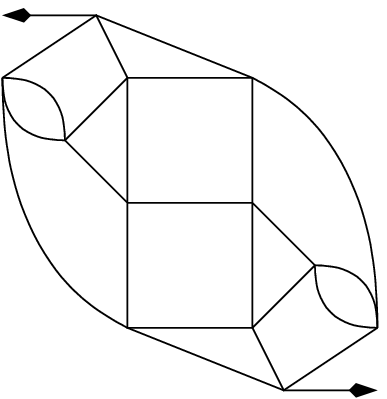}\par
{{\bf Nr.12-5} \quad $C_{2}$\\ $12_{1102}$ \quad $(24)$\\\vspace{3mm} }
\end{minipage}
\setlength{\unitlength}{1cm}
\begin{minipage}[t]{3.5cm}
\centering
\epsfxsize=2.5cm
\epsffile{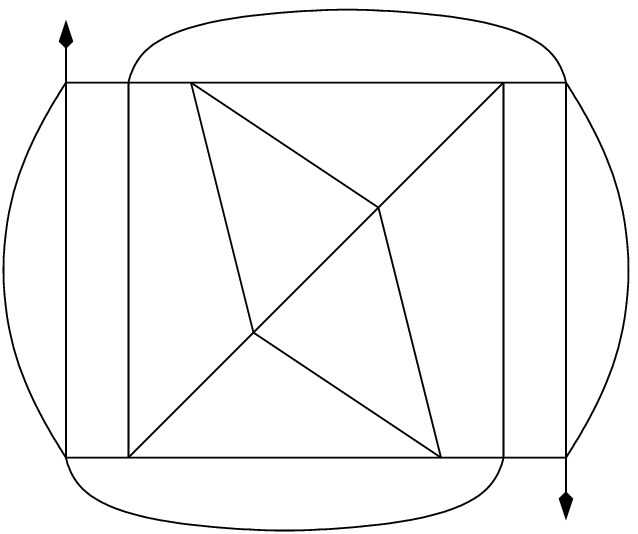}\par
{{\bf Nr.12-6} \quad $C_{2}$\\ $12_{1167}$ \quad $(24)$\\\vspace{3mm} }
\end{minipage}
\setlength{\unitlength}{1cm}
\begin{minipage}[t]{3.5cm}
\centering
\epsfxsize=2.5cm
\epsffile{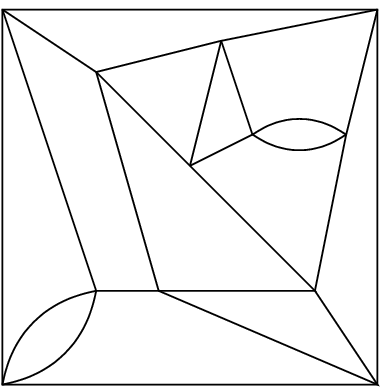}\par
{{\bf Nr.12-7} \quad $C_{1}$\\ $\sim 12_{626}$ \quad $(24)$\\\vspace{3mm} }
\end{minipage}
\setlength{\unitlength}{1cm}
\begin{minipage}[t]{3.5cm}
\centering
\epsfxsize=2.5cm
\epsffile{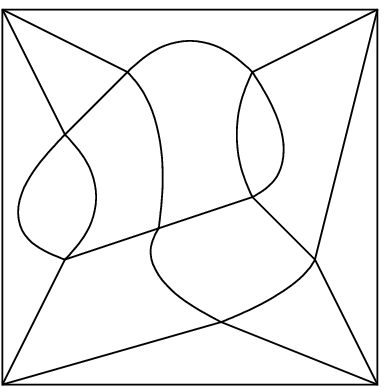}\par
{{\bf Nr.12-8} \quad $C_{1}$\\ $????$ \quad $(6;18)$\\\vspace{3mm} }
\end{minipage}
\setlength{\unitlength}{1cm}
\begin{minipage}[t]{3.5cm}
\centering
\epsfxsize=2.5cm
\epsffile{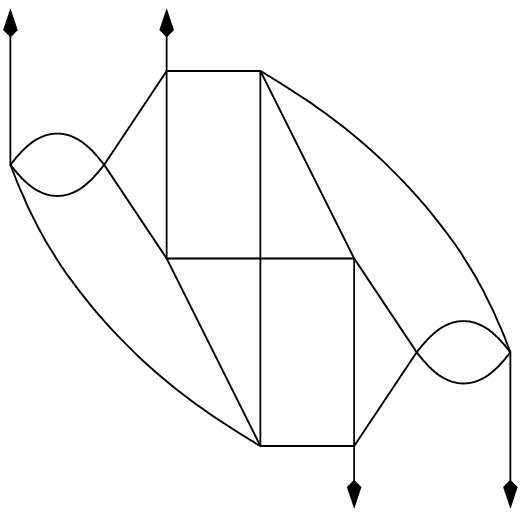}\par
{{\bf Nr.12-9} \quad $C_{2}$\\ $????$ \quad $(10,14)$\\\vspace{3mm} }
\end{minipage}
\setlength{\unitlength}{1cm}
\begin{minipage}[t]{3.5cm}
\centering
\epsfxsize=2.5cm
\epsffile{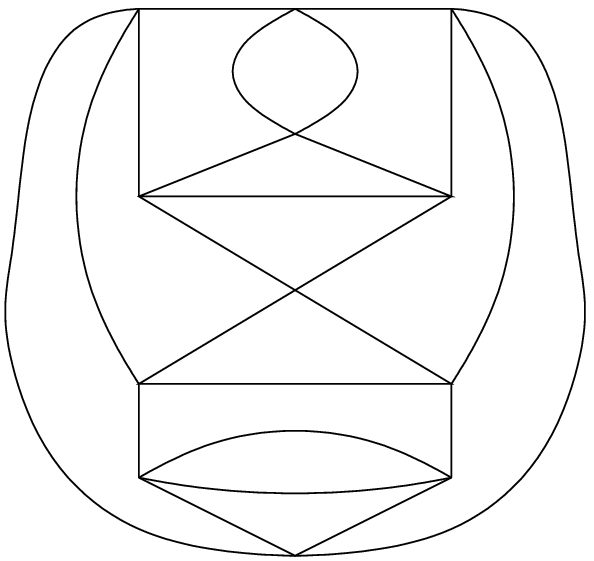}\par
{{\bf Nr.12-10} \quad $C_{s}$\\ $????$ \quad $(8,16)$\\\vspace{3mm} }
\end{minipage}
\setlength{\unitlength}{1cm}
\begin{minipage}[t]{3.5cm}
\centering
\epsfxsize=2.5cm
\epsffile{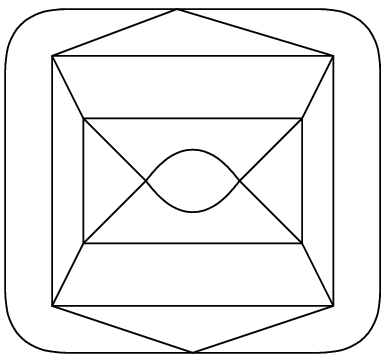}\par
{{\bf Nr.12-11} \quad $D_{2d}$\\ $????$ \quad $(4^2;16)$ red.\\\vspace{3mm} }
\end{minipage}
\setlength{\unitlength}{1cm}
\begin{minipage}[t]{3.5cm}
\centering
\epsfxsize=2.5cm
\epsffile{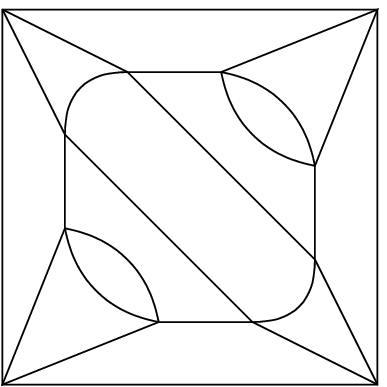}\par
{{\bf Nr.12-12} \quad $C_{2v}$\\ $????$ \quad $(8;8^2)$\\\vspace{3mm} }
\end{minipage}
\setlength{\unitlength}{1cm}
\begin{minipage}[t]{3.5cm}
\centering
\epsfxsize=2.5cm
\epsffile{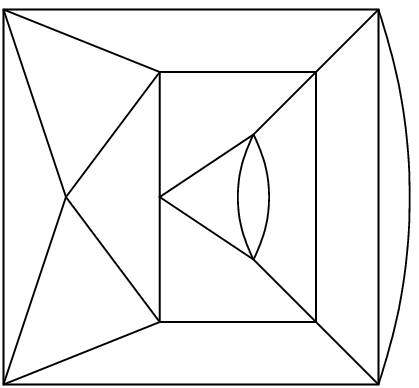}\par
{{\bf Nr.12-13} \quad $C_{s}$\\ $????$ \quad $(6;8,10)$\\\vspace{3mm} }
\end{minipage}
\setlength{\unitlength}{1cm}
\begin{minipage}[t]{3.5cm}
\centering
\epsfxsize=2.5cm
\epsffile{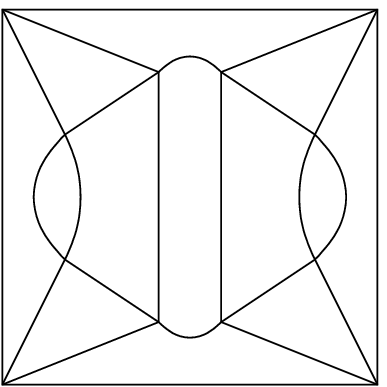}\par
{{\bf Nr.12-14} \quad $D_{2h}$\\ $????$ \quad $(4^2,8^2)$ red.\\\vspace{3mm} }
\end{minipage}
}

\subsection{$7$-hedrites}\label{subsection-7-hedrites}
{\small
\setlength{\unitlength}{1cm}
\begin{minipage}[t]{3.5cm}
\centering
\epsfxsize=2.5cm
\epsffile{7-hedrite7_1sec.eps}\par
{{\bf Nr.7-1} \quad $C_{2v}$\\ $7^2_{6}$ \quad $(4;10)$\\\vspace{3mm} }
\end{minipage}
\setlength{\unitlength}{1cm}
\begin{minipage}[t]{3.5cm}
\centering
\epsfxsize=2.5cm
\epsffile{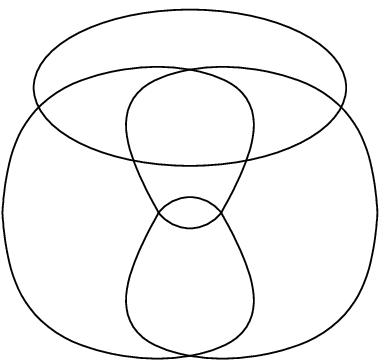}\par
{{\bf Nr.8-1} \quad $C_{s}$\\ $8^2_{13}$ \quad $(4;12)$\\\vspace{3mm} }
\end{minipage}
\setlength{\unitlength}{1cm}
\begin{minipage}[t]{3.5cm}
\centering
\epsfxsize=2.5cm
\epsffile{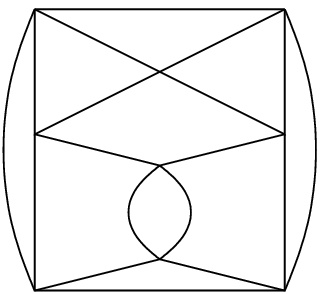}\par
{{\bf Nr.9-1} \quad $C_{s}$\\ $9_{34}$ \quad $(18)$\\\vspace{3mm} }
\end{minipage}
\setlength{\unitlength}{1cm}
\begin{minipage}[t]{3.5cm}
\centering
\epsfxsize=2.5cm
\epsffile{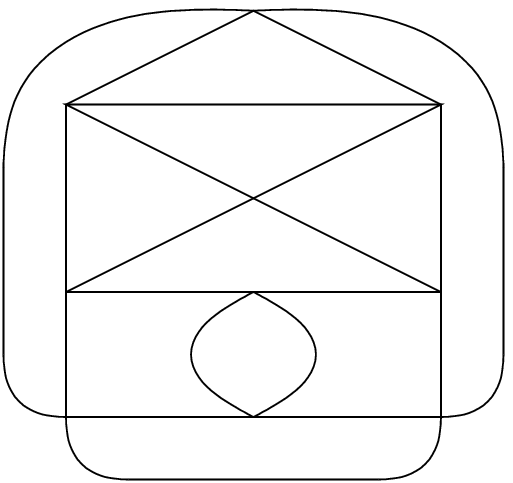}\par
{{\bf Nr.10-1} \quad $C_{s}$\\ $10_{121}$ \quad $(20)$\\\vspace{3mm} }
\end{minipage}
\setlength{\unitlength}{1cm}
\begin{minipage}[t]{3.5cm}
\centering
\epsfxsize=2.5cm
\epsffile{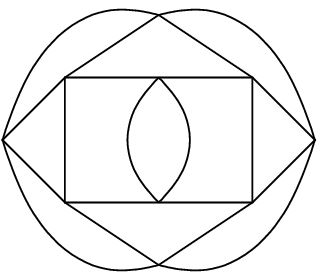}\par
{{\bf Nr.10-2} \quad $C_{2v}$\\ $10^2_{111}$ \quad $(10^2)$\\\vspace{3mm} }
\end{minipage}
\setlength{\unitlength}{1cm}
\begin{minipage}[t]{3.5cm}
\centering
\epsfxsize=2.5cm
\epsffile{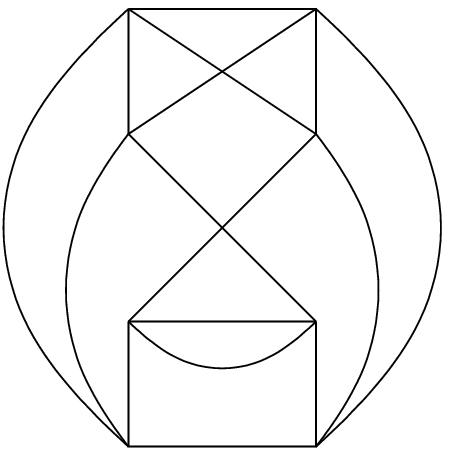}\par
{{\bf Nr.10-3} \quad $C_{s}$\\ $\sim 10^2_{69}$ \quad $(8,12)$\\\vspace{3mm} }
\end{minipage}
\setlength{\unitlength}{1cm}
\begin{minipage}[t]{3.5cm}
\centering
\epsfxsize=2.5cm
\epsffile{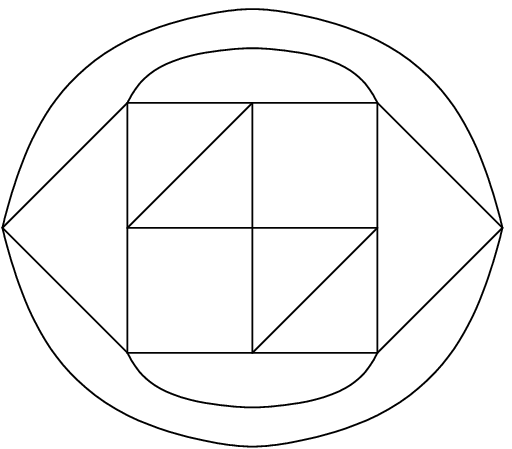}\par
{{\bf Nr.11-1} \quad $C_2$\\ $11_{288}$ \quad $(22)$\\\vspace{3mm} }
\end{minipage}
\setlength{\unitlength}{1cm}
\begin{minipage}[t]{3.5cm}
\centering
\epsfxsize=2.5cm
\epsffile{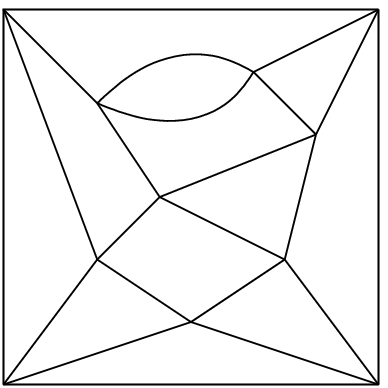}\par
{{\bf Nr.11-2} \quad $C_{1}$\\ $11_{301}$ \quad $(22)$\\\vspace{3mm} }
\end{minipage}
\setlength{\unitlength}{1cm}
\begin{minipage}[t]{3.5cm}
\centering
\epsfxsize=2.5cm
\epsffile{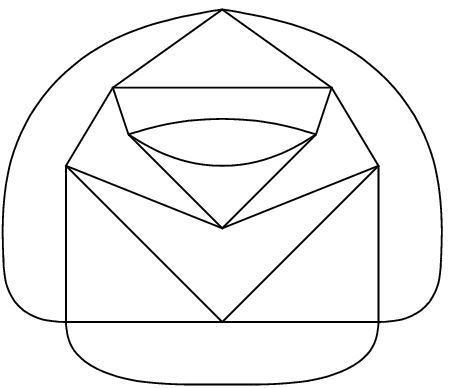}\par
{{\bf Nr.11-3} \quad $C_{s}$\\ $11^2_{150}$ \quad $(8,14)$\\\vspace{3mm} }
\end{minipage}
\setlength{\unitlength}{1cm}
\begin{minipage}[t]{3.5cm}
\centering
\epsfxsize=2.5cm
\epsffile{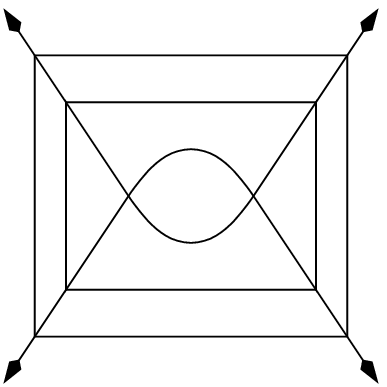}\par
{{\bf Nr.11-4} \quad $C_{2v}$\\ $11^{3}_{487}$ \quad $(4^2;14)$ red.\\\vspace{3mm} }
\end{minipage}
\setlength{\unitlength}{1cm}
\begin{minipage}[t]{3.5cm}
\centering
\epsfxsize=2.5cm
\epsffile{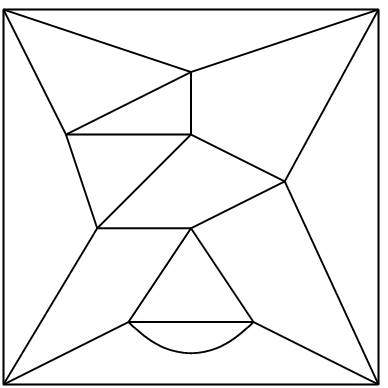}\par
{{\bf Nr.12-1} \quad $C_{1}$\\ $\sim 12_{361}$ \quad $(24)$\\\vspace{3mm} }
\end{minipage}
\setlength{\unitlength}{1cm}
\begin{minipage}[t]{3.5cm}
\centering
\epsfxsize=2.5cm
\epsffile{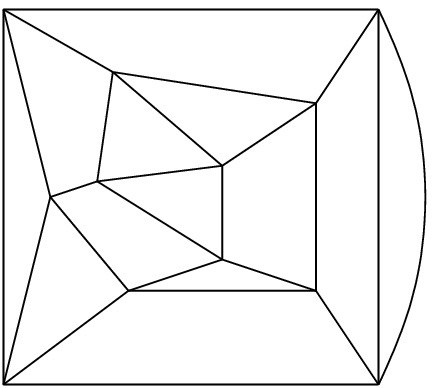}\par
{{\bf Nr.12-2} \quad $C_{1}$\\ $????$ \quad $(6;18)$\\\vspace{3mm} }
\end{minipage}
\setlength{\unitlength}{1cm}
\begin{minipage}[t]{3.5cm}
\centering
\epsfxsize=2.5cm
\epsffile{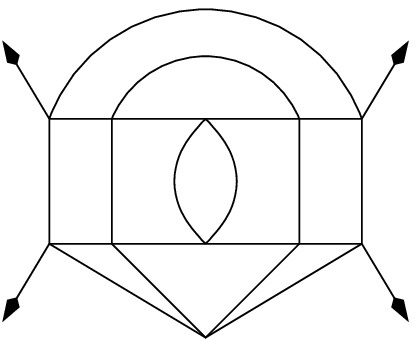}\par
{{\bf Nr.12-3} \quad $C_{s}$\\ $????$ \quad $(10,14)$\\\vspace{3mm} }
\end{minipage}
\setlength{\unitlength}{1cm}
\begin{minipage}[t]{3.5cm}
\centering
\epsfxsize=2.5cm
\epsffile{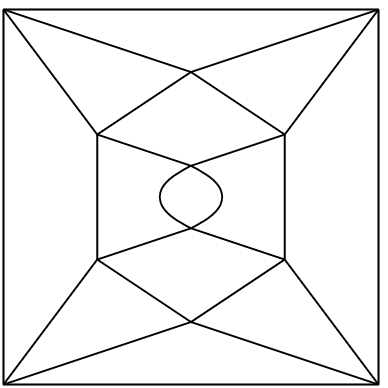}\par
{{\bf Nr.12-4} \quad $C_{2v}$\\ $????$ \quad $(6^2;12)$\\\vspace{3mm} }
\end{minipage}
\setlength{\unitlength}{1cm}
\begin{minipage}[t]{3.5cm}
\centering
\epsfxsize=2.5cm
\epsffile{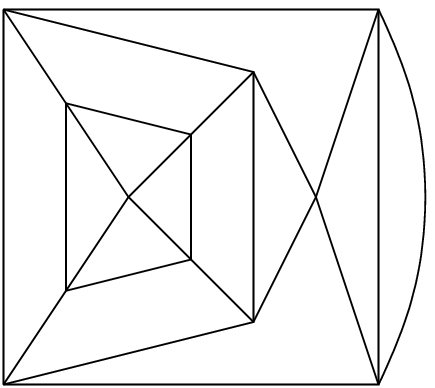}\par
{{\bf Nr.12-5} \quad $C_{s}$\\ $????$ \quad $(4^2;16)$ red.\\\vspace{3mm} }
\end{minipage}
}

\subsection{$8$-hedrites}\label{subsection-8-hedrites}
{\small
\setlength{\unitlength}{1cm}
\begin{minipage}[t]{3.5cm}
\centering
\epsfxsize=2.5cm
\epsffile{8-hedrite6-1.eps}\par
{{\bf Nr.6-1} \quad $O_h$\\ $6^3_2$ \quad $(4^3)$\\\vspace{3mm} }
\end{minipage}
\setlength{\unitlength}{1cm}
\begin{minipage}[t]{3.5cm}
\centering
\epsfxsize=2.5cm
\epsffile{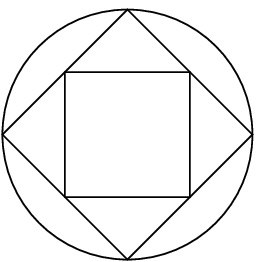}\par
{{\bf Nr.8-1} \quad $D_{4d}$\\ $8_{18}$ \quad $(16)$\\\vspace{3mm} }
\end{minipage}
\setlength{\unitlength}{1cm}
\begin{minipage}[t]{3.5cm}
\centering
\epsfxsize=2.5cm
\epsffile{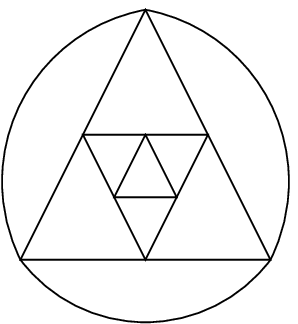}\par
{{\bf Nr.9-1} \quad $D_{3h}$\\ $9_{40}$ \quad $(18)$\\\vspace{3mm} }
\end{minipage}
\setlength{\unitlength}{1cm}
\begin{minipage}[t]{3.5cm}
\centering
\epsfxsize=2.5cm
\epsffile{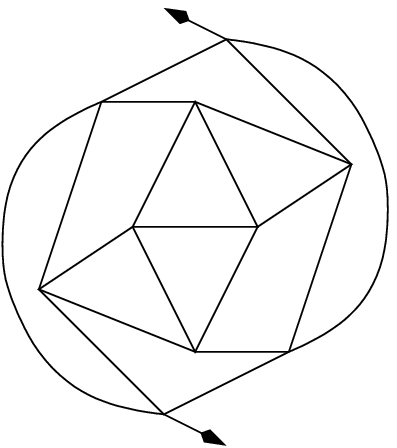}\par
{{\bf Nr.10-1} \quad $D_{2}$\\ $10^2_{56}$ \quad $(6;14)$\\\vspace{3mm} }
\end{minipage}
\setlength{\unitlength}{1cm}
\begin{minipage}[t]{3.5cm}
\centering
\epsfxsize=2.5cm
\epsffile{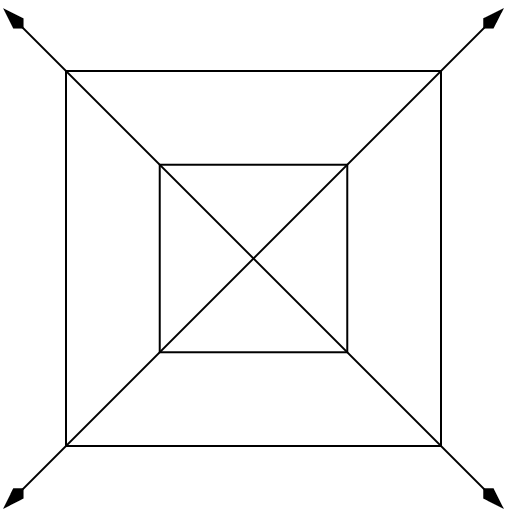}\par
{{\bf Nr.10-2} \quad $D_{4h}$\\ $10^4_{169}$ \quad $(4^2,6^2)$ red.\\\vspace{3mm} }
\end{minipage}
\setlength{\unitlength}{1cm}
\begin{minipage}[t]{3.5cm}
\centering
\epsfxsize=2.5cm
\epsffile{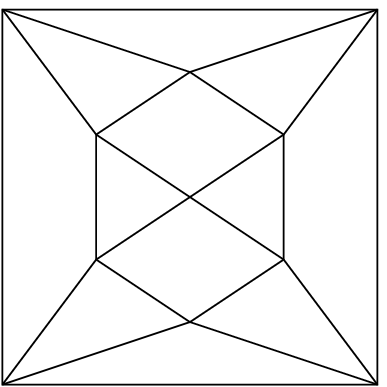}\par
{{\bf Nr.11-1} \quad $C_{2v}$\\ $11^3_{520}$ \quad $(6^2;10)$\\\vspace{3mm} }
\end{minipage}
\setlength{\unitlength}{1cm}
\begin{minipage}[t]{3.5cm}
\centering
\epsfxsize=2.2cm
\epsffile{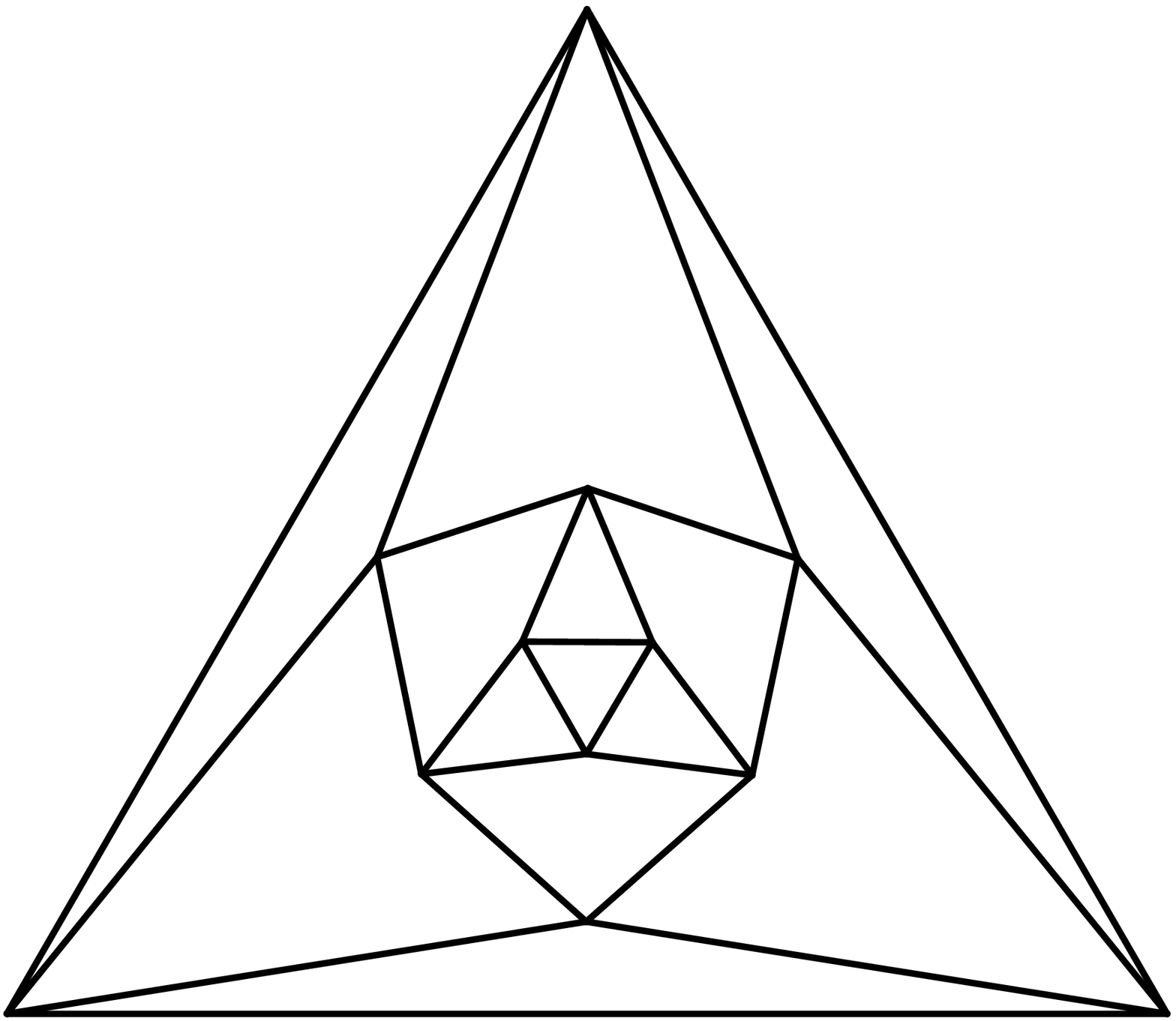}\par
{{\bf Nr.12-1} \quad $D_{3d}$\\ $12_{1019}$ \quad $(24)$\\\vspace{3mm} }
\end{minipage}
\setlength{\unitlength}{1cm}
\begin{minipage}[t]{3.5cm}
\centering
\epsfxsize=2.5cm
\epsffile{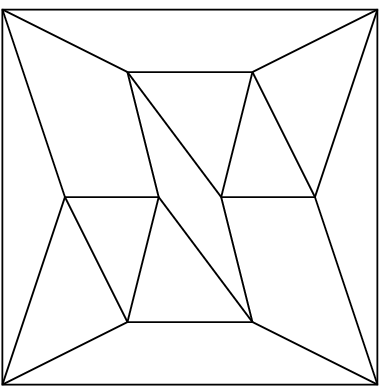}\par
{{\bf Nr.12-2} \quad $D_2$\\ $12_{868}$ \quad $(24)$\\\vspace{3mm} }
\end{minipage}
\setlength{\unitlength}{1cm}
\begin{minipage}[t]{3.5cm}
\centering
\epsfxsize=2.5cm
\epsffile{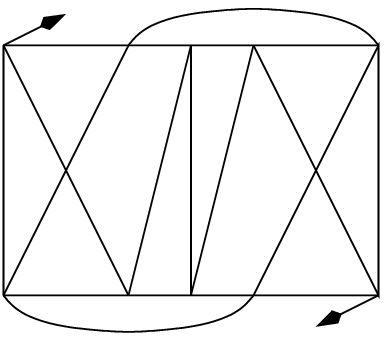}\par
{{\bf Nr.12-3} \quad $C_2$\\ $????$ \quad $(6;18)$\\\vspace{3mm} }
\end{minipage}
\setlength{\unitlength}{1cm}
\begin{minipage}[t]{3.5cm}
\centering
\epsfxsize=2.5cm
\epsffile{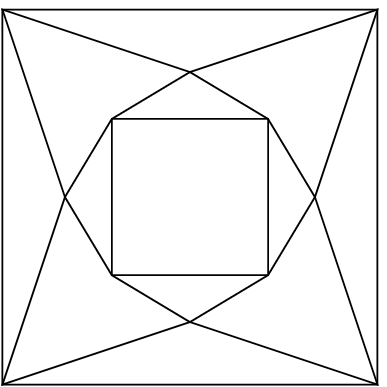}\par
{{\bf Nr.12-4} \quad $O_h$\\ $????$ \quad $(6^4)$\\\vspace{3mm} }
\end{minipage}
\setlength{\unitlength}{1cm}
\begin{minipage}[t]{3.5cm}
\centering
\epsfxsize=2.2cm
\epsffile{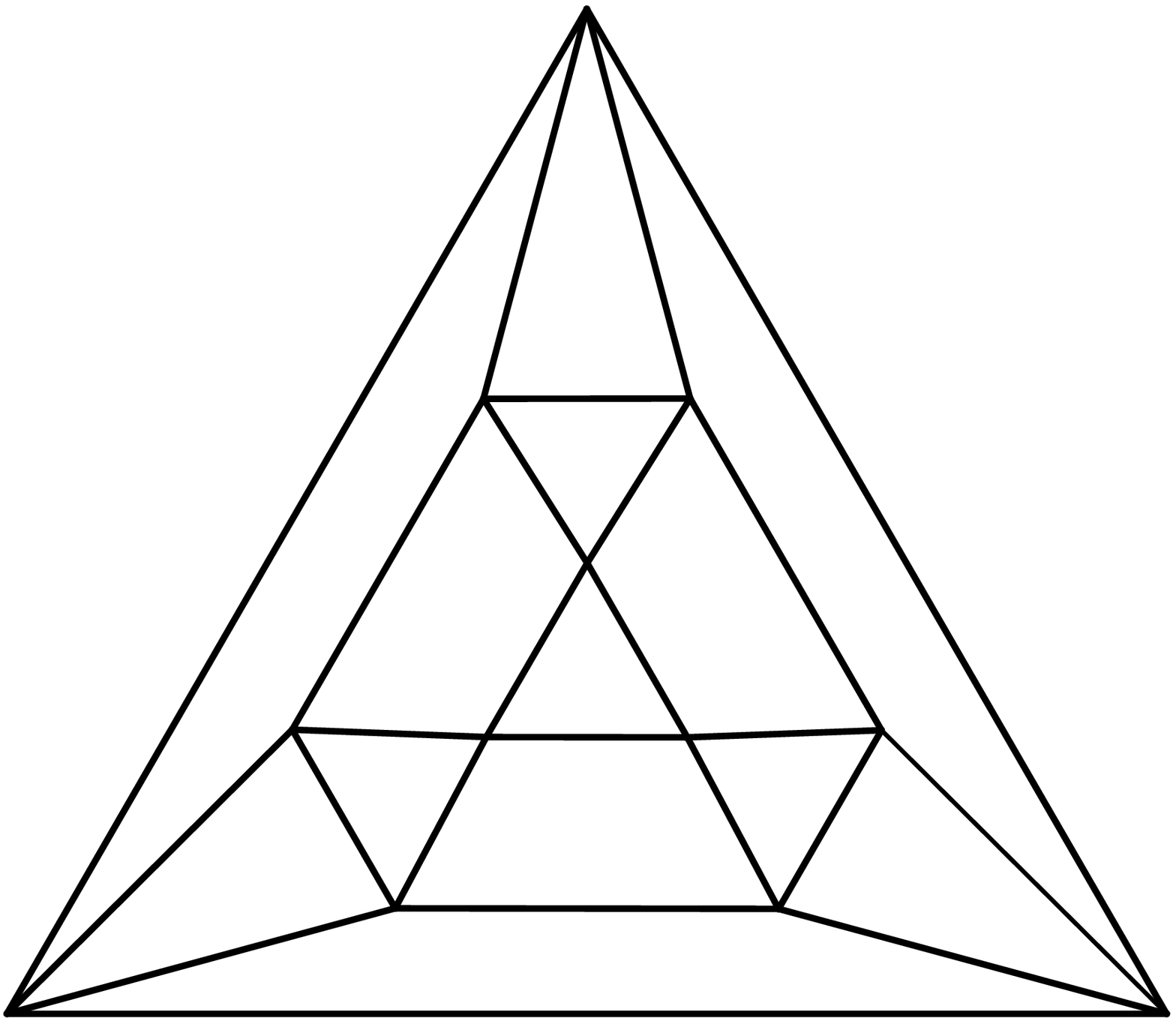}\par
{{\bf Nr.12-5} \quad $D_{3h}$\\ $????$ \quad $(6^4)$\\\vspace{3mm} }
\end{minipage}
}

\begin{table}
\begin{center}
{\scriptsize
\begin{minipage}{7cm}
\begin{tabular}{||l|l|l|l||}
\hline\hline
Nr.     &Group    &CC-vector      &alt.knot\\\hline\hline
\multicolumn{4}{||c||}{$4$-hedrites}\\\hline
{\bf 14-1${}^*$}&$D_{2d}$       &$14^2$         &\\
{\bf 14-2}      &$D_2$  &$14^2$         &\\
{\bf 14-3${}^*$}&$D_{2h}$       &$2^7, 14$ red.  &\\\hline\hline
\multicolumn{4}{||c||}{$5$-hedrites}\\\hline
{\bf 13-1${}^*$}&$C_{2v}$       &$26$           &$13_{3097}$\\
{\bf 13-2}      &$C_1$  &$26$           &$13_{4054}$\\
{\bf 13-3}      &$C_2$  &$6^2; 14$      &\\
{\bf 13-4}      &$C_s$  &$6; 8, 12$     &\\\hline
{\bf 14-1}      &$C_s$  &$28$           &$14_{16368}$\\
{\bf 14-2}      &$C_1$  &$6; 22$                &\\
{\bf 14-3}      &$C_1$  &$8; 20$                &\\
{\bf 14-4}      &$C_{2v}$       &$8^2; 12$      &\\
{\bf 14-5}      &$C_{2v}$       &$4^3; 16$ red. &\\
{\bf 14-6}      &$C_{2v}$       &$6^2; 8^2$ red.        &\\
{\bf 14-7}      &$C_{2v}$       &$4^3, 8^2$ red.        &\\\hline
{\bf 15-1}      &$D_3$  &$30$   &$15_{83814}$\\
{\bf 15-2}      &$C_{2v}$       &$30$   &$15_{54593}$\\
{\bf 15-3}      &$C_s$  &$30$   &$15_{83824}$\\
{\bf 15-4}      &$C_2$  &$30$   &$15_{20161}$\\
{\bf 15-5}      &$C_s$  &$12,18$        &\\
{\bf 15-6}      &$C_s$  &$14,16$        &\\
{\bf 15-7}      &$C_1$  &$8;22$ &\\
{\bf 15-8}      &$C_{2v}$       &       $8^2;14$&\\
{\bf 15-9}      &$C_s$  &$4^3; 18$ red. &\\
{\bf 15-10}     &$C_s$  &$4^3, 8; 10$ red.      &\\\hline\hline
\multicolumn{4}{||c||}{$6$-hedrites}\\\hline
{\bf 13-1}      &$C_2$  &$26$           &$\sim 13_{1739}$\\
{\bf 13-2}      &$C_2$  &$26$           &$13_{3586}$\\
{\bf 13-3}      &$C_2$  &$26$           &$\sim 13_{1345}$\\
{\bf 13-4}      &$C_s$  &$26$           &$13_{3811}$\\
{\bf 13-5}      &$C_1$  &$26$           &$13_{1485}$\\
{\bf 13-6}      &$C_1$  &$26$           &$13_{3957}$\\
{\bf 13-7}      &$C_1$  &$26$           &$\sim 13_{2957}$\\
{\bf 13-8}      &$C_1$  &$8; 18$                &\\
{\bf 13-9}      &$C_s$  &$12, 14$               &\\
{\bf 13-10}     &$C_s$  &$4^2; 18$ red. &\\
{\bf 13-11}     &$C_{2v}$       &$8^2, 10$ red. &\\
{\bf 13-12}     &$C_{2v}$       &$8^2; 10$      &\\
{\bf 13-13}     &$C_{2v}$       &$4^3; 14$ red. &\\
{\bf 13-14}     &$C_s$  &$4^2,8;10$ red.        &\\\hline
{\bf 14-1}      &$C_{2h}$       &$28$           &$14_{17173}$\\
{\bf 14-2}      &$C_{2}$        &$28$           &$14_{17079}$\\
{\bf 14-3}      &$C_2$  &$28$           &$14_{8767}$\\
{\bf 14-4}      &$C_2$  &$28$           &$14_{17734}$\\
{\bf 14-5}      &$C_2$  &$28$           &$14_{17148}$\\
{\bf 14-6}      &$C_1$  &$28$           &$14_{17309}$\\
{\bf 14-7}      &$C_1$  &$28$           &$14_{5570}$\\
{\bf 14-8}      &$C_{2}$        &$6; 22$                &\\
{\bf 14-9}      &$C_2$  &$6; 22$                &\\
{\bf 14-10}     &$C_2$  &$10; 18$               &\\
{\bf 14-11}     &$C_2$  &$10; 18$               &\\
{\bf 14-12}     &$C_s$  &$10, 18$               &\\
{\bf 14-13${}^{*}$}     &$D_{2h}$       &$14^2$         &\\
{\bf 14-14}     &$C_{2v}$       &$14^2$         &\\
{\bf 14-15}     &$C_2$  &$12, 16$               &\\
{\bf 14-16}     &$C_{s}$        &$12, 16$               &\\
{\bf 14-17}     &$C_2$  &$4^2; 20$ red. &\\
{\bf 14-18}     &$C_2$  &$8; 10^2$      &\\
{\bf 14-19}     &$C_2$  &$6^2; 16$      &\\
{\bf 14-20}     &$D_{2h}$       &$6^2, 8^2$     &\\
{\bf 14-21}     &$C_{2v}$       &$6^3, 10$ red. &\\
\hline\hline
\end{tabular}
\end{minipage}
\begin{minipage}[t]{7cm}
\begin{tabular}{||l|l|l|l||}
\hline\hline
{\bf 14-22}     &$C_{2h}$       &$4^2; 10^2$ red.       &\\
{\bf 14-23}     &$C_{2v}$       &$4^2, 8; 12$ red.      &\\\hline
{\bf 15-1}      &$C_2$  &$30$   &$\sim 15_{39533}$\\
{\bf 15-2}      &$C_2$  &$30$   &$15_{66949}$\\
{\bf 15-3}      &$C_2$  &$30$   &$15_{83008}$\\
{\bf 15-4}      &$C_1$  &$30$   &$15_{45248}$\\
{\bf 15-5}      &$C_1$  &$30$   &$\sim 15_{20975}$\\
{\bf 15-6}      &$C_1$  &$30$   &$15_{64488}$\\
{\bf 15-7}      &$C_1$  &$30$   &$\sim 15_{45357}$\\
{\bf 15-8}      &$C_1$  &$6; 24$        &\\
{\bf 15-9}      &$C_1$  &$6; 24$        &\\
{\bf 15-10}     &$C_{2v}$       &$8; 22$     &\\
{\bf 15-11}     &$C_1$  &$8; 22$        &\\
{\bf 15-12}     &$C_1$  &$8; 22$        &\\
{\bf 15-13}     &$C_s$  &$10, 20$       &\\
{\bf 15-14}     &$C_1$  &$10,20$        &\\
{\bf 15-15}     &$C_{2v}$       &$8;10,12$      &\\
{\bf 15-16}     &$C_{2v}$       &$6^2,8;10$     &\\
{\bf 15-17}     &$C_s$  &$6^3; 12$ red. &\\\hline\hline
\multicolumn{4}{||c||}{$7$-hedrites}\\\hline
{\bf 13-1}      &$C_s$  &$26$           &$13_{3861}$\\
{\bf 13-2}      &$C_1$  &$26$           &$13_{3769}$\\
{\bf 13-3}      &$C_1$  &$6; 20$                &\\
{\bf 13-4}      &$C_1$  &$10, 16$               &\\
{\bf 13-5}      &$C_s$  &$10, 16$               &\\
{\bf 13-6}      &$C_{2v}$       &$6^2; 14$      &\\
{\bf 13-7}      &$C_{s}$        &$6^2; 14$      &\\\hline
{\bf 14-1}      &$C_1$  &$28$           &$14_{13725}$\\
{\bf 14-2}      &$C_1$  &$28$           &$14_{10841}$\\
{\bf 14-3}      &$C_1$  &$28$           &$14_{5714}$\\
{\bf 14-4}      &$C_1$  &$28$           &$14_{14207}$\\
{\bf 14-5}      &$C_1$  &$6; 22$                &\\
{\bf 14-6}      &$C_s$  &$10, 18$               &\\
{\bf 14-7}      &$C_2$  &$14^2$         &\\
{\bf 14-8}      &$C_s$  &$6^2; 16$      &\\
{\bf 14-9}      &$C_{2v}$       &$6^2; 16$      &\\\hline
{\bf 15-1}      &$C_2$  &$30$   &$15_{82225}$\\
{\bf 15-2}      &$C_1$  &$30$   &$15_{60207}$\\
{\bf 15-3}      &$C_1$  &$30$   &$15_{80242}$\\
{\bf 15-4}      &$C_s$  &$6;24$ &\\
{\bf 15-5}      &$C_1$  &$6;24$ &\\
{\bf 15-6}      &$C_1$  &$6;24$ &\\
{\bf 15-7}      &$C_1$  &$10,20$        &\\
{\bf 15-8}      &$C_1$  &$10,20$        &\\
{\bf 15-9}      &$C_{2v}$  &$14,16$ &\\
{\bf 15-10}     &$C_1$  &$14,16$        &\\
{\bf 15-11}     &$C_{2v}$       &$6^2; 18$      &\\
{\bf 15-12}     &$C_{2v}$       &$4^3; 18$ red. &\\\hline
\hline
\multicolumn{4}{||c||}{$8$-hedrites}\\\hline
{\bf 13-1}      &$C_2$  &$26$           &$13_{3478}$\\
{\bf 13-2}      &$C_{2v}$       &$6^2; 14$      &\\\hline
{\bf 14-1}      &$C_2$  &$28$           &$14_{17895}$\\
{\bf 14-2}      &$C_s$  &$6; 22$                &\\
{\bf 14-3}      &$D_2$  &$6; 22$                &\\
{\bf 14-4}      &$D_{2d}$       &$14^2$         &\\
{\bf 14-5}      &$C_2$  &$6^2; 16$      &\\
{\bf 14-6}      &$D_2$  &$8; 10^2$      &\\
{\bf 14-7}      &$D_{4h}$       &$6^2, 8^2$     &\\
{\bf 14-8}      &$D_{4h}$       &$4^3, 8^2$ red.        &\\\hline
{\bf 15-1}      &$C_2$  &$30$   &$15_{82477}$\\
{\bf 15-2}      &$C_s$  &$6; 24$        &\\
{\bf 15-3}      &$C_s$  &$6;24$ &\\
{\bf 15-4}      &$C_2$  &$8; 22$        &\\
{\bf 15-5}      &$D_{3h}$       &$10^3$ &\\\hline\hline
\end{tabular}
\end{minipage}
}
\end{center}
\caption{All $i$-hedrites with $13$, $14$ and $15$ vertices}
\label{tab:i-hedrite13_14}
\end{table}


\newpage










\begin{thebibliography}{99}



\bibitem[Cox71]{Cox71}
H.S.M.Coxeter, {\em Virus macromolecules and geodesic domes}, in {\em A spectrum of mathematics}; ed. by J.C.Butcher, Oxford University Press/Auckland University Press: Oxford, U.K./Auckland New-Zealand, (1971) 98--107.


\bibitem[DDF02]{DDF}
M.Deza, M.Dutour and P.W.Fowler,
{\em Zigzags, Rail-roads and Knots in Fullerenes},
submitted.


\bibitem[DeGr99]{DG2}
M.Deza and V.P.Grishukhin,
{\em $l_1$-embeddable polyhedra},
in: Algebras and Combinatorics, Int. Congress CAC '97 Hong Kong,
ed. by K.P. Shum et al., Springer-Verlag (1999) 189--210.


\bibitem[DHL02]{DHL}
M.Deza, T.Huang and K-W.Lih,
{\em Central Circuit Coverings of Octahedrites and Medial Polyhedra},
Journal of Math. Research \& Exposition {\bf 22-1} (2002) 49--66.


\bibitem[DeSt02]{DSt}
M.Deza and M.Shtogrin,
{\em Octahedrites}, 
Symmetry, Special Issue ``Polyhedra and Science and Art'', 2002.


\bibitem[Dut]{Dut}
M.Dutour, {\em PlanGraph, a gap package for Planar Graph}, in preparation.


\bibitem[DD02]{DD02}
M.Dutour, M.Deza, {\em Zigzag Structure of Simple Bifaced Polyhedra}, submitted.


\bibitem[DD03]{DD03}
M.Dutour, M.Deza, {\em Goldberg-Coxeter Construction for convex polyhedra}, in preparation.


\bibitem[GaKe94]{GK}
M.L.Gargano and J.W.Kennedy,
{\em Gaussian graphs and digraphs}, Congressus Numerantium {\bf 101}
(1994) 161--170.


\bibitem[GoRo01]{God}
C.Godsil and G.Royle, {\em Algebraic Graph Theory}, Graduate Texts in 
Mathematics {\bf 207}, Springer-Verlag, Berlin - New York, 2001.


\bibitem[Gold37]{Gold37}
M.Goldberg, {\em A class of multisymmetric polyhedra}, Tohoku Math.
Journal, {\bf 43} (1937) 104--108.


\bibitem[Gr\"{u}n67]{Gr}
B.Gr\"{u}nbaum, {\em Convex polytopes}, Interscience, New York, 1967.


\bibitem[Gr\"{u}n72]{Gr2}
B.Gr\"{u}nbaum, {\em Arrangements and Spreads}, Regional Conference Series in
Mathematics {\bf 10}, American Mathematical Society, 1972.


\bibitem[Gr\"{u}nMo63]{GrMo}
B.Gr\"{u}nbaum and T.S.Motzkin, {\em The number of hexagons and the simplicity
of geodesics on certain polyhedra}, Canadian Journal of Mathematics {\bf 15} (1963) 744--751.


\bibitem[Harb97]{Ha}
H.Harborth, {\em Eulerian straight ahead cycles in drawings of complete
bipartite graphs}, Bericht 97/23, Institute f\"{u}r Mathematik, Tech. 
Universit\"{a}t
Braunschweg, 1997.


\bibitem[Heid98]{He}
O.Heidemeier, {\em Die Erzeugung von 4-regul\"{a}ren, planaren,
simplen, zusammenh\"{a}ngenden Graphen mit vorgegebenen Fl\"{a}chentypen},
Diplomarbeit, Universit\"{a}t Bielefeld, Fakult\"{a}t f\"{u}r Wirtschaft und
Mathematik, 1998. 


\bibitem[Jeo95]{Je}
D.Jeong, {\em Realizations with a cut-through Eulerian circuit},
Discrete Mathematics {\bf 137} (1995) 265--275.




\bibitem[Kaw96]{Kaw}
A.Kawauchi, {\em A survey of knot theory}, Birkh\"{a}user, 1996.


\bibitem[Kir85]{Kirk}
T. Kirkman, {\em The enumeration, description, and construction of knots with fewer than $10$ crossings}, Trans. Roy. Soc. Edin. {\bf 32} (1885), 281--309.


\bibitem[Kot69]{Ko}
A.Kotzig, {\em Eulerian lines in finite 4-valent graphs and their 
transformations}, in: Theory of Graphs, Proceedings of a colloquium, 
Tihany 1966, ed. by P.Erdos and G.Katona, Academic Press, 
New York (1969) 219--230.


\bibitem[Liu98]{Liu}
Liu Yanpei, {\em Embedding in Graphs}, Kluwer, Dodrecht, 1998.


\bibitem[PTZ96]{PTZ}
T.Pisanski, T.Tucker and A.Zitnik, {\em Eulerian Embedding of Graphs},
University of Ljubljana, IMMF Preprint Series {\bf 34}
(1996) 531.


\bibitem[Rol76]{Rolf}
D.Rolfsen, {\em Knots and Links}, Mathematics Lecture Series 7, Publish or
Perish, Berkeley, 1976;
second corrected printing: Publish or Perish, Houston, 1990.


\bibitem[Sha75]{Sh}
H.Shank, {\em The Theory of Left-Right Paths}, in: Combinatorial 
Mathematics III,
Proceedings of 3rd Australian Conference, St Lucia 1974, Lecture Notes in
Mathematics {\bf 452}, Springer-Verlag, Berlin - New York (1975),  pp. 42--54.


\bibitem[Thi]{T}
M.Thistlewaite, {\em Homepage}, \url{http://www.math.utk.edu/~morwen}.


\end{thebibliography}
\end{document}